\documentclass[preprint,11pt]{elsarticle}

\usepackage{latexsym,amssymb,amsmath}
\usepackage{bm}
\usepackage{graphicx}
\usepackage[normal]{subfigure}

\newcommand{\R}{\mathbb{R}}
\newcommand{\BB}{\mathbf{B}}
\newcommand{\bb}{\mathbf{b}}
\newcommand{\vv}{\mathbf{v}}
\newcommand{\OO}{\mathbf{0}}
\newcommand{\Bcal}{\mathcal{B}}

\newdefinition{remark}{Remark}

\begin{document}

\begin{frontmatter}

\title{On a class of two-dimensional incomplete Riemann solvers}

\author{Jos\'e M. Gallardo\corref{cor}} \ead{jmgallardo@uma.es}
\author{Kleiton A. Schneider} \ead{ks@uma.es}
\author{Manuel J. Castro} \ead{mjcastro@uma.es}

\cortext[cor]{Corresponding author}

\address{Departamento de An\'alisis Matem\'atico, Estad\'{\i}stica e Investigaci\'on Operativa, y Matem\'atica Aplicada, Universidad de M\'alaga, Campus de Teatinos s/n, M\'alaga 29080, Spain}

\tnotetext[t1]{This research has been partially supported by the Spanish Government Research project MTM2015-70490-C2-1R. The second author has been partially funded by UFMS, Brazil. The numerical computations have been performed at the Laboratory of Numerical Methods of the University of M\'alaga.}

\begin{abstract}
We propose a general class of genuinely two-dimensional incomplete Riemann solvers for systems of conservation laws. In particular, extensions of Balsara's multidimensional HLL scheme \cite{Balsara2012} to two-dimensional PVM/RVM (Polynomial/Rational Viscosity Matrix) finite volume methods are considered. The numerical flux is constructed by assembling, at each edge of the computational mesh, a one-dimensional PVM/RVM flux with two purely two-dimensional PVM/RVM fluxes at vertices, which take into account transversal features of the flow.

The proposed methods are applicable to general hyperbolic systems, although in this paper we focus on applications to magnetohydrodynamics. In particular, we propose an efficient technique for divergence cleaning of the magnetic field that provides good results when combined with our two-dimensional solvers. Several numerical tests including genuinely two-dimensional effects are presented to test the performances of the proposed schemes.
\end{abstract}


\end{frontmatter}

%%%%%

\section{Introduction} \label{sec:intro}

Since the pioneering work of Godunov \cite{Godunov1959}, Riemann solvers have been an important ingredient in the design of robust and accurate numerical methods for hyperbolic conservation laws. Usually, the exact solution of a Riemann problem contains many complex features which makes it unsuitable in practice. For this reason, a number of \emph{incomplete} Riemann solvers have been devised in the literature (Lax-Friedrichs, Rusanov, HLL, FORCE, etc.; see \cite{PVM2012} and the references therein), which takes into account only some of the waves appearing in the Riemann fan.

In the paper \cite{PVM2012} the authors introduced the so-called PVM methods, which constitute a general class of incomplete Riemann solvers defined in terms of viscosity matrices based on the polynomial evaluation of a given Roe matrix or the Jacobian of the flux at some average state. Indeed, they showed that many classical schemes in the literature can be viewed as particular cases of PVM solvers. The idea behind PVM methods was later extended in \cite{RVM2014} to the case of rational functions, giving rise to the class of RVM methods, which provide a much more precise representation of internal waves. In what follows, we will use the term AVM (Approximate Viscosity Matrix) to refer to both PVM and RVM solvers. In \cite{RVM2014} it was also shown that the choice of a precise first-order solver is important when designing higher order methods in terms of computational efficiency, at least for solutions involving complex structures. In particular, the choice of appropriate underlying functions for an AVM solver allows to control the amount of numerical difussion of the resulting scheme. An important feature of AVM methods is that they do not need the whole spectral decomposition of the system, but only a bound on its spectral radius. This makes this kind of solvers particularly efficient for treating complex problems: see \cite{PVM2012,RVM2014} for applications to magnetohydrodynamics (MHD) and multilayer shallow water systems. Further applications to relativistic MHD where considered in \cite{RMHD2017}, where a Jacobian-free version of some PVM solvers was also introduced. Furthermore, in \cite{RVM_Osher2016} the idea of AVM solvers was applied to build approximate versions of the DOT method \cite{DOT2011}, which in turn can be viewed as approximations of the classical Osher-Solomon scheme.

When solving multidimensional problems, it is a common practice to consider one-dimensional projected Riemann solvers on the edges of the cells of the computational domain. However, many researchers consider that within this approach, one-dimensional solvers lose much of their efficiency, mainly due to the fact that they do not take into account features of the solution propagating transversally to the cells boundaries. Moreover, the biasing introduced by one-dimensional solvers produces a reduction of the permissible Courant number when simulating multidimensional flows. For these reasons, there has been many attempts in the literature to build genuinely multidimensional Riemann solvers (see \cite{Balsara2010,Balsara2012,VNA2015} and the references therein; see also \cite{Roe2017} for a more detailed account on multidimensional upwinding), although most of them are specifically designed for the Euler or MHD equations.

In \cite{Wendroff1999} Wendroff proposed a two-dimensional Riemann solver for the Euler equations which extended the one-dimensional HLL method \cite{HLL1983}, in which the two-dimensional interactions at cell corners were taken into account through the approximate solutions of two-dimensional Riemann problems. However, a drawback of Wendroff's solver is that it does not possess explicit expressions allowing its direct implementation; also, high-order extensions of the solver do not seem straighforward to build. Some years later, Balsara \cite{Balsara2010} modified Wendroff's formulation and proposed a two-dimensional HLL solver for the Euler and MHD equations on structured meshes, which included closed forms of the fluxes and allowed for an easier high-order extension. A more robust version of Balsara's solver was later proposed in \cite{Balsara2012}, and an extension to unstructured meshes was presented in \cite{BDA2014}. Another interesting reformulation of Wendroff's approach has been recently proposed in \cite{VNA2015}.

In this paper we introduce new types of genuinely two-dimensional incomplete Riemann solvers, which can be viewed as general AVM extensions of Balsara's HLL solver \cite{Balsara2012}. To achieve this, first we propose a reinterpretation of Balsara's solver in the particular case in which a simple four-wave model is considered for the two-dimensional corrections at cell vertices. These corrections can be viewed as suitable combinations of one-dimensional numerical fluxes of HLL type. Once the role of the underlying one-dimensional HLL solver is clearly identified, it can be changed by an appropriate one-dimensional AVM numerical flux. This leads to two-dimensional contributions of AVM type at the cell corners, which take into account the transversal features of the flow. It is noteworthy that the proposed solvers are applicable to general hyperbolic systems of conservation laws; in particular, they are not restricted to the Euler or MHD equations. An additional feature of these multidimensional AVM solvers is that they are theoretically stable up to a CFL number of unity. 

Several families of basis functions for AVM methods were considered in \cite{PVM2012,RVM2014}; in particular, good results were obtained with Chebyshev polynomials and Newman rational functions. In this paper we will consider PVM solvers based on the internal polynomials previously studied in \cite{RMHD2017}, together with new families of RVM solvers based on suitable Pad\'e approximations, which enjoy good stability properties. An additional advantage of these solvers is that no entropy fix is needed in the presence of sonic points.

We will focus in this paper on the comparison of the proposed methods in the first-order case, as done previously for one-dimensional RVM methods in \cite{RVM2014}. As commented before, the choice of accurate first-order solvers is important when analizing complex scenarios. On the other hand, the use of our solvers as building blocks for higher order methods will be considered in a forthcoming paper.

In order to compare with other methods in the literature, applications to MHD equations will be considered. An additional difficulty that appears in the numerical resolution of MHD is the divergence-free constraint on the magnetic field, that has to be imposed in order to ensure the accuracy and stability of the numerical schemes. Several divergence cleaning techniques has been proposed in the literature: see \cite{Toth2000} for a complete account on this topic. Regarding this, we present a novel technique for divergence cleaning based on Powell's \cite{Powell1994} formulation of MHD equations. This technique has been included in our two-dimensional AVM solvers, modifying the numerical fluxes accordingly to \cite{MCP2014}. Numerical experiments show that the results of our methods are comparable to those obtained with the well-known projection method \cite{BB1980}, both in precision and computational time.

As remarked in \cite{Gosse2014}, the inclusion of source terms in Balsara's solver seems to be a nontrivial matter. A way to handle source or coupling terms is to reformulate the system in nonconservative form, and then to apply the theory of path-conservative numerical schemes \cite{Pares2006}. Taking into account the form of our multidimensional AVM methods, it is possible to extend them to the nonconservative case following the guidelines in \cite{PVM2012,RVM2014}. This extension will be analized in a future paper, where we will construct high-order multidimensional AVM schemes that will be applied, in particular, for solving multilayer shallow water systems including depth variations.

The paper is organized as follows. After some preliminaries in Section 2, Balsara's multidimensional HLL solver is recalled in Section 3. Then, a brief review of AVM solvers is given in Section 4, together with the definitions of internal polynomials and Pad\'e approximants. The core of the paper is Section 5, where multidimensional AVM solvers are introduced. Applications to MHD, including the proposed divergence cleaning technique, is considered in Section 6. Numerical experiments are then presented in Section 7. Finally, some conclusions are drawn in Section 8.

%%%%%

\section{Preliminaries} \label{sec:preliminaries}

Consider a two-dimensional hyperbolic system of conservation laws of the form
\begin{equation}\label{eq:2dsystem}
\partial_t U +\partial_x F(U)+\partial_y G(U) = 0,
\end{equation}
where $U$ is defined on $\Omega\times [0,T]$, being $\Omega\subset\R^2$ a domain, and takes values on an open convex set $\mathcal{O}\subset\R^N$; the physical fluxes $F, G\colon\mathcal{O}\to\R^N$ are assumed to be smooth functions. 

We are interested in the numerical solution of \eqref{eq:2dsystem} by means of finite volume methods on structured meshes. The spatial domain is divided into rectangular cells $C_{ij}=[x_{i-1/2}, x_{i+1/2}] \times [y_{j-1/2}, y_{j+1/2}]$ of size $\Delta x\times\Delta y$ and center $(x_i, y_j)$, and the time step is denoted as $t^n=n\Delta t$. Integrating \eqref{eq:2dsystem} over the control volume $C_{ij} \times [t^n, t^{n+1}]$ yields
\begin{equation} \label{eq:fv2d}
U_{ij}^{n+1} = U_{ij}^n-\frac{\Delta t}{\Delta x}(F_{i+1/2,j}^n-F_{i-1/2,j}^n) - \frac{\Delta t}{\Delta y}(G_{i,j+1/2}^n-G_{i,j-1/2}^n),
\end{equation}
where
\[
U_{ij}^n \approx \frac{1}{|C_{ij}|} \int_{C_{ij}}U(x,y,t^n)\,dx,
\]
and the numerical fluxes verify
\[
F_{i+1/2,j}^n \approx \frac{1}{\Delta t} \int_{t^n}^{t^{n+1}} \int_{y_{j-1/2}}^{y_{j+1/2}} F(U(x_{i+1/2}, y, t))\,dy\,dt,
\]
and
\[
G_{i, j+1/2}^n \approx \frac{1}{\Delta t} \int_{t^n}^{t^{n+1}} \int_{x_{i-1/2}}^{x_{i+1/2}} G(U(x, y_{j+1/2}, t))\,dx\,dt.
\]
In what follows, time dependence will be dropped unless necessary.

Usual definitions of the numerical fluxes are based on one-dimensional Riemann solvers in the  coordinate directions, which lead to expressions of the form $F_{i+1/2,j}=\mathcal{F}(U_{ij}, U_{i+1,j})$ and $G_{i,j+1/2}=\mathcal{G}(U_{ij}, U_{i,j+1})$ for certain flux functions $\mathcal{F}(U,V)$ and $\mathcal{G}(U,V)$. However, this approach does not take into account the genuinely two-dimensional effects that can be produced at each of the four vertices $(x_{i\pm 1/2}, y_{j\pm 1/2})$ of a cell. A way to introduce these contributions in the definition of the numerical fluxes was proposed by Balsara in \cite{Balsara2010} and \cite{Balsara2012} for a HLL-type scheme. Specifically, the two-dimensional interactions at a given vertex are taken into account by considering the approximate solution of the two-dimensional Riemann problem that has as initial condition the four constant neighboring states at the vertex. Assuming that the strong interaction area in the solution of the two-dimensional Riemann problem is rectangular, Balsara was able to compute closed-form expressions for the numerical fluxes, thus improving previous work done in this direction by Wendroff \cite{Wendroff1999}. For the sake of completeness, and to introduce some notation, Balsara's solver will be reviewed in the next section.

%%%%%

\section{Balsara's multidimensional HLL solver} \label{sec:Balsara_solver}

The purpose of this section is to give a brief review of Balsara's HLL multidimensional solver. Although all the details can be found in \cite{Balsara2010,Balsara2012} (see also \cite{VNA2015}), we give here the information that will be needed later to construct our schemes in Section \ref{sec:AVM2d}.

\begin{figure}[!ht]
	\begin{center}
		\includegraphics[width=0.7\textwidth]{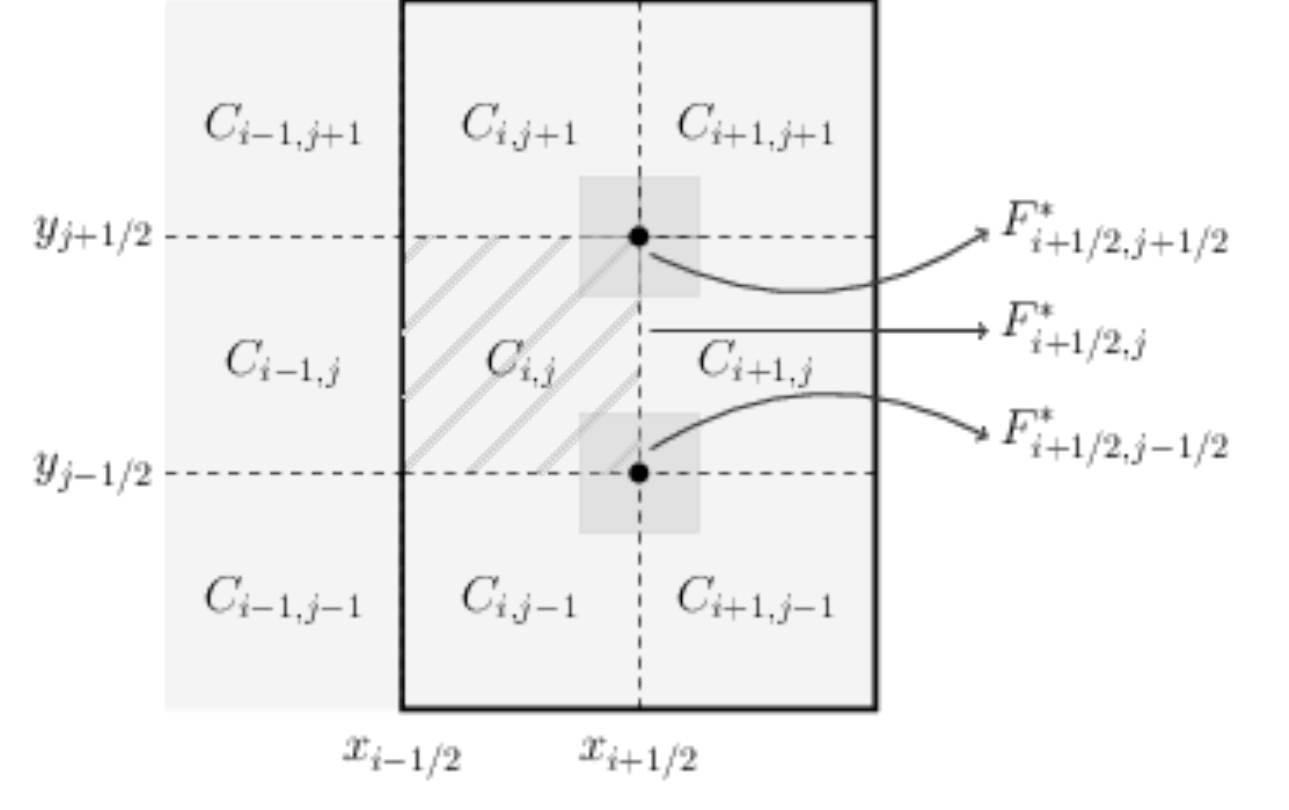}
		\caption{Stencil used to build the numerical flux $F_{i+1/2, j}$.} \label{fig:fig1}
	\end{center}
\end{figure}

To define the numerical flux at a given edge, Balsara's solver considers a linear convex combination of a one-dimensional HLL flux and two multidimensional contributions at the vertices. To be more precise, let us consider for a given cell $C_{ij}$ the six-cell stencil depicted in Figure \ref{fig:fig1}. The flux in the $x$-direction through the right edge can be expressed as
\begin{equation} \label{eq:2dassembledflux}
F_{i+1/2,j} = \alpha F_{i+1/2,j+1/2}^*+\beta F_{i+1/2,j}^*+\gamma F_{i+1/2,j-1/2}^*,
\end{equation}
where $F_{i+1/2,j}^*$ is the usual one-dimensional HLL flux between the cells $C_{ij}$ and $C_{i+1,j}$, and $F_{i+1/2,j\pm 1/2}^*$ represent the two-dimensional contributions at the vertices $(x_{i+1/2}, y_{j\pm 1/2})$; proper choices of the coefficients $\alpha$, $\beta$ and $\gamma$ will be given later. Similar expressions can be stated for the $y$-flux $G_{i, j+1/2}$.

\begin{figure}[!ht]
	\begin{center}
		\includegraphics[width=0.384\textwidth]{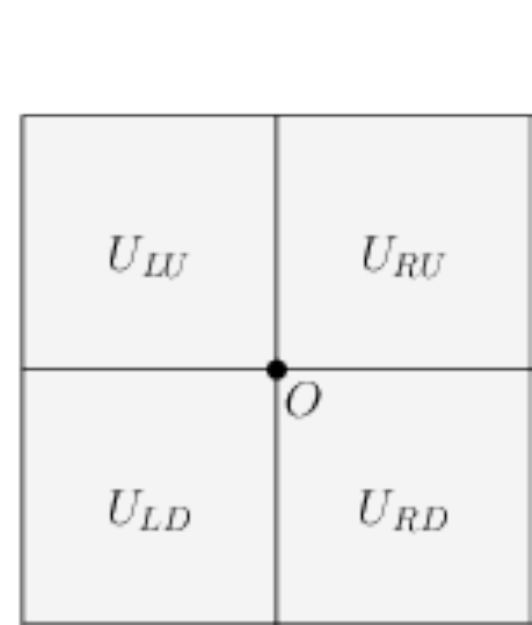} \qquad
		\includegraphics[width=0.5\textwidth]{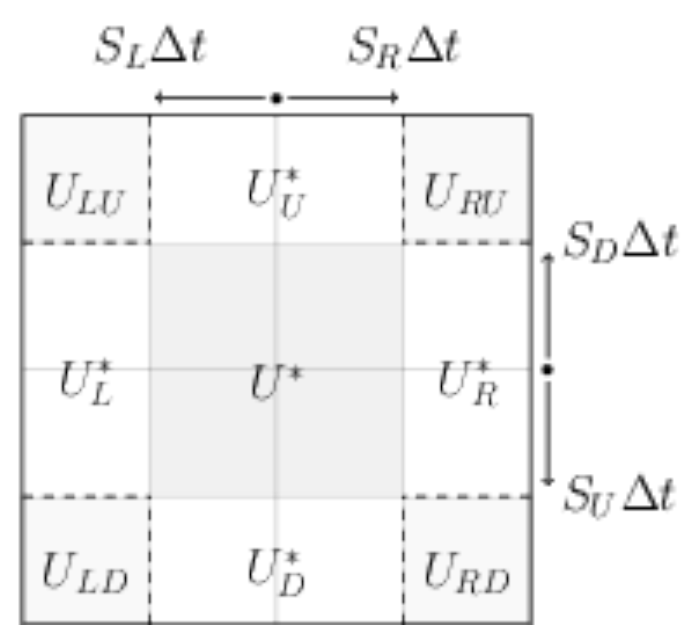}
		\caption{Left: Local stencil. Right: Structure of the solution of the Riemann problem.} \label{fig:fig2}
	\end{center}
\end{figure}

To define the two-dimensional flux $F_{i+1/2,j+1/2}^*$, Balsara considered a two-dimensional Riemann problem with initial condition given by the constant states at each of the four cells which have $(x_{i+1/2}, y_{j+1/2})$ as a common vertex. To simplify the notation, let us consider the local stencil depicted in Figure \ref{fig:fig2} (left); for example, for $O\equiv (x_{i+1/2}, y_{j+1/2})$ the cell values would be given by $U_{LD}=U_{ij}$, $U_{RD}=U_{i+1, j}$, $U_{LU}=U_{i, j+1}$ and $U_{RU}=U_{i+1, j+1}$. The initial condition for the local Riemann problem is then
\begin{equation}\label{eq:iniRP}
U(x,y,t^0)=
\begin{cases}
U_{LD} & \text{if } x<0,\  y<0, \\
U_{RD} & \text{if } x>0,\  y<0, \\
U_{LU} & \text{if } x<0,\  y>0, \\
U_{RU} & \text{if } x>0,\  y>0.
\end{cases}
\end{equation}

It is known that problem \eqref{eq:2dsystem}-\eqref{eq:iniRP} has a self-similar solution roughly consisting of four directional one-dimensional Riemann problems arising at common edges, along with a region of strong interaction where complex structures appear; see, for example, \cite{PL1998, SCG1993, ZZ1990} for more detailed descriptions. In the design of numerical methods, however, it is impractical to resolve all the fine structures arising in the strong interaction region. Thus, this region is averaged into a single constant state, in a similar way to the intermediate state arising in the one-dimensional HLL solver. In Balsara's formulation, the region of strong interaction is assumed to be rectangular, although other options are possible (see, e. g., \cite{VNA2015}).

To build a two-dimensional Riemann solver of HLL type, first we need to identify the maximal speeds in the coordinate directions. Following the notation in \cite{Balsara2012}, let $S_L^U$ and $S_R^U$ denote the fastest left and right speeds at the edge between $U_{LU}$ and $U_{RU}$ (see Figure \ref{fig:fig2}, right). An usual choice for $S_R^U$ and $S_L^U$ are given by
\begin{equation} \label{eq:SRU_SLU}
\begin{split}
& S_R^U = \max(\lambda_x^N(U_{RU}), \bar{\lambda}_x^N(U_{LU}, U_{RU})), \\
& S_L^U = \min(\lambda_x^1(U_{LU}), \bar{\lambda}_x^1(U_{LU}, U_{RU})),
\end{split}
\end{equation}
where $\lambda_x^1(U_{LU})$ and $\lambda_x^N(U_{RU})$ denote the maximal left- and right-going wave speeds associated to the states $U_{LU}$ and $U_{RU}$, respectively, while $\bar{\lambda}_x^1(U_{LU}, U_{RU})$ and $\bar{\lambda}_x^N(U_{LU}, U_{RU})$ are the maximal left- and right-going speeds arising from a linearized Riemann solver between the states $U_{LU}$ and $U_{RU}$. Similar definitions apply for $S_L^D$ and $S_R^D$. Speeds in the $y$-direction, $S_D^L$, $S_D^R$, $S_U^L$ and $S_U^R$, are defined in an analogous way. Finally, the strong interaction region is assumed to be bounded by the maximal left and right wave speeds given by
\[
S_L=\min(S_L^D, S_L^U), \quad S_R=\max(S_R^D, S_R^U),
\]
and the maximal downward and upward speeds
\[
S_D=\min(S_D^L, S_D^R), \quad S_U=\max(S_U^L, S_U^R).
\]
A simplifying assumption can be made at this point, which consists in considering only the four speeds $S_L$, $S_R$, $S_D$ and $S_U$ in the wave structure of the solution. That is, we assume that $S_\alpha=S_\alpha^D=S_\alpha^U$ for $\alpha=L,R$, and $S_\beta=S_\beta^L=S_\beta^R$ for $\beta=D,U$. This leads to compact expressions for the strongly interacting state and its associated numerical fluxes. As later in Section \ref{sec:AVM2d} we will be interested on a four-wave model, this assumption will be made in what follows.

Assume now that we are in the subsonic case in both directions, that is, $S_L<0<S_R$ and $S_D<0<S_U$. Figure \ref{fig:fig2} (right) shows the proposed wave structure at a time $t=\Delta t$. The one-dimensional Riemann problems arising at the edges give rise to four HLL states; thus, for 
$\alpha=D, U$, the states $U_{L\alpha}$ and $U_{R\alpha}$ give
\begin{equation} \label{eq:UDstates}
U_\alpha^* = \frac{S_R U_{R\alpha}-S_L U_{L\alpha}+F_{L\alpha}-F_{R\alpha}}{S_R-S_L},
\end{equation}
while for $\beta=L, R$, the states $U_{\beta D}$ and $U_{\beta U}$ are
\begin{equation} \label{eq:LRstates}
U_\beta^* = \frac{S_U U_{\beta U}-S_D U_{\beta D}+G_{\beta D}-G_{\beta U}}{S_U-S_D},
\end{equation}
where $F_{\beta\alpha}$ and $G_{\beta\alpha}$ denote, respectively, the fluxes $F(U_{\beta\alpha})$ and $G(U_{\beta\alpha})$. The corresponding HLL fluxes are then given by
\[
F_\alpha^* = \frac{S_R F_{L\alpha}-S_L F_{R\alpha}+S_L S_R (U_{R\alpha}-U_{L\alpha})}{S_R-S_L}, \quad \alpha=D, U,
\]
and
\[
G_\beta^* = \frac{S_U G_{\beta D}-S_D G_{\beta U}+S_D S_U(U_{\beta U}-U_{\beta D})}{S_U-S_D}, \quad \beta=L, R.
\]

\begin{figure}[!ht]
	\begin{center}
		\includegraphics[width=0.8\textwidth]{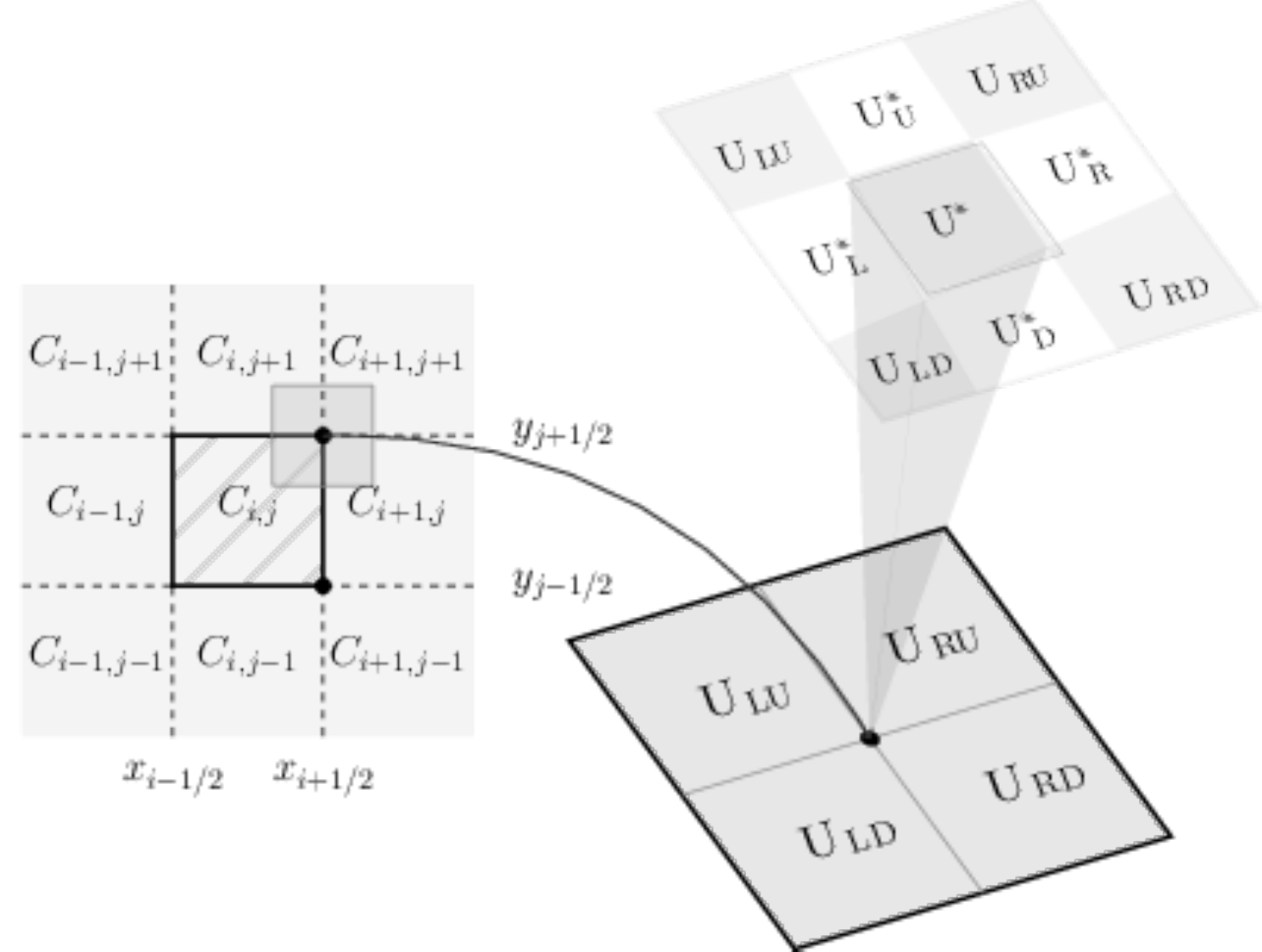}
		\caption{Structure of the solution of the 2d Riemann problem at a vertex.} \label{fig:fig3}
	\end{center}
\end{figure}

Now, integration of \eqref{eq:2dsystem} on the control volume $(S_L\Delta t, S_R\Delta t) \times (S_D\Delta t, S_U\Delta t) \times (0, \Delta t)$ leads to an explicit expression for the constant resolved state $U^*$, namely
\[
\begin{split}
U^* = & \frac{S_R S_U U_{RU}+S_L S_D U_{LD}-S_R S_D U_{RD}-S_L S_U U_{LU}}{(S_R-S_L)(S_U-S_D)} \\
& -\frac{(F_{RU}-F_{LU})S_U-(F_{RD}-F_{LD})S_D}{2(S_R-S_L)(S_U-S_D)} \\
& -\frac{(G_{RU}-G_{RD})S_R-(G_{LU}-G_{LD})S_L}{2(S_R-S_L)(S_U-S_D)} \\
& -\frac{F_R^*-F_L^*}{2(S_R-S_L)} - \frac{G_U^*-G_D^*}{2(S_U-S_D)}.
\end{split}
\]
Figure \ref{fig:fig3} shows the considered structure of the solution of the two-dimensional Riemann problem at a vertex. Notice that all the terms appearing in the definition of $U^*$ are known, except the \emph{transverse fluxes} appearing in the last two terms, which have to be defined taking into account the corresponding intermediate state and the normal flux. For example, $F_R^*$ is constructed using the state $U_R^*$ and the flux $G_R^*$ previously defined; in particular, notice that $F_R^*\neq F(U_R^*)$ in general. Balsara offers an explicit solution in \cite{Balsara2012} for the Euler equations, although similar constructions may be done for other systems: see Section \ref{sec:MHD} for the case of MHD equations.

To compute the resolved flux $F^*$, equation \eqref{eq:2dsystem} is integrated on the subvolume $(0, S_R\Delta t) \times (S_D\Delta t, S_U\Delta t) \times (0, \Delta t)$, yielding
\begin{equation} \label{eq:2dfluxF}
\begin{split}
F^* = & \frac{S_U}{S_U-S_D}F_U^* - \frac{S_D}{S_U-S_D}F_D^* \\
& - \frac{S_L S_R}{(S_R-S_L)(S_U-S_D)}(G_{RU}-G_{LU}+G_{LD}-G_{RD}) \\
& - \frac{S_L S_U (F_{RU}-F_R^*)-S_R S_U (F_{LU}-F_L^*)}{(S_R-S_L)(S_U-S_D)} \\
& - \frac{S_R S_D (F_{LD}-F_L^*)-S_L S_D (F_{RD}-F_R^*)}{(S_R-S_L)(S_U-S_D)}.
\end{split}
\end{equation}
Notice again that, in general, $F^*\neq F(U^*)$. Similarly, integration of \eqref{eq:2dsystem} on $(S_L\Delta t, S_R\Delta t) \times (0, S_U\Delta t) \times (0, \Delta t)$ leads to
\begin{equation} \label{eq:2dfluxG}
\begin{split}
G^* = & \frac{S_R}{S_R-S_L}G_R^* - \frac{S_L}{S_R-S_L}G_L^* \\
& - \frac{S_D S_U}{(S_R-S_L)(S_U-S_D)}(F_{RU}-F_{LU}+F_{LD}-F_{RD}) \\
& - \frac{S_R S_D (G_{RU}-G_U^*)-S_L S_D (G_{LU}-G_U^*)}{(S_R-S_L)(S_U-S_D)} \\
& - \frac{S_L S_U (G_{LD}-G_D^*)-S_R S_U (G_{RD}-G_D^*)}{(S_R-S_L)(S_U-S_D)}.
\end{split}
\end{equation}
As stated before, the fluxes \eqref{eq:2dfluxF} and \eqref{eq:2dfluxG} are used when the flow is subsonic in both spatial directions. When the flow is supersonic in some of the spatial directions, the corresponding expressions for $F^*$ and $G^*$ may be defined as the corresponding upwinded fluxes; the details can be found in \cite[Sec. 2.1]{Balsara2012}.

\begin{figure}[!ht]
	\begin{center}
		\includegraphics[width=0.65\textwidth]{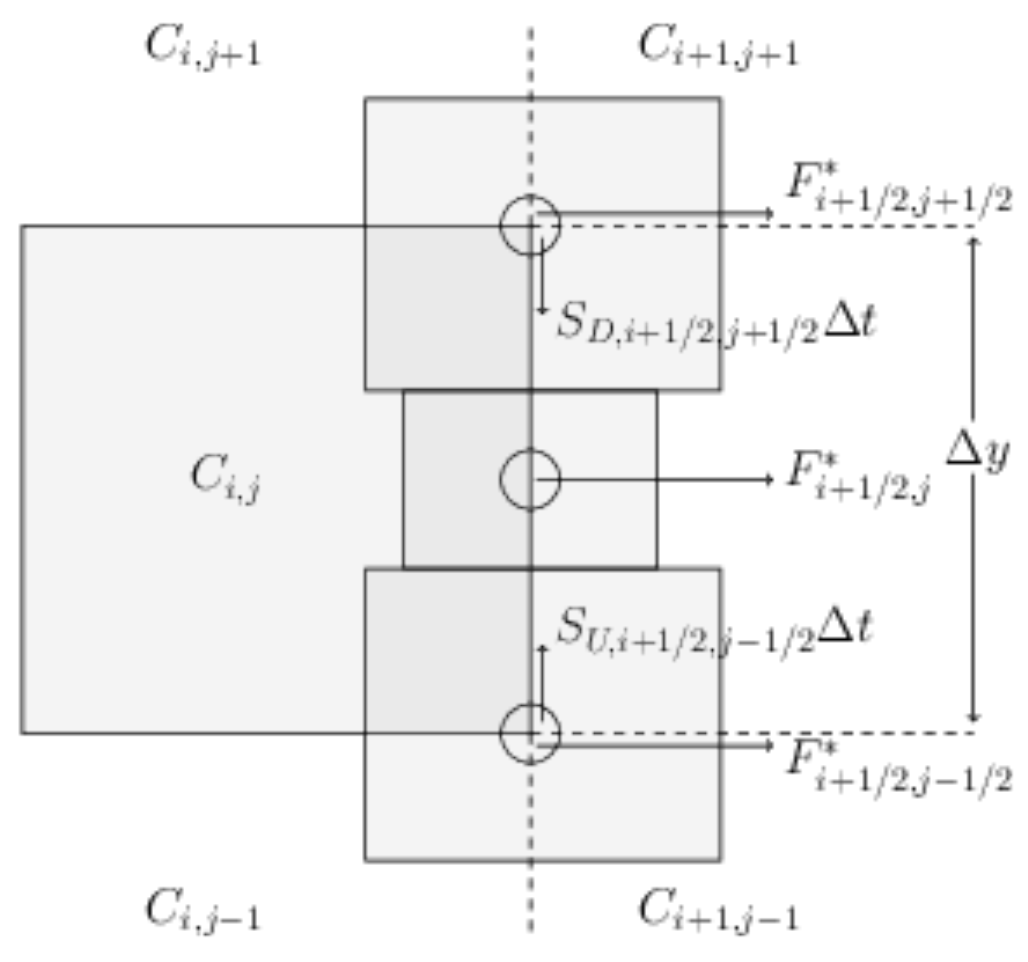}
		\caption{Assembling of one- and two-dimensional fluxes at an edge.} \label{fig:fig4}
	\end{center}
\end{figure}

Finally, the assembling \eqref{eq:2dassembledflux} is obtained from a time-space average of the normal flux through the edge (see Figure \ref{fig:fig4}). This leads to the coefficients
\[
\alpha = -S_{D, i+1/2, j+1/2}\frac{\Delta t}{2\Delta y}, \quad \gamma = S_{U,i+1/2,j-1/2}\frac{\Delta t}{2\Delta y}
\]
and
\[
\beta = 1-(S_{U,i+1/2,j-1/2}-S_{D, i+1/2, j+1/2})\frac{\Delta t}{2\Delta y},
\]
with the obvious notations. Strictly speaking, these values are valid only for the subsonic case in both directions. However, supersonic cases can be accommodated by simply substituting $S_{U,R}$ by $S_{U,R}^+=\max(S_{U,R}, 0)$ and $S_{D,L}$ by $S_{D,L}^-=\min(S_{D,L}, 0)$. Another option is to use Simpson's rule to approximate the integral, which leads to the simpler choice
\[
\alpha=\gamma=\frac{1}{6}, \quad \beta=\frac{2}{3}.
\]
We have chosen this approach in the numerical experiments, although no significative differences have been found when the first option is used. More accurate approaches are possible: see, e. g., \cite[Sect. 3.3.2]{VNA2015}.

With respect to the CFL condition needed to ensure stability, the multidimensional corrections allow for a maximum CFL number of unity: see \cite{Balsara2010} for a complete account on this point.

%%%%%

\section{Some notes on AVM solvers} \label{sec:AVMsolvers}

In Section \ref{sec:AVM2d} we will see how to construct more precise versions of Balsara's multidimensional scheme by changing the underlying HLL solver in a suitable way. To this end, and for the sake of completeness, we give in this section a brief survey on AVM solvers.

For simplicity, in this part will be enough to focus on one-dimensional systems. Therefore, we are interested in the numerical solution of the hyperbolic system
\[
\partial_t U+\partial_x F(U) = 0
\]
by means of a finite volume method of the form
\[
U_i^{n+1} = U_i^n-\frac{\Delta t}{\Delta x}(F_{i+1/2}-F_{i-1/2}),
\]
with numerical flux given by
\begin{equation} \label{eq:AVMnumflux}
F_{i+1/2} = \frac{F(U_i)+F(U_{i+1})}{2}-\frac{1}{2}Q_{i+1/2}(U_{i+1}-U_i).
\end{equation}
Here $Q_{i+1/2}$ denotes the numerical \emph{viscosity matrix}, which determines the numerical diffusion of the scheme. It is worth noticing that Roe's method \cite{Roe1981} can be written in the form \eqref{eq:AVMnumflux} by considering the viscosity matrix $Q_{i+1/2}=|A_{i+1/2}|$, being $A_{i+1/2}$ a Roe matrix for the system. The idea of PVM (Polynomial Viscosity Matrix) solvers introduced in \cite{PVM2012} consists in approximating $|A_{i+1/2}|$ by means of a suitable polynomial evaluation of $A_{i+1/2}$. More precisely, let $p(x)$ be a polynomial approximation of $|x|$ in the interval $[-1, 1]$, and let $\lambda_{i+1/2,\text{max}}$ be the eigenvalue of $A_{i+1/2}$ with maximum modulus (or an upper bound of it). Then the viscosity matrix of the PVM method associated to $p(x)$ is given by
\begin{equation}\label{eq:viscosity_matrix}
Q_{i+1/2} = |\lambda_{i+1/2,\text{max}}|p(|\lambda_{i+1/2,\text{max}}|^{-1}A_{i+1/2}).
\end{equation}
Notice that the best $p(x)$ approaches $|x|$, the closer the behaviour of the associated PVM scheme will be to that of Roe's method. A fundamental issue is that the spectral decomposition of $A_{i+1/2}$ is not needed in the construction of a PVM method, but only a bound on its spectral radius. This feature makes PVM methods greatly efficient and applicable to systems in which the eigenstructure is not known or is difficult to obtain. In the cases in which a Roe matrix is not available or is costly to compute, $A_{i+1/2}$ can be taken as the Jacobian matrix of the system evaluated at some average state. Several well-known schemes in the literature belong to the family of PVM methods: Lax-Friedrichs, Rusanov, HLL, FORCE, Roe, etc. (see \cite{PVM2012}). 

It is known that rational functions may provide much more precise approximations to $|x|$ than polynomials. The RVM methods introduced in \cite{RVM2014} take advantage of this fact, following the PVM idea but using rational functions as basis instead of polynomials. In what follows we will encompass both kind of solvers, PVM and RVM, under the common term AVM (Approximate Viscosity Matrix) solvers. AVM solvers have also been used in \cite{RVM_Osher2016} to construct efficient versions of the Osher-Solomon method, following the ideas in \cite{DOT2011}. Recently, in \cite{RMHD2017} some types of AVM solvers have been implemented in Jacobian-free form; this is particularly interesting when solving complex systems as, for example, relativistic MHD.  

Finally, we must remark that the stability of an AVM scheme strongly depends on the properties of the underlying function $f(x)$. In particular, it must verify the \emph{stability condition}
\begin{equation}\label{eq:stability_cond}
|x| \leq f(x) \leq 1, \quad \forall\, x\in [-1, 1].
\end{equation}

With respect to the choice of solvers, in \cite{RVM2014} the authors proposed a family of PVM methods based on Chebyshev polynomials, which provide optimal approximations to the absolute value function, and could also be implemented in Jacobian-free form. In the same paper, RVM solvers based on Newman rational functions \cite{Newman1964} where found to be the more efficient choice: they provided similar results as Roe's method for complex problems in MHD and multilayer shallow water systems, but with a much smaller computational cost. Despite their efficiency, a drawback of Chebyshev and Newman AVM solvers is that they do not satisfy the stability condition \eqref{eq:stability_cond} strictly, so a slight modification has to be performed which may perturb the computation of the maximal speeds of propagation. On the other hand, in \cite{RMHD2017} a new family of internal polynomial approximations to $|x|$ verifying \eqref{eq:stability_cond} was proposed. Here we will consider this family of internal polynomials together with a class of rational Pad\'e approximations based on \cite{KL1991}, which also satisfy \eqref{eq:stability_cond}. For the sake of completeness, we briefly describe them in what follows.

\begin{figure}[!ht]
\begin{center}
\includegraphics[width=0.49\textwidth]{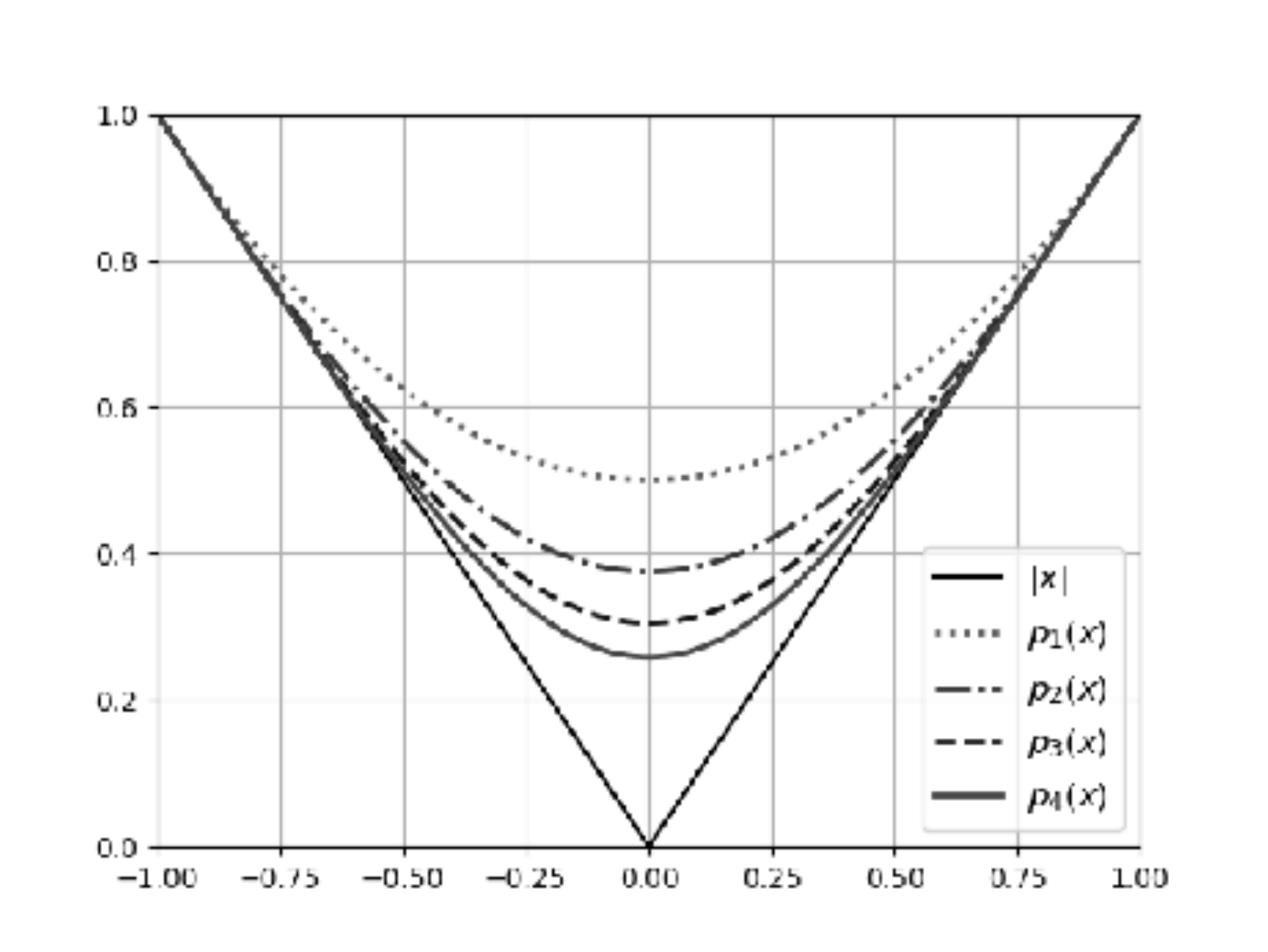}
\includegraphics[width=0.49\textwidth]{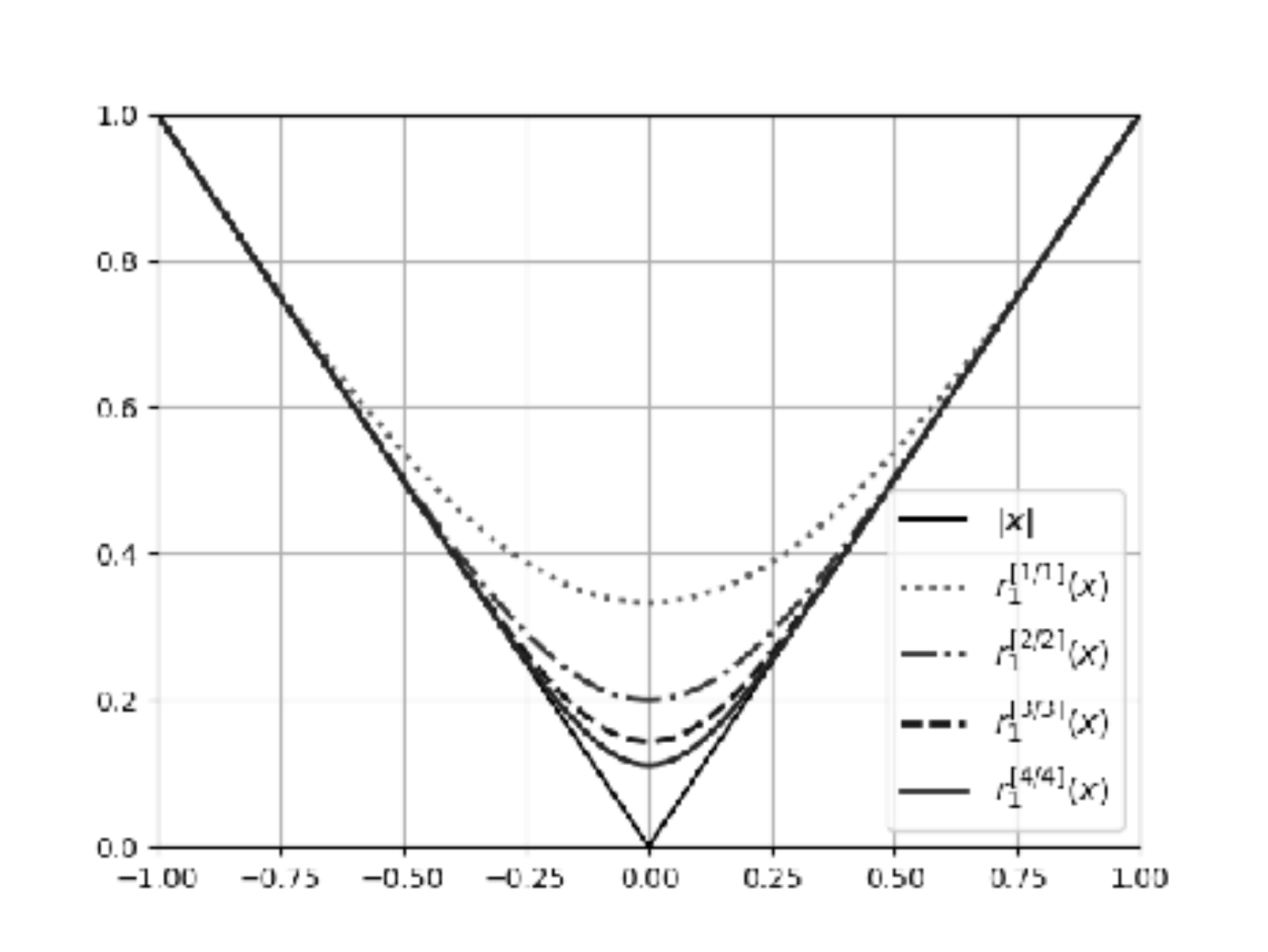} \\
\includegraphics[width=0.49\textwidth]{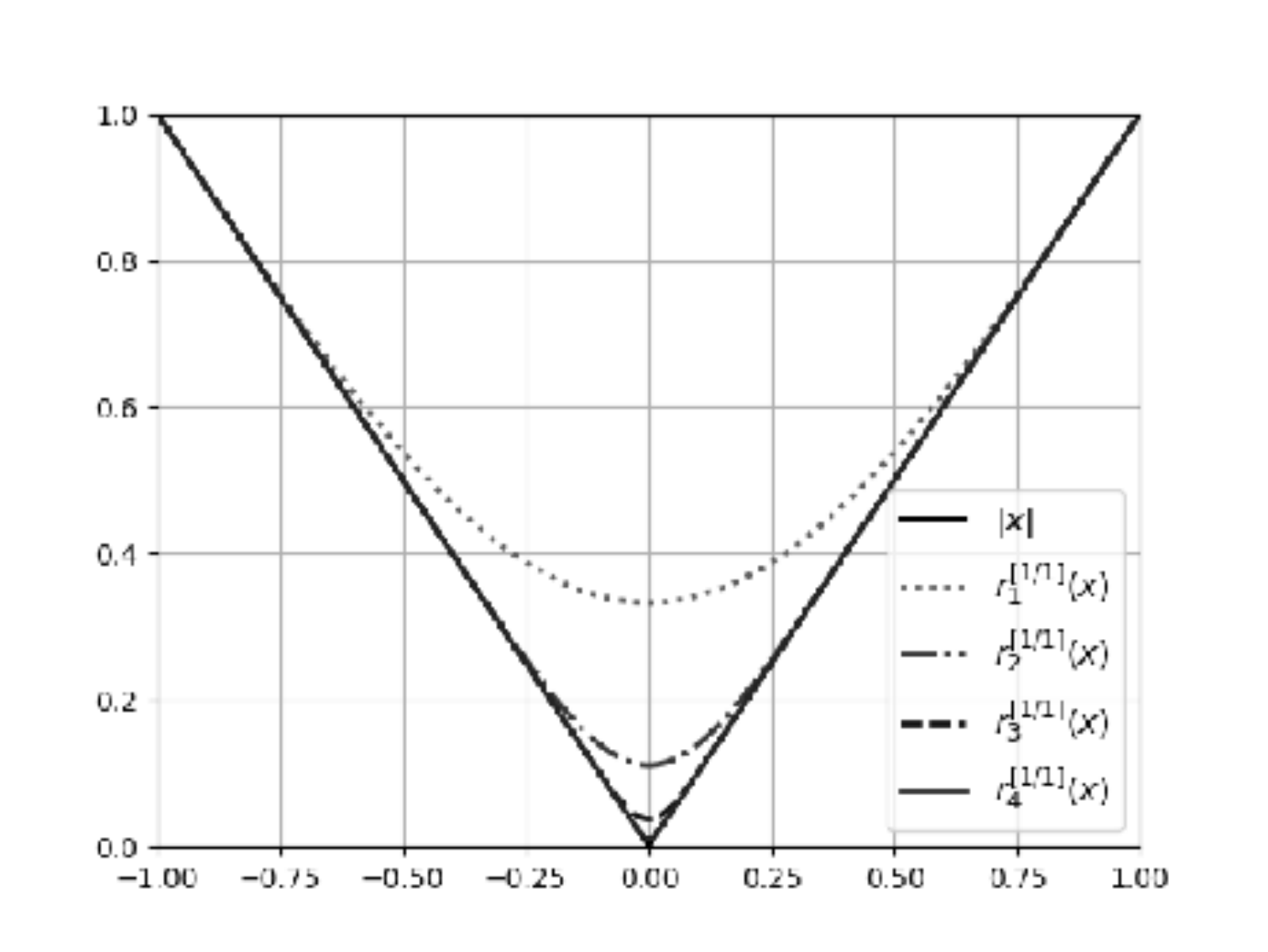}
\includegraphics[width=0.49\textwidth]{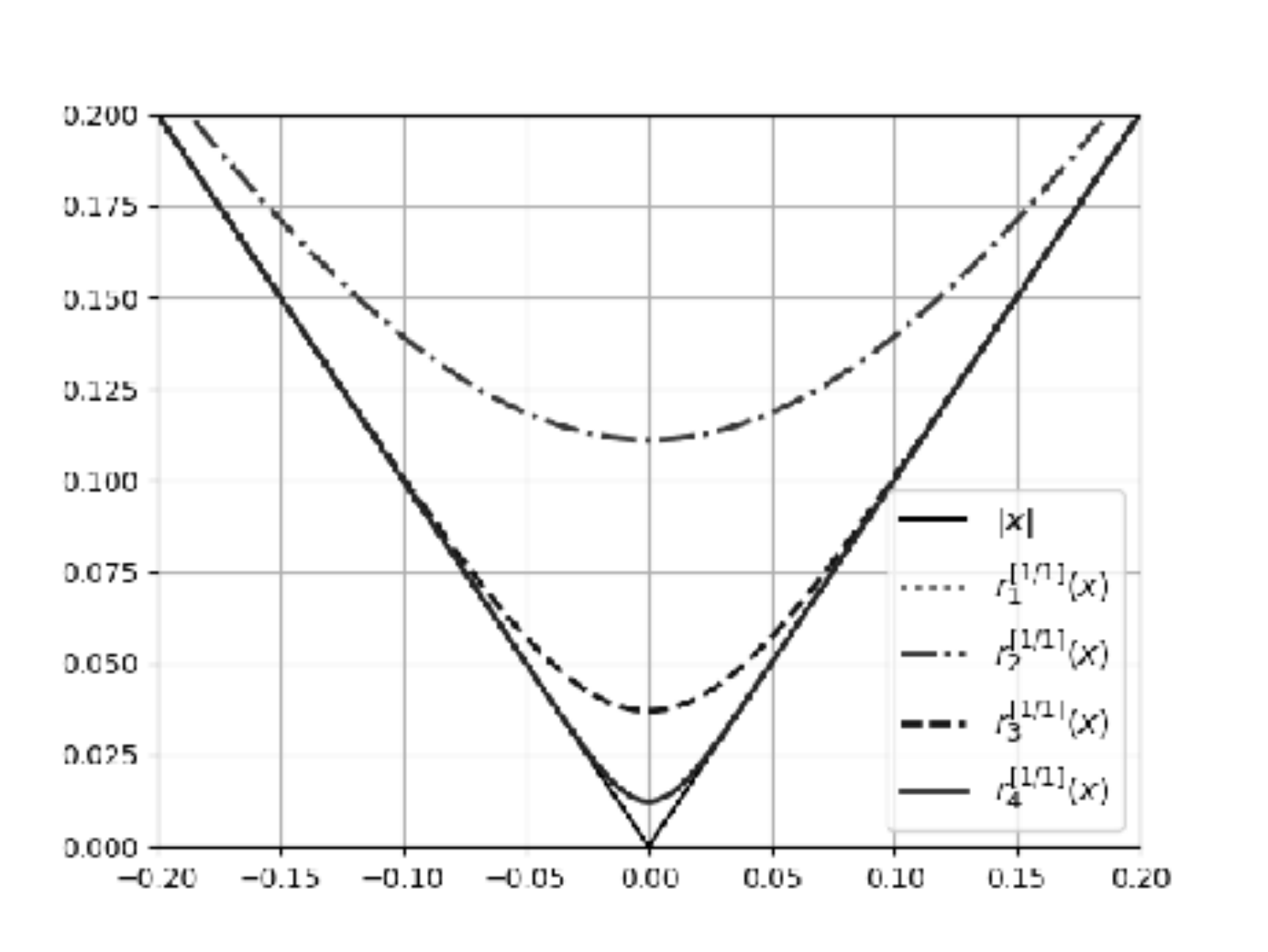}
\caption{Top: Internal polynomials $p_n(x)$ for $n=1, 2, 3, 4$, and Pad\'e rational approximants $r_1^{[m/k]}(x)$ for $m=k=1,2,3,4$. Bottom: Pad\'e approximants $r_n^{[1/1]}(x)$ for $n=1, 2, 3, 4$, and zoom near the origin.}
\label{fig:approxs}
\end{center}
\end{figure}

The internal polynomial approximations are recursively defined as
\[
\begin{split}
& p_0(x)\equiv 1, \\
& p_{n+1}(x)=\frac{1}{2}\big(2p_n(x)-p_n(x)^2+x^2\big),\quad  n=0, 1, 2,\dots
\end{split}
\]
The following straightforward properties hold:
\begin{enumerate}
\item $p_n(x)$ is even and $\deg(p_n)=2^n$.
\item $|x|<p_n(x)<1$ for $x\in (-1, 1)$, and $p_n(\pm 1)=1$.
\item $p_n'(1)=1$ and $p_n'(-1)=-1$.
\item $\min_{-1\leq x\leq 1}p_n(x) = p_n(0)>0$.
\item The sequence $\{p_n(x)\}_{n\in\mathbb{N}}$ converges uniformly to $|x|$.
\end{enumerate}
Regarding efficiency issues, we have found that a direct implementation of $p_n(x)$ is a better option than using its recursive form. In Table \ref{table:internal_coeffs} are given the coefficients of the polynomials
\[
p_n(x) = \alpha_0x^{2^n}+\alpha_1x^{2^{n-1}}+\cdots\alpha_{n-1}x^2+\alpha_n
\]
for $n=1, 2, 3, 4$. Figure \ref{fig:approxs} (top left) shows the polynomials $p_n(x)$, for $n=1, 2, 3, 4$.

\begin{table}[th]
\caption{Coefficients of the internal approximations $p_n(x)$, for $n=1, 2, 3, 4$.}
\begin{center}
\begin{tabular}{lccccccccc}
$n$ & $\alpha_0$ & $\alpha_1$ & $\alpha_2$ & $\alpha_3$ & $\alpha_4$ & $\alpha_5$ & $\alpha_6$ & $\alpha_7$ & $\alpha_8$ \\
\hline
1 & $\frac{1}{2}$ & $\frac{1}{2}$ \\
2 & $-\frac{1}{8}$ & $\frac{3}{4}$ & $\frac{3}{8}$ \\
3 & $-\frac{1}{128}$ & $\frac{3}{32}$ & $-\frac{23}{64}$ & $\frac{31}{32}$ & $\frac{39}{128}$ \\
4 & $-\frac{1}{32768}$ & $\frac{3}{4096}$ & $-\frac{59}{8192}$ & $\frac{169}{4096}$ & $-\frac{2635}{16384}$ & $\frac{1693}{4096}$ & $-\frac{5891}{8192}$ & $\frac{4807}{4096}$ & $\frac{8463}{32768}$ \\
\hline
\end{tabular}
\end{center}
\label{table:internal_coeffs}
\end{table}

On the other hand, for given integers $k, m\geq 0$, we consider the polynomials
\[
P_{km}(\xi)=\sum_{n=0}^{k}\frac{(\tfrac{1}{2})_n(\tfrac{1}{2}-m)_m(n-k-m)_m}{n!(-k-m)_m(n+\tfrac{1}{2}-m)_m}\xi^n
\]
and
\[
Q_{km}(\xi)=\sum_{n=0}^{m}\frac{(-m)_n(-\tfrac{1}{2}-k)_n}{n!(-k-m)_n}\xi^n,
\]
where we have used the notation $(a)_n=a(a+1)\ldots(a+n-1)$ and $(a)_0=1$. Then, we follow \cite{KL1991} to build a sequence of Pad\'e approximations of the form
\[
\begin{split}
& x_0 = 1, \\
& x_{n+1} = x_n\frac{Q_{km}(1-x^2/x_n^2)}{P_{km}(1-x^2/x_n^2)}, \quad  n=0, 1, 2,\dots
\end{split}
\]
which converges to $|x|$, for any given $x\in [-1, 1]$. Thus, the Pad\'e approximants of order $[m/k]$ are defined recursively as
\[
\begin{split}
& r_0^{[m/k]}(x)\equiv 1, \\
& r_{n+1}^{[m/k]}(x)=r_n(x) \frac{Q_{km}(1-x^2/r_n(x)^2)}{P_{km}(1-x^2/r_n(x)^2))}, \quad  n=0, 1, 2,\dots
\end{split}
\]
(the superscript $[m/k]$ will be dropped unless necessary). It is not difficult to prove that each $r_n(x)$ verifies the properties 1-5 listed before for $p_n(x)$. Figure \ref{fig:approxs} (top right) shows the functions $r_1^{[m/k]}(x)$ for $m=k=1, 2, 3, 4$.

For the ease of coding, we give in Table \ref{table:pade_coeffs} the coefficients of the rational function $r_{n+1}(x)$ written in the form
\[
r_n(x)\frac{\alpha_0r_n(x)^{2m}+\alpha_1r_n(x)^{2(m-1)}x^2+\cdots+\alpha_{m-1}r_n(x)^2x^{2(m-1)}+\alpha_mx^{2m}}{\beta_0r_n(x)^{2k}+\beta_1r_n(x)^{2(k-1)}x^2+\cdots+\beta_{k-1}r_n(x)^2x^{2(k-1)}+\beta_kx^{2k}},
\]
%\[
%r_{n+1}(x) = r_n(x)\frac{\displaystyle \sum_{i=0}^m %\alpha_ir_n(x)^{2(m-i)}x^{2i}}{\displaystyle \sum_{i=0}^k \beta_ir_n(x)^{2(k-i)}x^{2i}},
%\]
for $0\leq k, m\leq 3$. For example, the family of Pad\'e approximants of order $[1/1]$ would be given by
\[
r_{n+1}(x) = r_n(x)\frac{r_n(x)^2+3x^2}{3r_n(x)^2+x^2},  \quad  n=0, 1, 2,\dots
\]
while those of order $[2/2]$ would be
\[
r_{n+1}(x) = r_n(x)\frac{r_n(x)^4+10r_n(x)^2x^2+5x^4}{5r_n(x)^4+10r_n(x)^2x^2+x^4},  \quad  n=0, 1, 2,\dots
\]
Figure \ref{fig:approxs} (bottom) depicts the approximants of order $[1/1]$ for $n=1, 2, 3, 4$, together with a zoom near the origin to check the precision of the approximations.

\begin{table}[th]
\caption{Coefficients of the rational function $r_1(x)$, for $0\leq k, m\leq 3$.}
\begin{center}
\begin{tabular}{c|cccc}
    & $k=0$ & $k=1$ & $k=2$ & $k=3$ \\ \hline
$m=0$ & $\frac{[1]}{[1]}$ & $\frac{[-2]}{[-3,1]}$ & $\frac{[8]}{[15,-10,3]}$ & $\frac{[-16]}{[-35,35,-21,-5]}$ \\  
$m=1$ & $\frac{[1,1]}{[2]}$  &  $\frac{[1,3]}{[3,1]}$ & $\frac{[4,20]}{[15,10,-1]}$ & $\frac{[8,56]}{[35,35,-7,1]}$  \\ 
$m=2$ & $\frac{[3,6,-1]}{[8]}$ &  $\frac{[1,6,1]}{[4,4]}$ & $\frac{[1,10,5]}{[5,10,1]}$ & $\frac{[6,84,70]}{[35,105,,21,-1]}$  \\
$m=3$ & $\frac{[5,15,-5,1]}{[16]}$
   & $\frac{[5,45,15,-1]}{[24,40]}$ & $\frac{[1,15,15,1]}{[6,20,6]}$ & $ \frac{[1,21,35,7]}{[7,35,21,1]}$ \\
\end{tabular}
\end{center}
\label{table:pade_coeffs}
\end{table}

Finally, we would like to remark that an additional advantage of the considered approximations $p_n(x)$ and $r_n(x)$ is that they provide an automatic entropy fix to handle sonic flow. The reason for this is that neither $p_n(x)$ nor $r_n(x)$ crosses the origin: see \cite[Sect. 4]{RMHD2017}.

\begin{remark}
	After extensive numerical investigation, our conclusion is that for relatively simple problems it is sufficient to consider a simple AVM solver, as for example HLL or the one based in the internal polynomial $p_1(x)$. However, when complex structures appear, the choice of more precise AVM solvers is a determinant factor in terms of efficiency (see \cite{RVM2014}). This is even true when AVM solvers are used as building blocks in the design of high-order schemes, as it was demonstrated in \cite{RVM_Osher2016}. \qed
\end{remark}

%%%%%

\section{Multidimensional AVM solvers} \label{sec:AVM2d}

The purpose of this section is to extend Balsara's multidimensional solver by changing, in an appropriate way, the underlying HLL flux by a more precise AVM solver. The key point to achieve this is the reinterpretation of Balsara's solver in such a way that the two-dimensional corrections are expressed as suitable combinations of one-dimensional numerical fluxes of HLL type, which in turn can be changed by another AVM flux. We will use the same notations as in Section \ref{sec:Balsara_solver} for the speeds, states and fluxes. In particular, we are interested in a simple four-wave model, so only the speeds $S_L$, $S_R$, $S_D$ and $S_U$ will be taken into account.

To begin with, consider the local stencil shown in Figure \ref{fig:fig2}, and let us see how the two-dimensional resolved flux $F^*$ is defined. Let $U_L^*$ and $U_R^*$ be the HLL states given by \eqref{eq:LRstates}, and let $F_L^*$ and $F_R^*$ be the transverse fluxes considered in Section \ref{sec:Balsara_solver}. Notice that the explicit form of these transverse fluxes is problem-dependent: see Section \ref{sec:MHD} for the case of the MHD equations. Next, consider the AVM-type flux \eqref{eq:AVMnumflux} in the $x$-direction given by
\begin{equation} \label{eq:AVM2dfluxF}
F^* = \frac{F_L^*+F_R^*}{2}-\frac{1}{2}Q_x^*(U_R^*-U_L^*),
\end{equation}
with $Q_x^*=|\lambda_\text{max}|f(|\lambda_\text{max}|^{-1}A^*)$, where $f(x)$ is a function verifying the stability condition \eqref{eq:stability_cond}, $A^*$ is a Roe-like matrix such that
\begin{equation} \label{eq:pseudoRoe}
A^*(U_R^*-U_L^*) = F_R^*-F_L^*,
\end{equation}
and $\lambda_\text{max}$ is a bound on its spectral radius.

\begin{remark}
	As stated in Section \ref{sec:Balsara_solver}, in general $F_\alpha^*\neq F(U_\alpha^*)$, $\alpha=L, R$, so $A^*$ is not exactly a Roe matrix. Its form will depend on how the transverse fluxes are built. In the cases in which such a matrix cannot be constructed, $A^*$ can be taken as the Jacobian of the flux $F$ evaluated at some average state, so \eqref{eq:pseudoRoe} holds only approximately. \qed
\end{remark}

\begin{remark} \label{rem:transverse_fluxes}
	Notice that the vertical contributions in the fluxes \eqref{eq:AVM2dfluxF} are included through the transverse fluxes $F_R^*\equiv F_R^*(U_R^*, G_R^*)$ and $F_L^*\equiv F_L^*(U_L^*, G_L^*)$, as commented in Section \ref{sec:Balsara_solver}.\qed
\end{remark}

It was shown in \cite{PVM2012} that the HLL flux can be viewed as an AVM solver based on the first-order polynomial
\begin{equation} \label{eq:HLLpoly}
f(x) = \alpha_0+\alpha_1 x,
\end{equation}
with
\[
\alpha_0 = \frac{S_R|S_L|-S_L|S_R|}{S_R-S_L}, \quad \alpha_1 = \frac{|S_R|-|S_L|}{S_R-S_L}.
\]
Thus, substituting in \eqref{eq:AVM2dfluxF} and using \eqref{eq:pseudoRoe}, we get the following expression:
\begin{equation} \label{eq:AVM2dflux}
\begin{split}
F^* = & \frac{S_U}{S_U-S_D}\bigg( \frac{F_L^*+F_R^*}{2}-\frac{1}{2}\big(\alpha_0(U_{RU}-U_{LU})+\alpha_1(F_R^*-F_L^*)\big)\bigg) \\
& -\frac{S_D}{S_U-S_D}\bigg( \frac{F_L^*+F_R^*}{2}-\frac{1}{2}\big(\alpha_0(U_{RD}-U_{LD})+\alpha_1(F_R^*-F_L^*)\big)\bigg) \\
& + \frac{\alpha_0}{2(S_U-S_D)}(G_{RU}-G_{LU}+G_{LD}-G_{RD}),
\end{split}
\end{equation}
where the explicit form \eqref{eq:LRstates} for the HLL states has been applied. Notice that the first two lines of \eqref{eq:AVM2dfluxF} can be interpreted as two one-dimensional HLL fluxes associated, respectively, to the upper and lower parts of the stencil, while the last line accounts for the cross derivative of the flux $G$.

Assuming now that we are in the subsonic case in both directions, so
\[
\alpha_0 = -2\frac{S_LS_R}{S_R-S_L}, \quad \alpha_1 = \frac{S_R+S_L}{S_R-S_L},
\]
substitution in \eqref{eq:AVM2dflux} leads to exactly the same formula \eqref{eq:2dfluxF} for Balsara's resolved flux. This means that the two-dimensional resolved flux  \eqref{eq:2dfluxF} can be reinterpreted as an AVM-type flux of the form \eqref{eq:AVM2dfluxF}. Analogous conclusions are obtained when the flux is supersonic in one or both directions.

Similar computations can be performed to express the resolved flux \eqref{eq:2dfluxG} in the $y$-direction as
\begin{equation} \label{eq:AVM2dfluxG}
G^* = \frac{G_D^*+G_U^*}{2}-\frac{1}{2}Q_y^*(U_U^*-U_D^*),
\end{equation}
where $Q_y^*=g(B^*)$; here, $g(x)=\beta_0+\beta_1x$ with 
\[
\beta_0 = \frac{S_U|S_D|-S_D|S_U|}{S_U-S_D}, \quad \beta_1 = \frac{|S_U|-|S_D|}{S_U-S_D},
\]
and $B^*$ should verify $B^*(U_U^*-U_D^*)=G_U^*-G_D^*$.

Once the two-dimensional fluxes \eqref{eq:2dfluxF} and \eqref{eq:2dfluxG} have been expressed as the AVM-type fluxes \eqref{eq:AVM2dfluxF} and \eqref{eq:AVM2dfluxG}, it is possible to extend these to build more precise two-dimensional solvers. This is achieved by changing the underlying HLL polynomial \eqref{eq:HLLpoly} by another suitable function satisfying \eqref{eq:stability_cond}, thus obtaining an AVM-type two-dimensional flux. It is important to remark that to build an AVM solver, only the maximum wave speeds $S_\alpha$ are needed; despite simplicity, this is another reason why we consider a four-wave model. In particular, we always assume that the wave structure of the solution is that shown in Figure \ref{fig:fig2} (right). On the other hand, notice that in \eqref{eq:AVM2dfluxF}-\eqref{eq:AVM2dfluxG} the resolved states $U_\alpha^*$ are precisely the HLL states \eqref{eq:UDstates}-\eqref{eq:LRstates}. We will conserve these resolved states even though the basis function $f(x)$ used to build the fluxes \eqref{eq:AVM2dfluxF}-\eqref{eq:AVM2dfluxG} is changed. A possible way to introduce more structure in the resolved states, in the case of a PVM solver, could be through the equivalence of PVM solvers and simple Riemann solvers stated in \cite{MCP2014}; this idea is currently under investigation.

An important consequence of the above considerations is that our two-dimensional AVM solvers could be adapted to define path-conservative schemes for nonconservative systems, in the spirit of \cite{Pares2006}. In particular, this allows to extend the two-dimensional AVM solvers to the case of hyperbolic systems with coupling and source terms. This will be the subject of a forthcoming paper, so we will not go into details here.

For the sake of completeness, we end this section by summarizing the construction of the horizontal numerical flux $F_{i+1/2,j}$ for a 2D AVM method. The vertical flux $G_{i,j+1/2}$ is defined in a totally analogous way.

First of all, for generic cell values $U_0$ and $U_1$, and maximal speeds of propagation $S_0$ and $S_1$, define the HLL states
\[
U^\text{HLL}_x(U_0,U_1,S_0,S_1)=\frac{S_1U_1-S_0U_0+F(U_0)-F(U_1)}{S_1-S_0}
\]
and
\[
U^\text{HLL}_y(U_0,U_1,S_0,S_1)=\frac{S_1U_1-S_0U_0+G(U_0)-G(U_1)}{S_1-S_0}.
\]
Given a function $f(x)$ satisfying \eqref{eq:stability_cond}, define the one-dimensional AVM fluxes
\[
F^\text{AVM}(U_0,U_1)=\frac{F(U_0)+F(U_1)}{2}-\frac{1}{2}Q_x(U_1-U_0)
\]
and
\[
G^\text{AVM}(U_0,U_1)=\frac{G(U_0)+G(U_1)}{2}-\frac{1}{2}Q_y(U_1-U_0),
\]
where the viscosity matrices $Q_x$ and $Q_y$ are defined as in \eqref{eq:viscosity_matrix}, in terms of suitable Roe or Jacobian matrices $A_x$ and $A_y$, respectively.

Consider now the six-cell stencil around the vertical edge depicted in Figure \ref{fig:fig1}. Then:
\begin{itemize}
	\item The 1d flux through the edge is defined by
	\[
	F_{i+1/2,j}^*=F^\text{AVM}(U_{ij},U_{i+1,j}).
	\]
	\item To build the 2d flux $F_{i+1/2,j+1/2}^*$ at the vertex $(x_{i+1/2}, y_{j+1/2})$, first compute the four local speeds of propagation $S_L$, $S_R$, $S_U$ and $S_D$ (we omit the indexes for clarity) as explained in Section \ref{sec:Balsara_solver}. 
	
	\underline{Subsonic case in both directions}: $S_L<0<S_R$ and $S_D<0<S_U$. Define the states
	\begin{equation} \label{eq:star_states}
	\begin{split}
	& U_L^*=U^\text{HLL}_y(U_{ij},U_{i,j+1},S_D,S_U), \\ & U_R^*=U^\text{HLL}_y(U_{i+1,j},U_{i+1,j+1},S_D,S_U),
	\end{split}
	\end{equation}
	and the 1d AVM fluxes
	\begin{equation} \label{eq:star_fluxes}
	\begin{split}
	& G_L^*=G^\text{AVM}(U_{ij},U_{i,j+1}), \\
	& G_R^*=G^\text{AVM}(U_{i+1,j},U_{i+1,j+1}).
	\end{split}
	\end{equation}
	Then build the transverse flux $F_L^*$ in terms of $U_L^*$ and $G_L^*$, and the transverse flux $F_R^*$ in terms of $U_R^*$ and $G_R^*$. Finally, the 2d flux is given by
	\[
	F_{i+1/2,j+1/2}^*=\frac{F_L^*+F_R^*}{2}-\frac{1}{2}Q_x^*(U_R^*-U_L^*),
	\]
	where $Q_x^*$ is defined as explained in Section \ref{sec:AVM2d} in terms of a Roe-like matrix $A_x^*$ satisfying \eqref{eq:pseudoRoe}; if such a matrix is not available, simply take the Jacobian of the physical flux $F$ evaluated at some average state.
	
	\underline{Supersonic case in some or both directions}: The flux $F_{i+1/2,j+1/2}^*$ is defined as the corresponding upwinded flux, exactly as in \cite[Sec. 2.1]{Balsara2012}.
	\item The 2d flux $F_{i+1/2,j-1/2}^*$ at the lower vertex is defined similarly.
	\item Finally, the 1d and 2d fluxes are assembled as
	\begin{equation} \label{eq:assembling}
	F_{i+1/2,j+1/2}=\frac{1}{6}F_{i+1/2,j+1/2}^*+\frac{2}{3}F_{i+1/2,j}^*+\frac{1}{6}F_{i+1/2,j-1/2}^*,
	\end{equation}
	where, for simplicity, Simpson's rule has been considered (see the end of Section \ref{sec:Balsara_solver}).
\end{itemize}

In short, the two-dimensional flux $F^*$ can be interpreted as follows. First, from the right states $U_{RD}$ and $U_{RU}$ in Figure \ref{fig:fig2}, we build the associated state $U_R^*$ and the AVM vertical flux $G_R^*$; from these we define the transverse flux $F_R^*$, which may be viewed as a flux in the $x$-direction involving the states $U_{RD}$ and $U_{RU}$ and the contributions $G_{RD}$ and $G_{RU}$ in the $y$-direction. Next, the flux $F_L^*$ is similarly defined from the left states $U_{LD}$ and $U_{LU}$. Once $U_L^*$, $U_R^*$, $F_L^*$ and $F_R^*$ have been constructed, we use them as building blocks to define $F^*$ as an AVM-type flux: $F^*\equiv F^\text{AVM}(U_L^*,U_R^*)$, where the fluxes $F(U_L^*)$ and $F(U_R^*)$ are substituted by $F_L^*$ and $F_R^*$ respectively. The AVM solver used is determined by the choice of the viscosity matrix $Q_x^*$.

Notice that the two-dimensional fluxes at vertices have been constructed as combination of several one-dimensional fluxes. Thus, for practical implementation only the definition of a one-dimensional AVM flux is really needed, what makes the computer coding simple and clear.

Finally, the choice of the time step is given by the CFL condition
\begin{equation} \label{eq:CFL}
\Delta t = 2 \delta \cdot \min_{ij} \Delta t_{ij}, \qquad  \delta\in (0, 1],
\end{equation}
where each local time step $\Delta t_{ij}$ is defined as
\[
\Delta t_{ij} = \bigg( \frac{\lambda_{ij}^x}{\Delta x}+\frac{\lambda_{ij}^y}{\Delta y} \bigg)^{-1},
\]
where $\lambda_{ij}^\alpha$ denotes the maximal speed of propagation in the cell $C_{ij}$ in the $\alpha$-direction, taking into account both the speeds at the edges and those arising in the two-dimensional Riemann problems at vertices. As remarked in Section \ref{sec:Balsara_solver} (see also \cite{Balsara2010}), the two-dimensional contributions allow the use of a maximal CFL number $\delta$ of unity. This represents an additional advantage with respect to standard methods in which a one-dimensional solver is applied dimension by dimension, for which the maximum CFL number is $0.5$.

%%%%%

\section{Application to ideal magnetohydrodynamics} \label{sec:MHD}

In this section we focus on applications to ideal Newtonian MHD, although we must remark that the two-dimensional AVM schemes introduced in the previous section are applicable to general hyperbolic systems. We will use the following notation: $\rho$ is density, $\vv=(v_x,v_y,v_z)^t$ is the velocity field, $\BB=(B_x,B_y,B_z)^t$ is the magnetic field, $E$ is the total energy and $P$ is the hydrostatic pressure. Then, the equations of MHD read as
\begin{equation} \label{eq:MHD}
\begin{cases}
\partial_t\rho +\nabla\cdot(\rho\vv)=0, \\
\partial_t(\rho\vv)+\nabla\cdot\big(\rho\vv\vv^t+(P+\frac{1}{2}\BB^2)I-\BB\BB^t\big)=0, \\
\partial_t\BB+\nabla\cdot(\BB\vv^t-\vv\BB^t)=0, \\
\partial_tE+\nabla\cdot\big((E+P+\frac{1}{2}\BB^2)\vv-\BB(\vv\cdot\BB)\big)=0,
\end{cases}
\end{equation}
together with the divergence-free condition
\begin{equation} \label{eq:divfree}
\nabla\cdot\BB=0.
\end{equation}
To close the system, an ideal equation of state $P=(\gamma-1)\rho\varepsilon$ is considered, where $\gamma$ is the adiabatic constant and $\varepsilon$ is the specific internal energy, which is related to the total energy by $E=\frac{1}{2}\rho\vv^2+\frac{1}{2}\BB^2+\rho\varepsilon$. The total pressure is defined as $P^*=P+\frac{1}{2}\BB^2$, where $\frac{1}{2}\BB^2$ represents the magnetic pressure.

On the other hand, define $\bb=\BB/\sqrt{\rho}$ and the acoustic sound speed as $a=\sqrt{\gamma P/\rho}$. Then the fastest speeds of propagation in the $\alpha$-direction, for $\alpha=x,y$, are given by
\begin{equation} \label{eq:speeds}
\lambda_\alpha^1 = v_\alpha-c_{\alpha,f}, \qquad \lambda_\alpha^8 = v_\alpha+c_{\alpha,f}, 
\end{equation}
where
\[
c_{\alpha,f}^2 = \frac{1}{2}\big(a^2+\bb^2+\sqrt{(a^2+\bb^2)^2-4a^2b_\alpha^2}\big).
\]
For a detailed description of the spectral structure of system \eqref{eq:MHD} see \cite{Brio-Wu,Serna}.

We remark that in the derivation of \eqref{eq:MHD} from physical principles, the condition \eqref{eq:divfree} is thoroughly used. However, it is also possible to derive a MHD model without imposing the divergence constraint in the process, as proposed in \cite{Powell1994} (see also \cite{FMR2009,PRLGZ1999}). This leads to the \emph{nonconservative form} of the MHD equations:
\begin{equation} \label{eq:NCMHD}
\begin{cases}
\partial_t\rho +\nabla\cdot(\rho\vv)=0, \\
\partial_t(\rho\vv)+\nabla\cdot\big(\rho\vv\vv^t+(P+\frac{1}{2}\BB^2)I-\BB\BB^t\big)=-(\nabla\cdot\BB)\BB, \\
\partial_t\BB+\nabla\cdot(\BB\vv^t-\vv\BB^t)=-(\nabla\cdot\BB)\vv, \\
\partial_tE+\nabla\cdot\big((E+P+\frac{1}{2}\BB^2)\vv-\BB(\vv\cdot\BB)\big)=-(\nabla\cdot\BB)(\vv\cdot\BB).
\end{cases}
\end{equation}
Notice that if condition \eqref{eq:divfree} is imposed, then \eqref{eq:NCMHD} reduces to \eqref{eq:MHD}. The spectral structure of problem \eqref{eq:NCMHD} was fully analized in \cite{Powell1994}; in particular, the maximal speeds of propagation are the same as in the conservative case.

In the design of numerical methods for solving \eqref{eq:MHD} it is fundamental to handle the divergence constraint \eqref{eq:divfree} in a proper way. In particular, in the presence of shocks, standard numerical methods can produce large divergence errors with may lead to negative densities or pressures. Several methods have been proposed in the literature to impose the divergence-free condition numerically: see \cite{Toth2000} for an in-depth review on this topic. 
A technique of particular interest to handle the divergence constraint was proposed by Powell in \cite{Powell1994}, which is based in the eight-waves model \eqref{eq:NCMHD}. Notice that taking divergence of the magnetic field equation in \eqref{eq:NCMHD} gives
\[
\partial_t(\nabla\cdot\BB)+\nabla\cdot((\nabla\cdot\BB)\vv)=0,
\]
so any possible non-zero divergence produced by the numerical scheme is advected by the velocity field and eventually leaves the domain through the boundaries. Thus, the nonconservative form \eqref{eq:NCMHD} of the MHD system include a mechanism to handle the divergence constraint; we will back to this point later.

From now on we will concentrate on the two-dimensional case. Considering the vector of conserved variables $U=(\rho,\rho\vv,\BB,E)^t$, system \eqref{eq:MHD} can be written in conservative form as
\[
\partial_tU+\partial_xF(U)+\partial_yG(U)=0,
\]
where the physical fluxes are given by
\[
F(U)=
\begin{pmatrix}
\rho v_x \\
\rho v_x^2+P^*-B_x^2 \\
\rho v_xv_y-B_xB_y \\
\rho v_xv_z-B_xB_z \\
0 \\
v_xB_y-v_yB_x \\
v_xB_z-v_zB_x \\
v_x(E+P^*) - B_x(\vv\cdot\BB)
\end{pmatrix}
,\quad 
G(U) = 
\begin{pmatrix}
\rho v_y \\
\rho v_xv_y-B_xB_y \\
\rho v_y^2+P^*-B_y^2 \\
\rho v_yv_z-B_yB_z \\
v_yB_x-v_xB_y \\
0 \\
v_yB_z-v_zB_y \\
v_y(E+P^*) - B_y(\vv\cdot\BB)
\end{pmatrix}
.
\]
The transverse fluxes (see Remark \ref{rem:transverse_fluxes}) can be defined using values from the associated intermediate state and the normal fluxes. Thus, denoting the components of the states \eqref{eq:star_states} as $U_\alpha^*=(u_{\alpha,i}^*)_{i=1}^8$ and the fluxes \eqref{eq:star_fluxes} as $G_\alpha^*=(g_{\alpha,i}^*)_{i=1}^8$, for $\alpha=L, R$, the transverse fluxes $F_\alpha^*$ are defined as
\[
F_\alpha^* = 
\begin{pmatrix}
u_{\alpha,2}^* \medskip \\
g_{\alpha,3}^*+\dfrac{(u_{\alpha,2}^*)^2-(u_{\alpha,3}^*)^2}{u_{\alpha,1}^*}+(u_{\alpha,6}^*)^2-(u_{\alpha,5}^*)^2 \medskip \\
g_{\alpha,2}^* \medskip \\
\dfrac{u_{\alpha,2}^*u_{\alpha,4}^*}{u_{\alpha,1}^*}-u_{\alpha,5}^*u_{\alpha,7}^* \medskip \\
0 \medskip \\
\dfrac{u_{\alpha,2}^*u_{\alpha,6}^*-u_{\alpha,3}^*u_{\alpha,5}^*}{u_{\alpha,1}^*} \medskip \\
\dfrac{u_{\alpha,2}^*u_{\alpha,7}^*-u_{\alpha,4}^*u_{\alpha,5}^*}{u_{\alpha,1}^*} \medskip \\
\dfrac{u_{\alpha,2}^*}{u_{\alpha,1}^*}(u_{\alpha,8}^*+g_{\alpha,3}^*-\dfrac{(u_{\alpha,3}^*)^2}{u_{\alpha,1}^*}+{u_{\alpha,6}^*}^2) - \dfrac{u_{\alpha,5}^*}{u_{\alpha,1}^*}(u_{\alpha,2}^*u_{\alpha,5}^*+u_{\alpha,3}^*u_{\alpha,6}^*+u_{\alpha,4}^*u_{\alpha,7}^*)
\end{pmatrix}
,
\]
taking into account the form of the fluxes $F$ and $G$. The transverse fluxes $G_\beta^*$, for $\beta=D, U$, can be constructed in an analogous way from $U_\beta^*$ and $F_\beta^*$. 

As commented before, some kind of divergence cleaning mechanism has to be added to the numerical scheme. For doing this we consider the nonconservative system \eqref{eq:NCMHD}, that can be expressed as
\begin{equation} \label{eq:ncsystem}
\partial_tU+\partial_xF(U)+\partial_yG(U)=-(\nabla\cdot\BB)S(U),
\end{equation}
where the source term is given by
\[
S(U)=(0, \BB, \vv, \vv\cdot\BB)^t.
\]
Alternatively, \eqref{eq:ncsystem} can be written in the form
\begin{equation} \label{eq:ncsystem_alt}
\partial_tU+\partial_xF(U)+\partial_yG(U)=-\Bcal_1(U)\partial_xU-\Bcal_2(U)\partial_yU,
\end{equation}
where $\Bcal_1(U)$ and $\Bcal_2(U)$ are the $8\times 8$ matrices
\[
\Bcal_1(U)=(\OO\vert\OO\vert\OO\vert\OO\vert S(U)\vert\OO\vert\OO\vert\OO\vert), \quad \Bcal_2(U)=(\OO\vert\OO\vert\OO\vert\OO\vert\OO\vert S(U)\vert\OO\vert\OO\vert),
\]
with $\OO$ the zero column vector in $\R^8$. Now, following \cite{MCP2014}, both the HLL states and the numerical fluxes must be properly modified in order to take into account the extra terms $\Bcal_1(U)$ and $\Bcal_2(U)$ appearing in \eqref{eq:ncsystem_alt}. In particular, the HLL states should be modified as
\[
\widetilde{U}^\text{HLL}_x(U_0,U_1,S_0,S_1)=\frac{S_1U_1-S_0U_0+F(U_0)-F(U_1)-\Bcal_1(\widetilde{U})(U_1-U_0)}{S_1-S_0}
\]
and
\[
\widetilde{U}^\text{HLL}_y(U_0,U_1,S_0,S_1)=\frac{S_1U_1-S_0U_0+G(U_0)-G(U_1)-\Bcal_2(\widetilde{U})(U_1-U_0)}{S_1-S_0},
\]
being $\widetilde{U}$ some average state. Accordingly, the numerical fluxes now should read as
\[
\widetilde{F}^\text{AVM}(U_0,U_1)=\frac{F(U_0)+F(U_1)}{2}-\frac{1}{2}(Q_x-\Bcal_1(\widetilde{U}))(U_1-U_0)
\]
and
\[
\widetilde{G}^\text{AVM}(U_0,U_1)=\frac{G(U_0)+G(U_1)}{2}-\frac{1}{2}(Q_y-\Bcal_2(\widetilde{U}))(U_1-U_0).
\]

\begin{remark}
	Notice that for MHD the vectors $\Bcal_i(\widetilde{U})(U_1-U_0)$ are simply
	\[
	\Bcal_1(\widetilde{U})=((B_x)_1-(B_x)_0)S(\widetilde{U}), \quad \Bcal_2(\widetilde{U})=((B_y)_1-(B_y)_0)S(\widetilde{U}). \tag*{\qed}
	\]
\end{remark}

Finally, the definition of a two-dimensional AVM flux for the nonconservative system \eqref{eq:ncsystem_alt} follows exactly the same guidelines explained at the end of Section \ref{sec:AVM2d}, substituting $U_x^\text{HLL}$, $U_y^\text{HLL}$, $F^\text{AVM}$ and $G^\text{AVM}$ by $\widetilde{U}_x^\text{HLL}$, $\widetilde{U}_y^\text{HLL}$, $\widetilde{F}^\text{AVM}$ and $\widetilde{G}^\text{AVM}$ respectively.

Notice that the proposed schemes are designed for solving the eight-waves model \eqref{eq:NCMHD}, so they automatically include a divergence cleaning procedure. Thus, if both the initial and boundary conditions satisfy the divergence constraint, it is expected that the scheme will keep it up to the accuracy of the truncation error. Notice also that the proposed nonconservative divergence cleaning technique could be applied, in a similar way, to any arbitrary numerical method.% In Section \ref{sec:numerical} we consider two versions of the two-dimensional AVM methods introduced so far: the first one, for solving the conservative MHD equations \eqref{eq:MHD} using the well-known projection method \cite{BB1980} for divergence cleaning, and the second one for solving the nonconservative MHD system \eqref{eq:NCMHD} as explained above. We will see that both kind of methods produce similar results, although the second type has a smaller computational cost, as the projection method requires the solution of a collection of elliptic problems at each time step.

%%%%%

\section{Numerical results} \label{sec:numerical}

The performances of the proposed multidimensional AVM schemes are analyzed in this section. To this end, we have chosen a number of tests that constitute standard references in MHD, in order that our schemes could be confronted with other methods in the literature.

An important issue in the construction of the schemes is the proper choice of the speeds of propagation. For this we have followed the guidelines stated in Section \ref{sec:Balsara_solver}, with the expressions \eqref{eq:speeds} for the maximal wave speeds. In particular, for the intermediate speeds $\bar{\lambda}_\alpha(U_0, U_1)$ appearing in \eqref{eq:SRU_SLU} we have simply chosen the speeds $\lambda_\alpha(\widetilde{U})$ at the intermediate state $\widetilde{U}=\frac{U_0+U_1}{2}$. Another possible but more expensive choice is to consider $\widetilde{U}$ as the Roe state associated to the states $U_0$ and $U_1$. In our computational experiments, both choices have led to similar results. Once the maximal speeds of propagation have been computed, the time step is determined from the CFL condition \eqref{eq:CFL}. In particular, we remark that for the 2D methods the theoretical admissible CFL number is $1$, while for the 1D methods is $0.5$.

We will denote as Int-$n$ and Pad\'e-$[m/k]$ the AVM methods based on the internal polynomial $p_n(x)$ and the Pad\'e approximation $r_1^{[m/k]}(x)$, respectively. For each AVM method, 2D denotes its multidimensional version, where the numerical flux at an edge is given by \eqref{eq:assembling}; on the other hand, 1D refers to the corresponding one-dimensional solver applied dimension by dimension, which is equivalent to taking $F_{i+1/2,j+1/2}=F_{i+1/2,j}^*$ in \eqref{eq:assembling}.

Finally, unless otherwise stated, the divergence cleaning technique introduced in Section \ref{sec:MHD} will be used in all the numerical experiments.

\subsection{Accuracy test} \label{test:accuracy}

To study the accuracy of the proposed schemes, we consider a test from \cite{LW2003} for the ideal Euler equations with exact smooth solution given by
\[
\rho(x, y, t) = 1+0.2\sin(\pi(x+y-t/2)), \quad v_x=1,\quad v_y=-1/2, \quad P=1,
\]
and adiabatic constant $\gamma=1.4$. The computational domain is $[-1, 1]\times [-1, 1]$, with periodic boundary conditions in both directions. The schemes have been run until a final time $T=4$, which corresponds to a full period of the wave.

\begin{figure}[!ht]
\begin{center}
\includegraphics[width=0.7\textwidth]{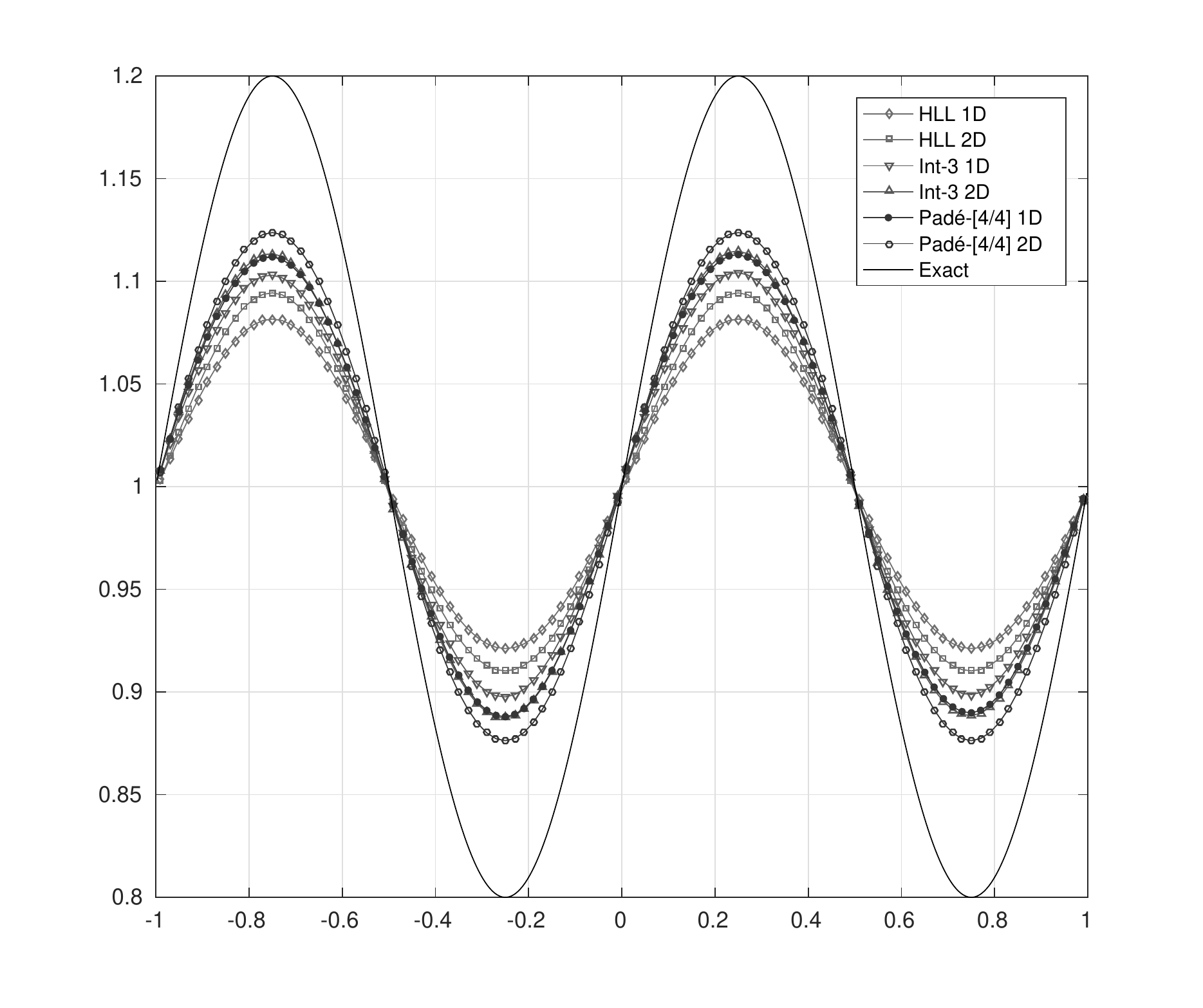}
\caption{Test \ref{test:accuracy}. Comparison of several 1D-2D AVM methods on a $100\times 100$ mesh. Diagonal cut along the main diagonal for the density variable.}
\label{fig:comparison}
\end{center}
\end{figure}

Figure \ref{fig:comparison} shows a cut along the main diagonal of the final densities obtained on a $100\times 100$ mesh with the 1D and 2D versions of HLL, Int-3 and Pad\'e-[4/4] schemes. Notice that the Int-3 method is based on an polynomial of eighth degree, while the Pad\'e-[4/4] method is based on a rational function of degree $[8/8]$. All the computations have been done using a common CFL number of $0.5$, although the 2D methods perform well with CFL numbers up to $0.8$. As it can be observed, for each scheme the two-dimensional contributions provide a more precise computation of the solution than their one-dimensional counterparts.

\begin{figure}[!ht]
\begin{center}
\includegraphics[width=0.7\textwidth]{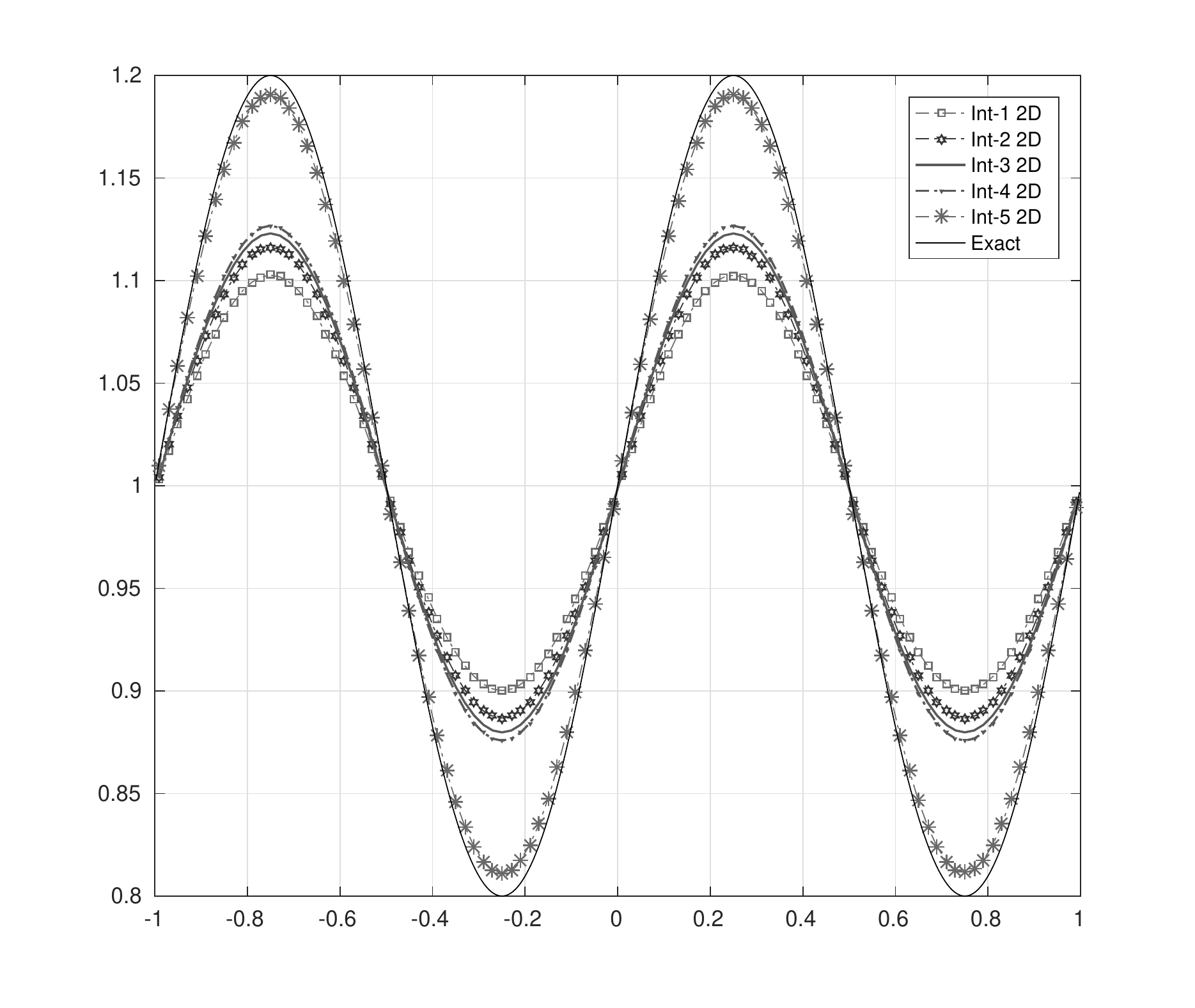}
\caption{Test \ref{test:accuracy}. Solutions obtained with the Int-$n$ 2D schemes, for $n=1, 2, 3, 4, 5$. Diagonal cut along the main diagonal for the density variable.}
\label{fig:comparison_Int_n}
\end{center}
\end{figure}

To check the behavior of the solutions in terms of the basis function chosen for an AVM 2D method, we have represented in Figure \ref{fig:comparison_Int_n} the results obtained with the Int-$n$ 2D schemes, for $n=1, 2, 3, 4, 5$. It is clear that the precision of the solutions increases as the degree of the basis polynomial $p_n(x)$ increases. A similar remark is also valid for the Pad\'e-$[m/k]$ 2D schemes.

\begin{figure}[!ht]
\begin{center}
\includegraphics[width=0.49\textwidth]{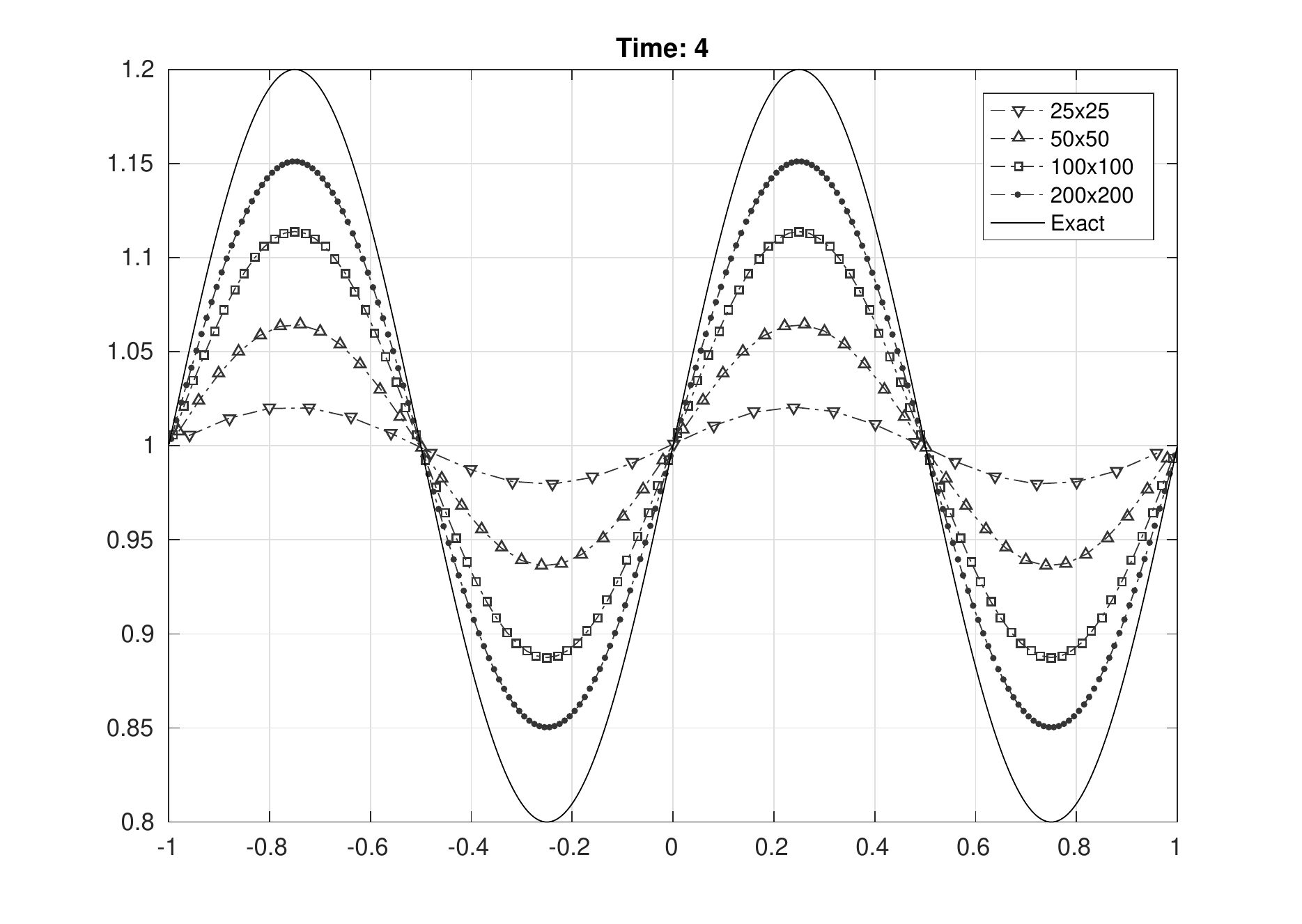}
\includegraphics[width=0.49\textwidth]{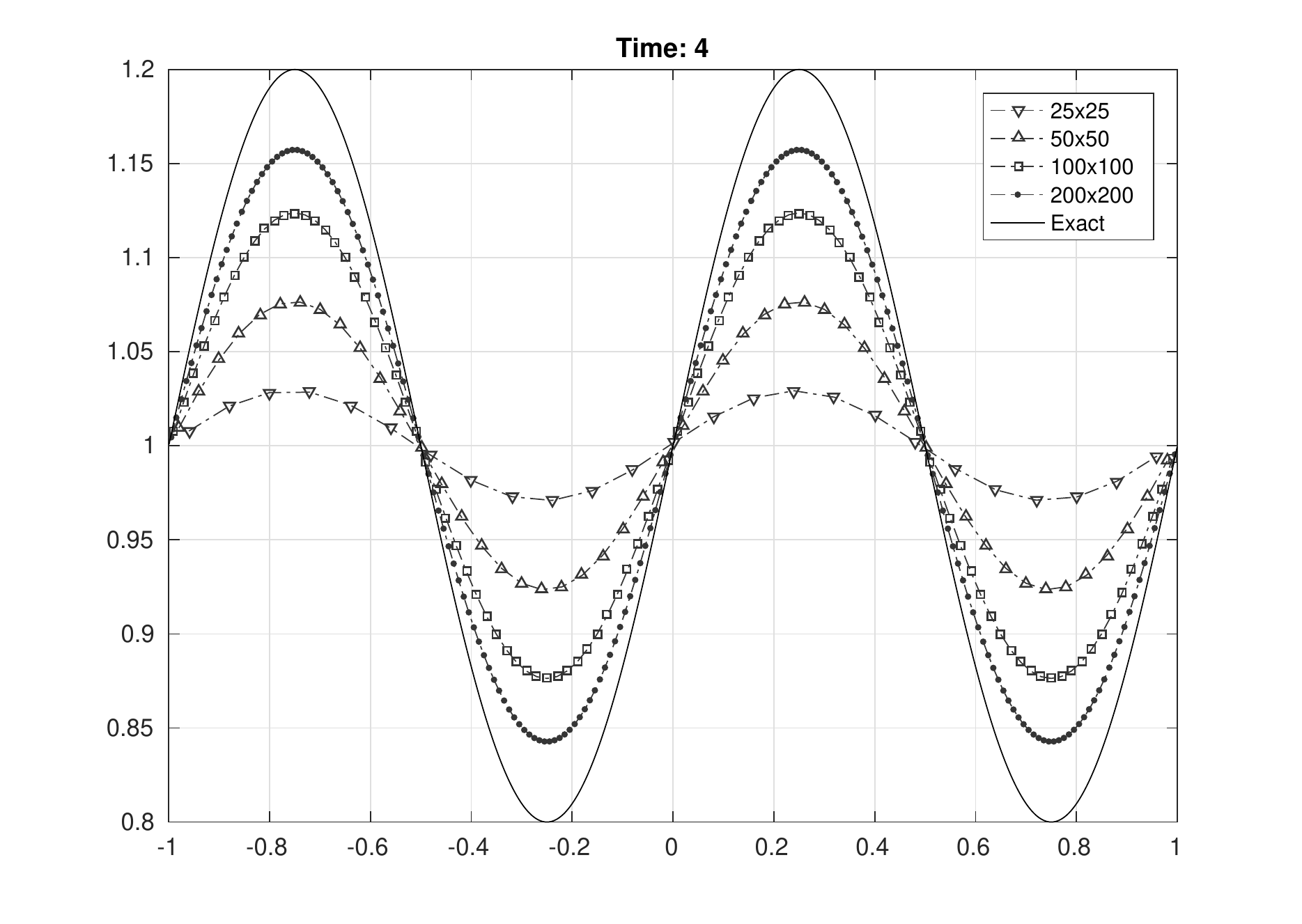}
\caption{Test \ref{test:accuracy}. Diagonal cut along the main diagonal for the density variable on several meshes. Left: Int-3 2D. Right: Pad\'e-[4/4] 2D.}
\label{fig:meshes}
\end{center}
\end{figure}

\begin{figure}[!ht]
\begin{center}
\includegraphics[width=0.7\textwidth]{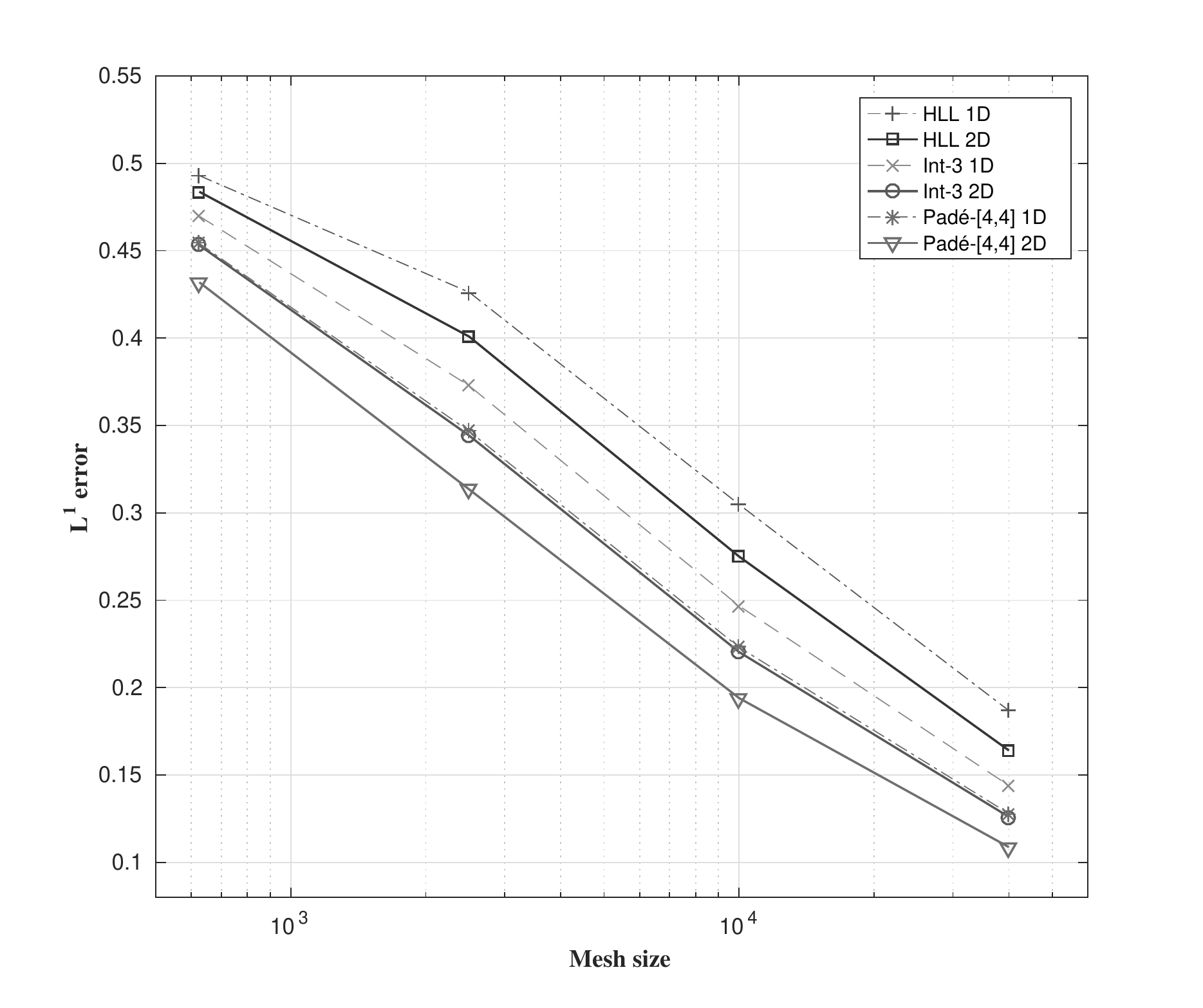}
\caption{Test \ref{test:accuracy}. Error curves for several 1D-2D AVM methods.}
\label{fig:errorcurves}
\end{center}
\end{figure}

On the other hand, to check the convergence to the exact solution as the mesh is refined, in Figure \ref{fig:meshes} are represented the solutions obtained with different meshes of size $N\times N$, for $N=25, 50, 100, 200$. Table \ref{tab:errors} shows the $L^1$ errors in the density variable, while the corresponding error curves are depicted in Figure \ref{fig:errorcurves}. Notice that the Int-3 2D scheme provides very similar result as the Pad\'e-[4/4] 1D method.

\begin{table}[!htb]
\centering
\begin{tabular}{lll|lll}
 Scheme & Mesh size & $L^1$ error &   Scheme & Mesh size & $L^1$ error\\ \hline
 HLL 1D & $25\times 25$   & 4.9299e-01 &  HLL 2D & $25\times 25$   & 4.8377e-01   \\
 			& $50\times 50$   & 4.2625e-01 &             & $50\times 50$   & 4.0087e-01  \\
 			& $100\times 100$ & 3.0486e-01 &    		 & $100\times 100$ & 2.7504e-01 \\ 
 			& $200\times 200$ & 1.8684e-01 &   			 & $200\times 200$ & 1.6402e-01\\ 
 			\hline
 Int-3 1D & $25\times 25$   & 4.6978e-01 &  Int-3 2D & $25\times 25$   & 4.5315e-01    \\
 			& $50\times 50$   & 3.7274e-01 &  			   & $50\times 50$   & 3.4422e-01   \\
 			& $100\times 100$ & 2.4677e-01 &    		   & $100\times 100$ & 2.2038e-01 \\ 
 			& $200\times 200$ & 1.4393e-01 &  			   & $200\times 200$ & 1.2603e-01   \\ 
 			\hline
 Pad\'e-[4/4] 1D & $25\times 25$   & 4.5418e-01 & Pad\'e-[4/4] 2D & $25\times 25$   & 4.3194e-01    \\
 			& $50\times 50$   & 3.4670e-01 &  			& $50\times 50$   & 3.1343e-01   \\
 			& $100\times 100$ & 2.2308e-01 &  			& $100\times 100$ & 1.9403e-01   \\ 
 			& $200\times 200$ & 1.2791e-01 &  			& $200\times 200$ & 1.0871e-01   \\ 
\end{tabular}
\caption{Comparison of the $L^1$ errors obtained with several 1D and 2D AVM schemes.}
\label{tab:errors}
\end{table}

\subsection{Orszag-Tang vortex} \label{test:OT}

The Orszag-Tang vortex \cite{Orszag-Tang} is a well-known model of complex MHD flow containing many significant features. Departing from an initial smooth state, the system develops a series of complex interactions between the different shock waves generated as the system evolves in the transition to turbulence.

Specifically, we consider the initial conditions proposed in \cite{FMR2009}:
\[
\begin{split}
& \rho(x, y, 0) = \gamma^2, \quad v_x(x, y, 0) = -\sin(y), \quad v_y(x, y, 0) = \sin(x), \\
& B_x(x, y, 0) = -\sin(y), \quad B_y(x, y, 0) = \sin(2x), \quad P(x, y, 0) = \gamma,
\end{split}
\]
with $\gamma=5/3$. The computational domain is given by $[0, 2\pi]\times [0, 2\pi]$, with periodic boundary conditions in both directions.

We use this problem to test the robustness of several AVM 2D schemes equipped with the divergence cleaning technique introduced in Section \ref{sec:MHD}. As it is well-known (see, e.g., \cite{FMR2009}), negative pressures may appear if the divergence is not properly controlled. In order to compare, we have considered the HLL 2D, Int-3 2D and Pad\'e-[4/4] 2D methods, with a common CFL number of $0.8$.

\begin{figure}[!ht]
\begin{center}
\includegraphics[height=0.24\textheight]{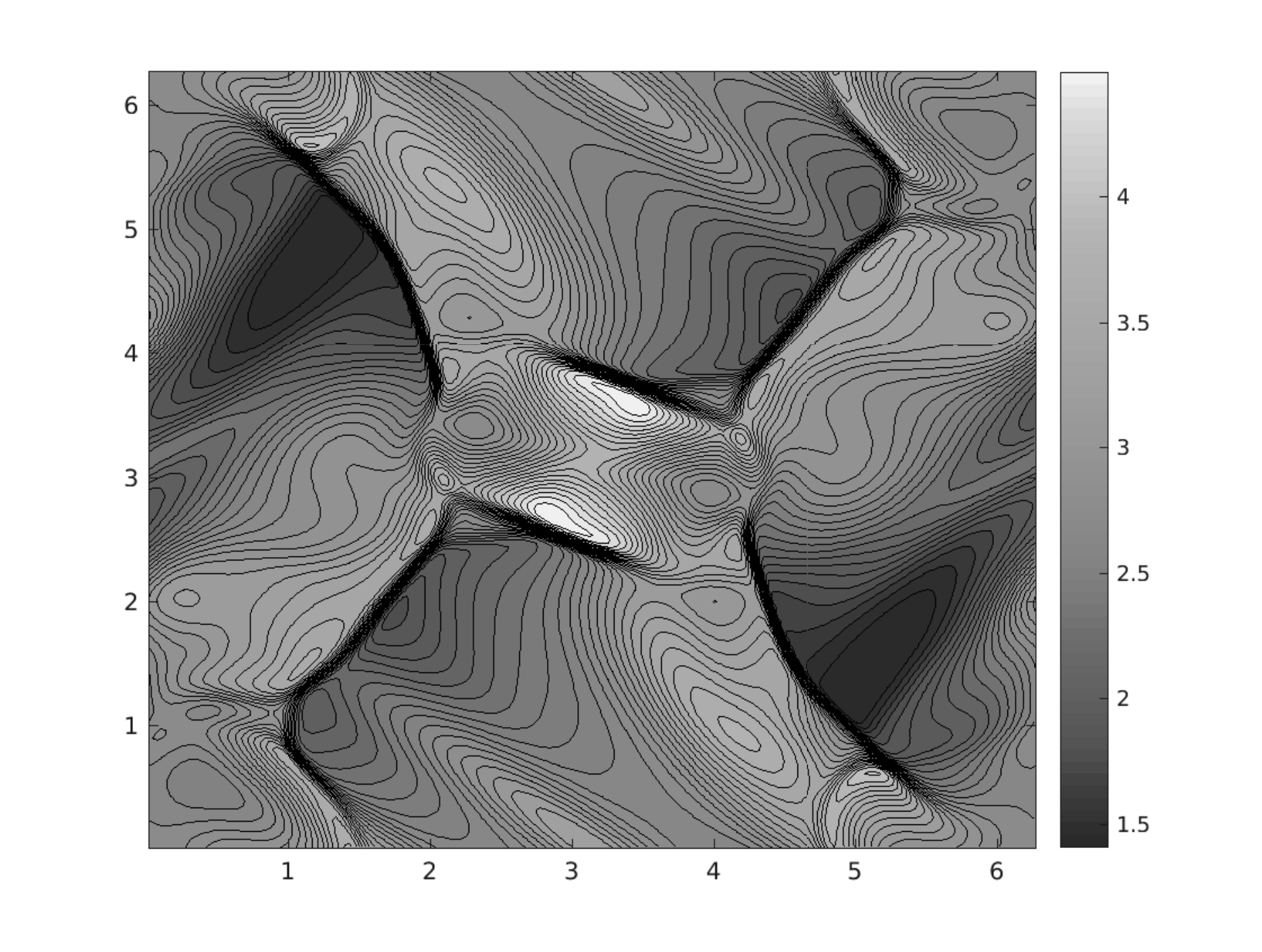}
\includegraphics[height=0.24\textheight]{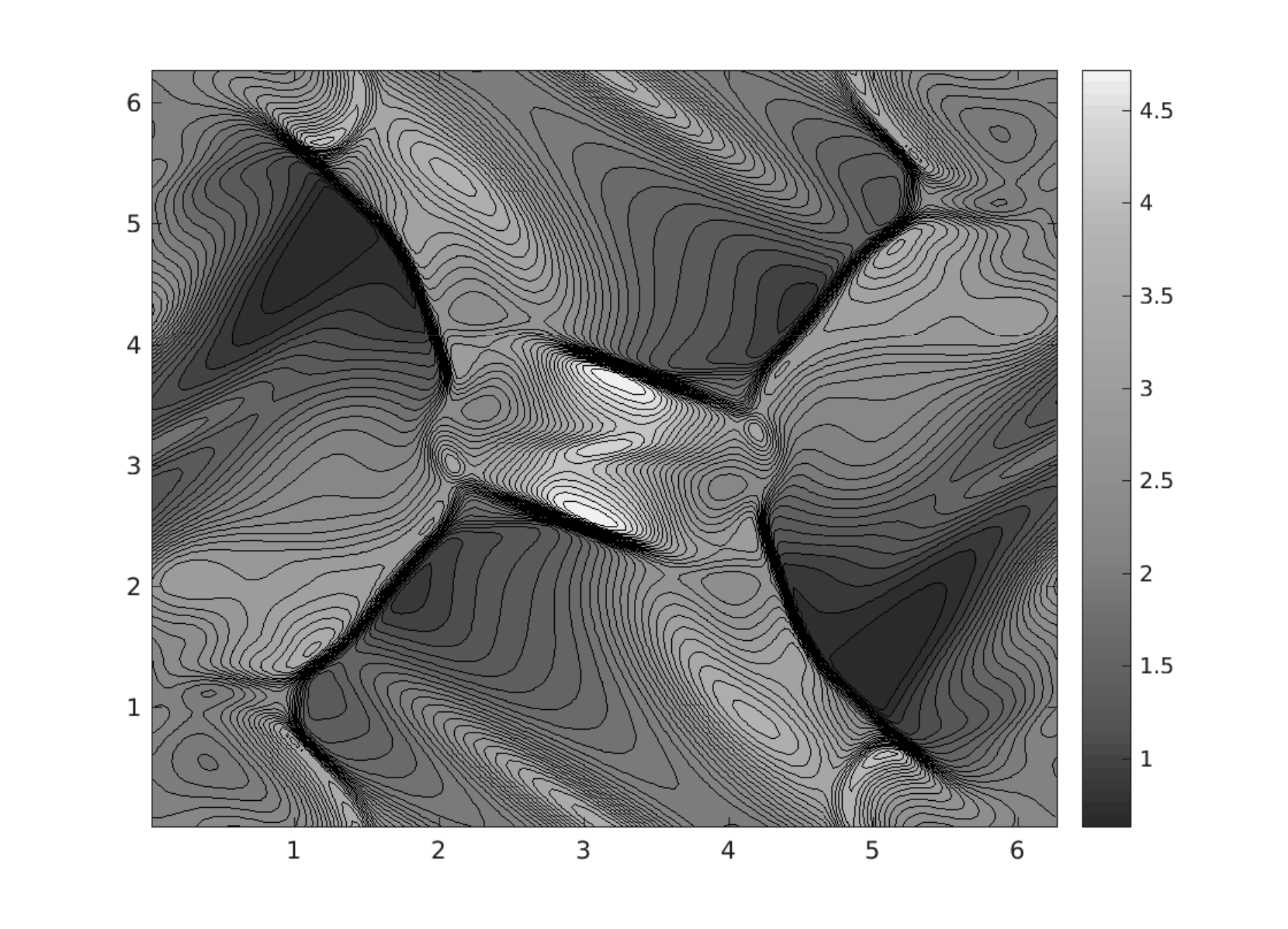} \\
\includegraphics[height=0.24\textheight]{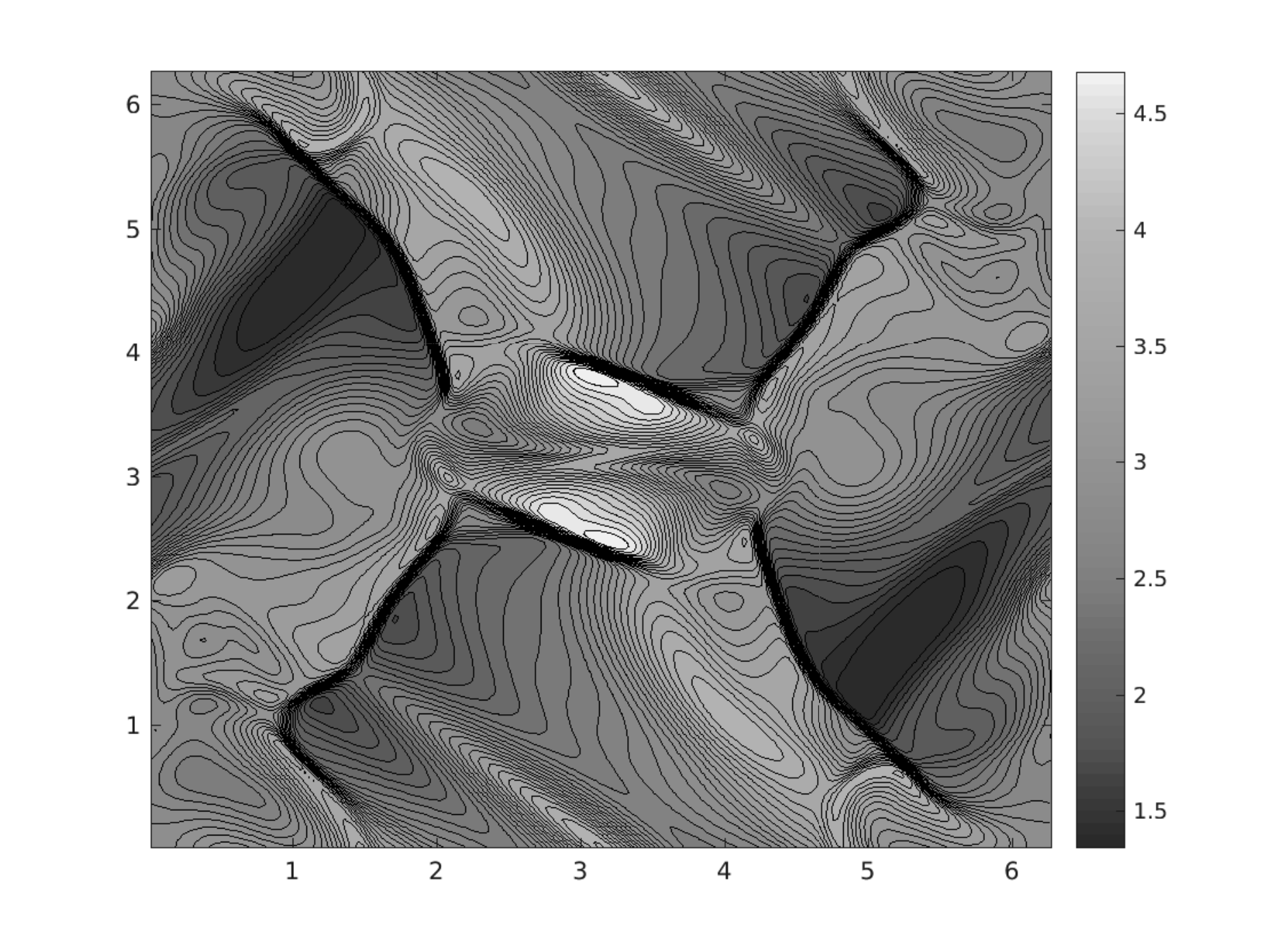}
\includegraphics[height=0.24\textheight]{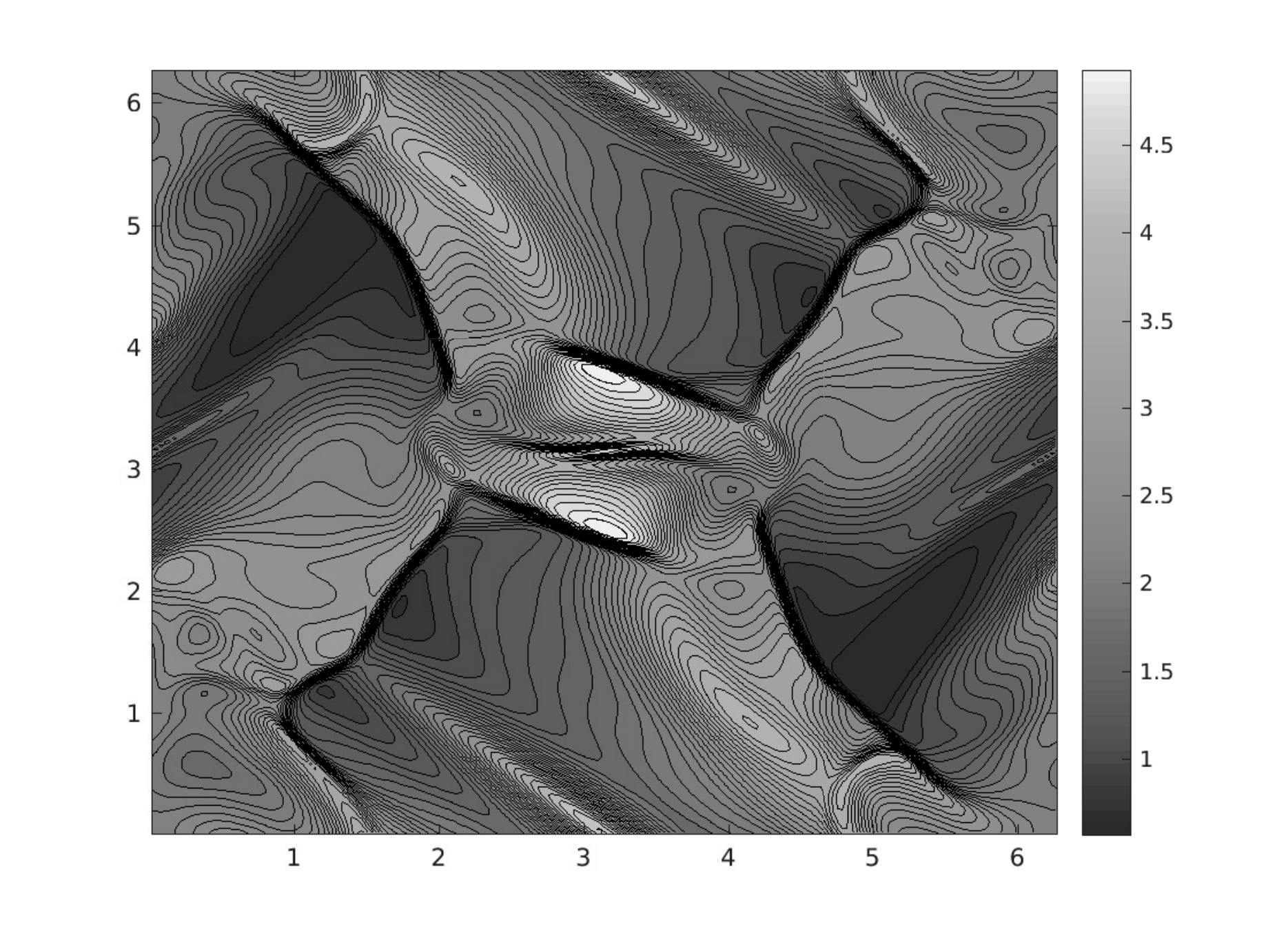} \\
\includegraphics[height=0.24\textheight]{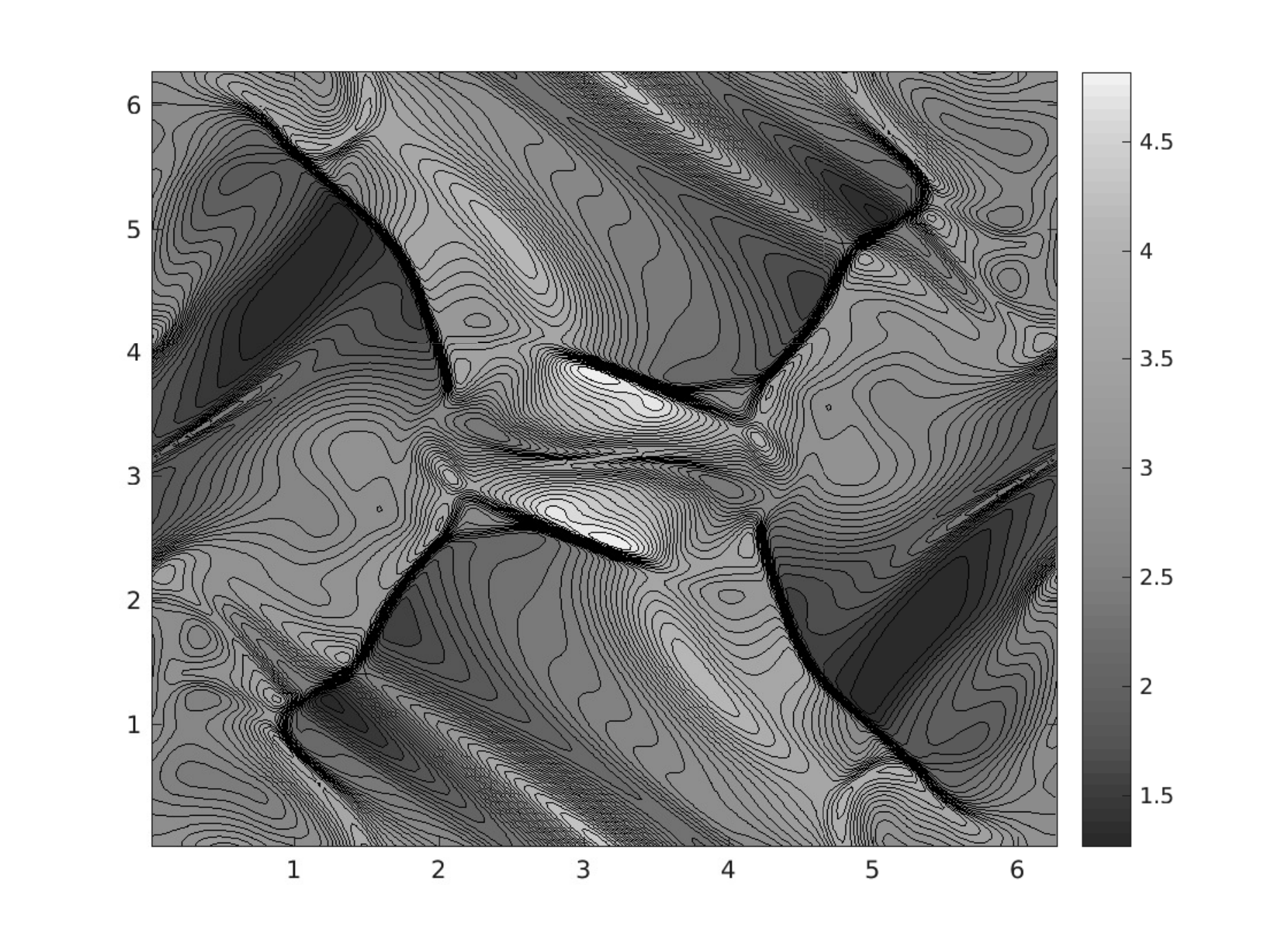}
\includegraphics[height=0.24\textheight]{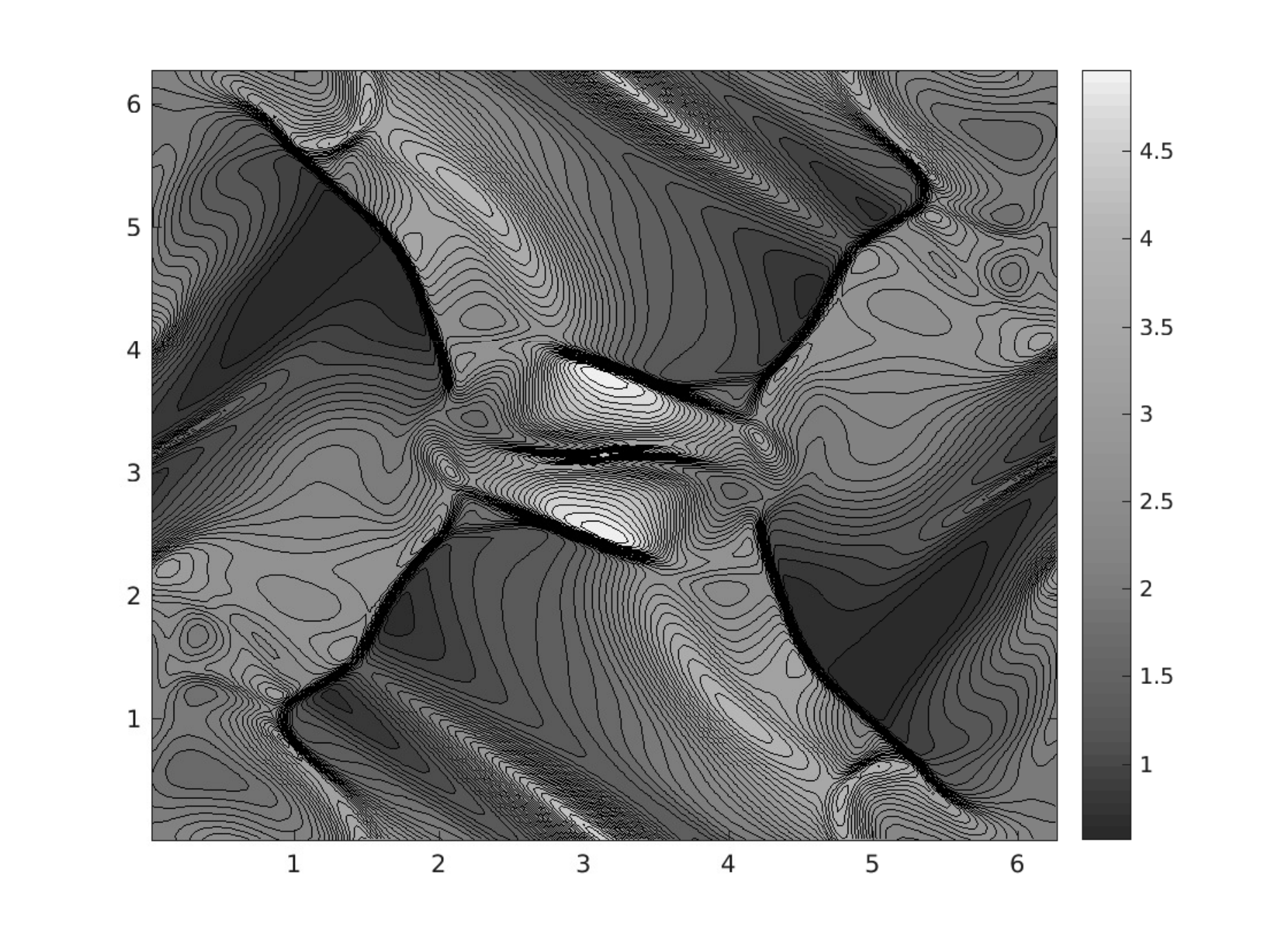}
\caption{Test \ref{test:OT}. Density (left) and pressure (right) computed at time $t=\pi$ on a $200\times 200$ mesh. From top to bottom: HLL 2D, Int-3 2D and Pad\'e-[4/4] 2D.}
\label{fig:OT}
\end{center}
\end{figure}

\begin{figure}[!ht]
\begin{center}
\includegraphics[width=0.48\textwidth]{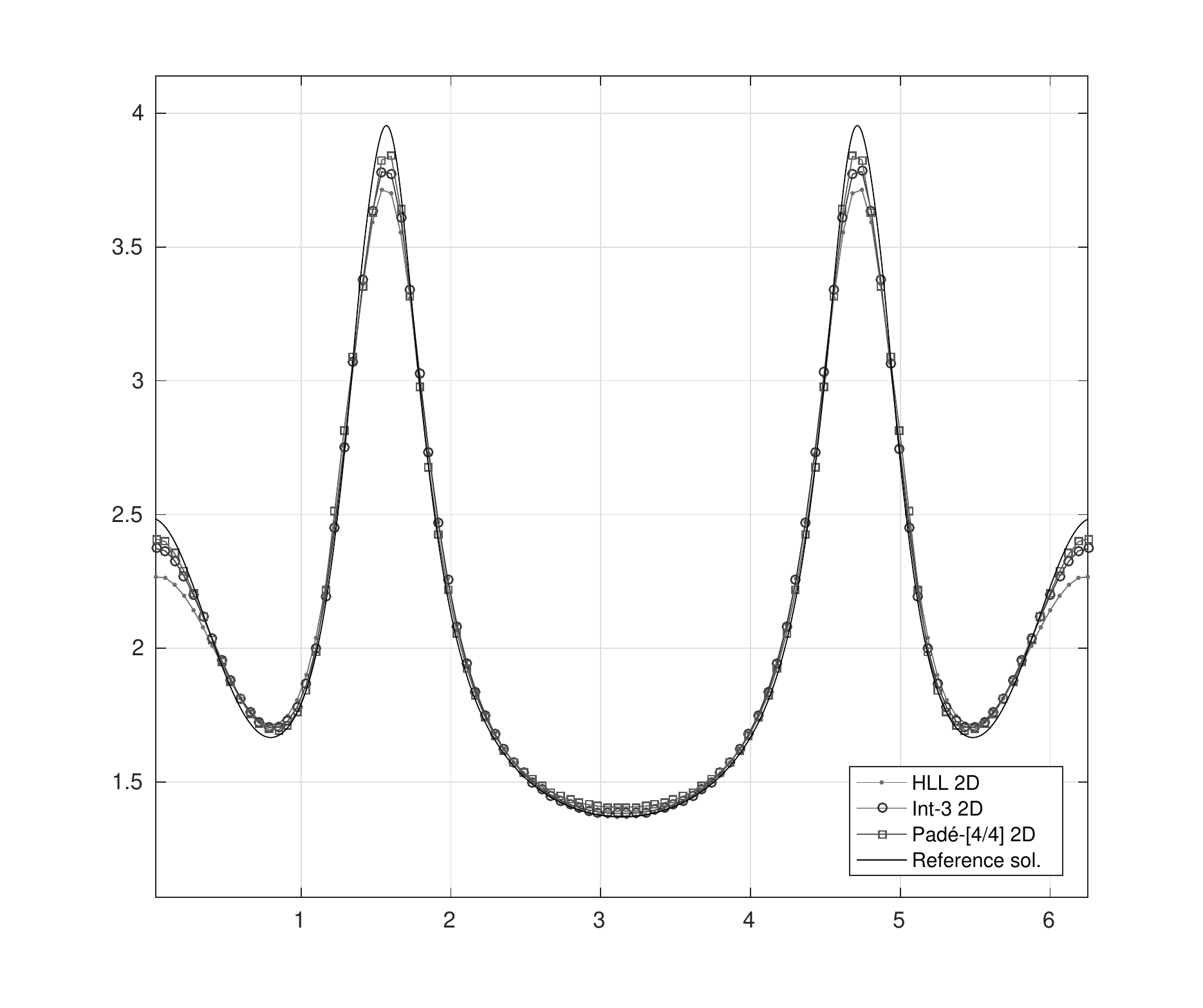}
\includegraphics[width=0.48\textwidth]{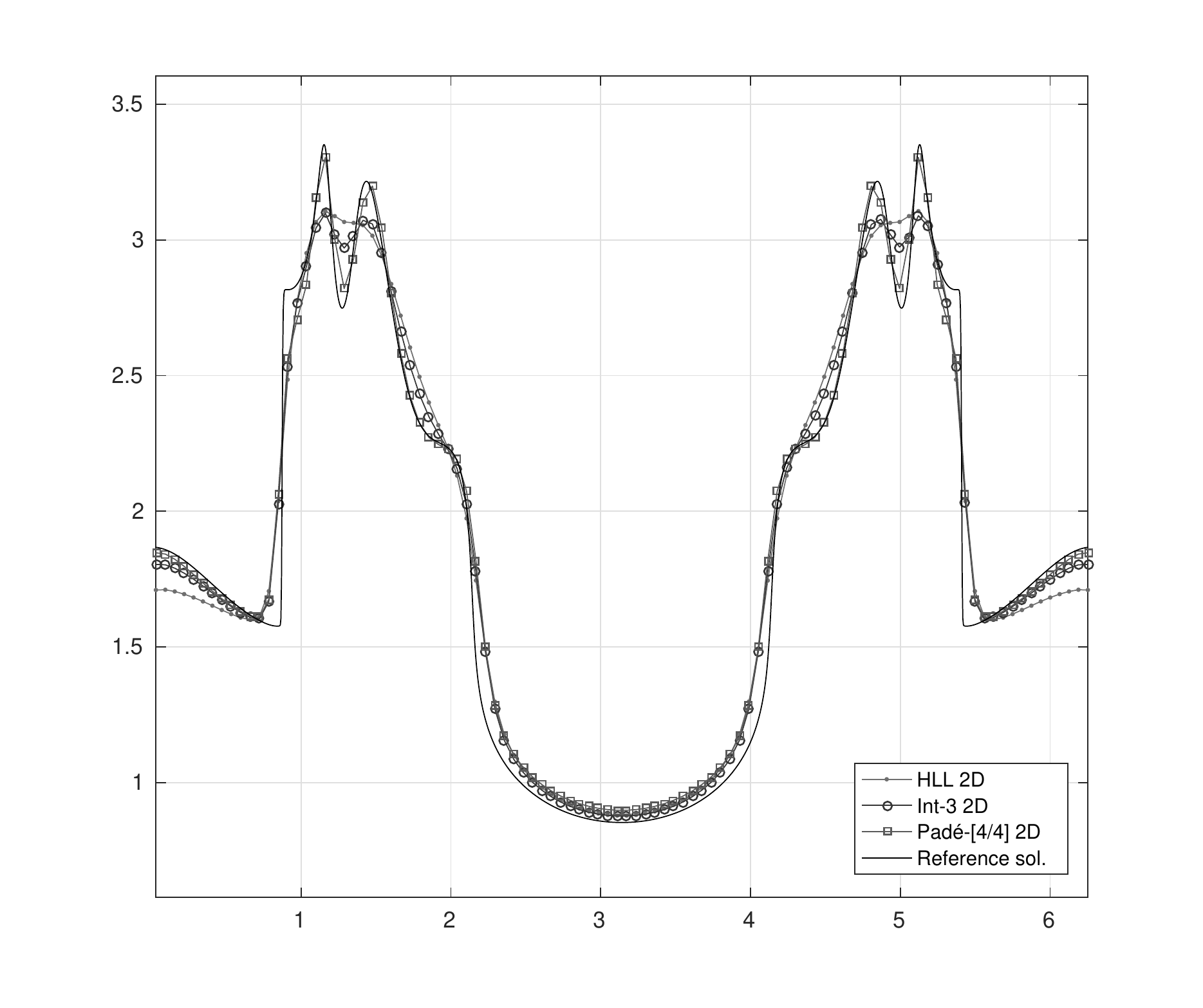} \\
\includegraphics[width=0.48\textwidth]{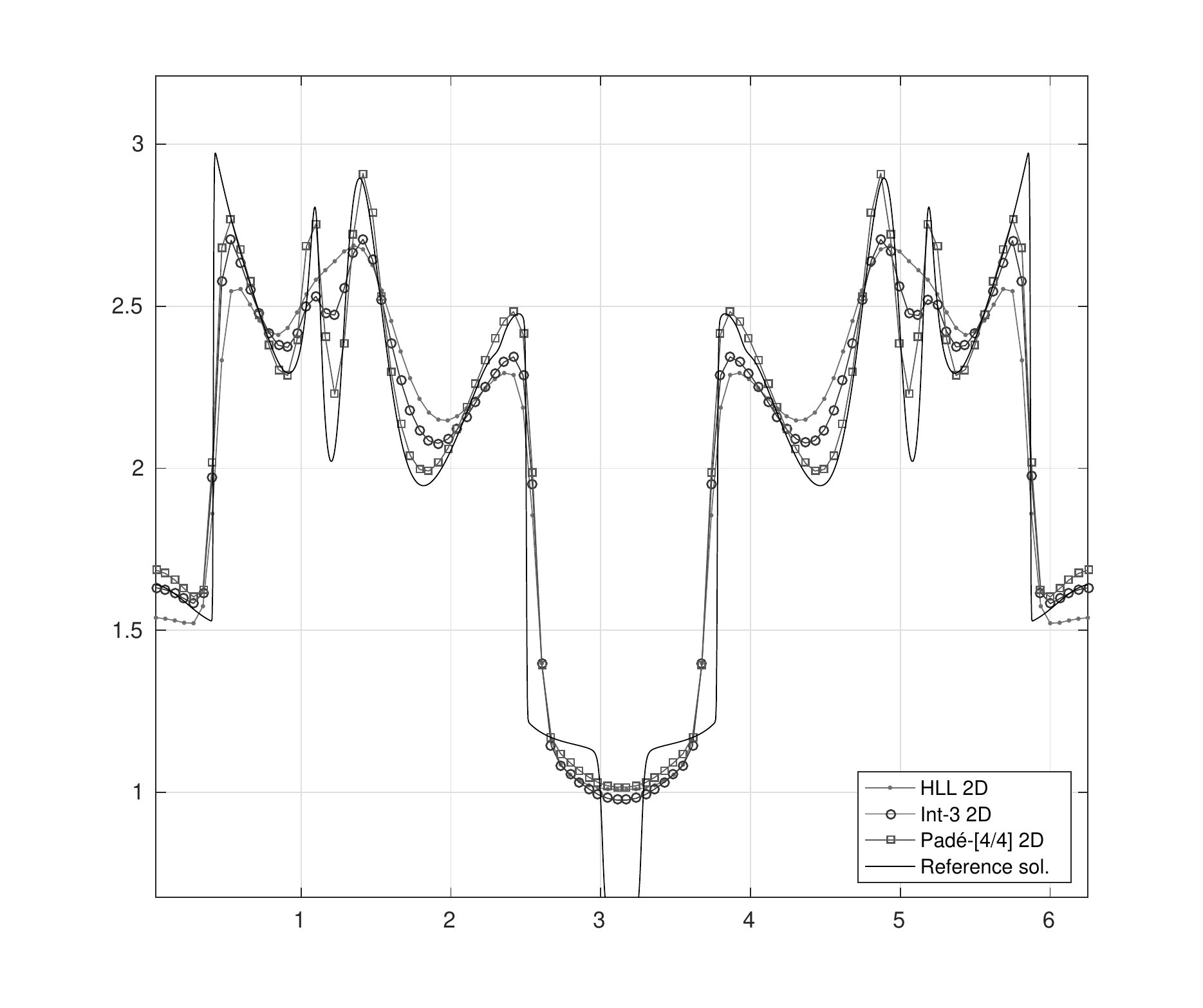}
\includegraphics[width=0.48\textwidth]{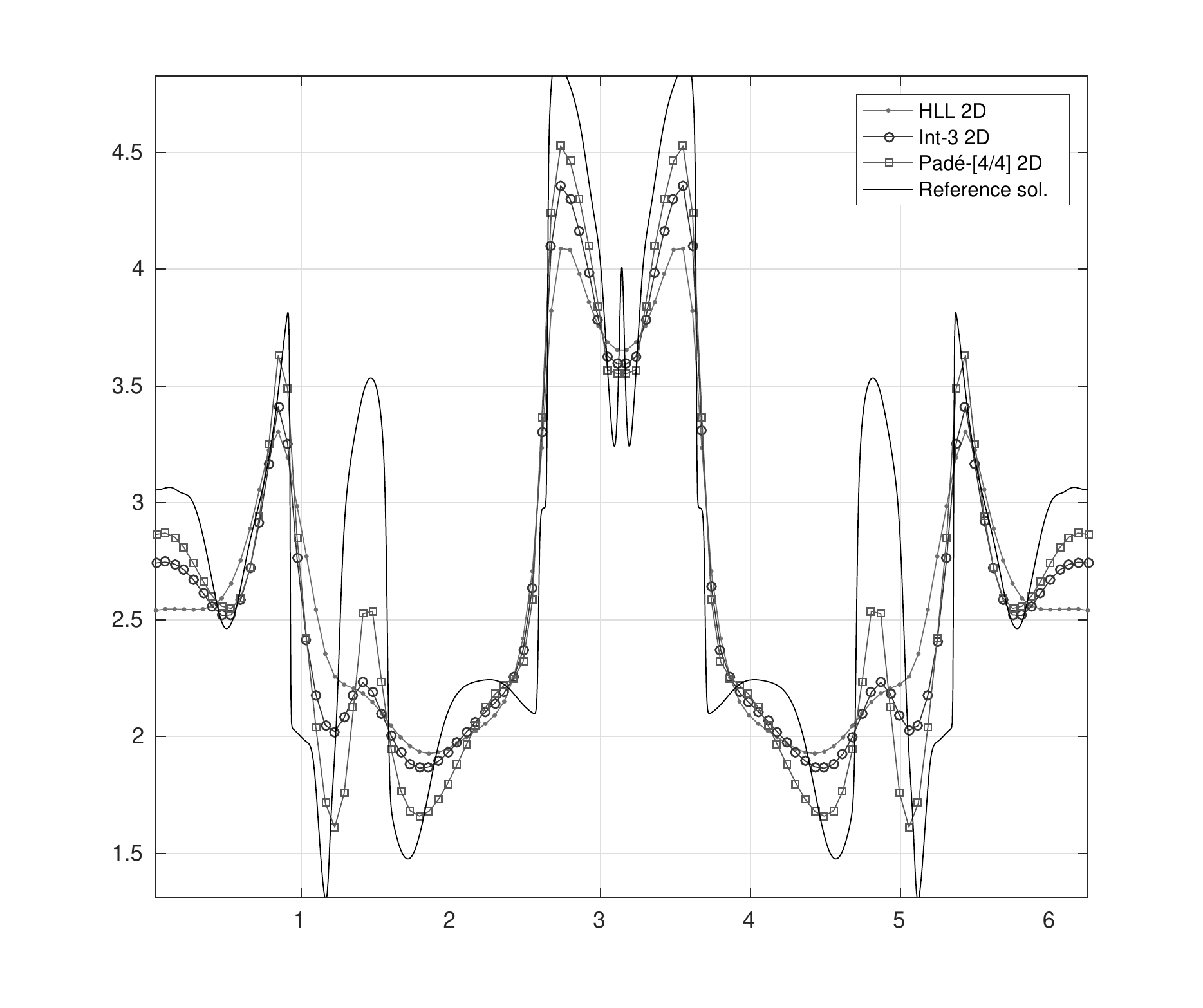}
\caption{Test \ref{test:OT}. Cuts along the main diagonal of the density at times $t=1$, $t=1.5$, $t=2$ and $t=\pi$.}
\label{fig:OTdiagonalcuts}
\end{center}
\end{figure}

Figure \ref{fig:OT} shows the densities and pressures obtained at time $t=\pi$. Although there is no accepted reference solution for this test, our results are in good agreement with those found in the literature: see, e.g., \cite{Balsara2010,FMR2009,Li-Shu}. As it can be observed, the Int-3 2D and Pad\'e-[4/4] 2D schemes provide a sharper resolution of the solution than Balsara's HLL 2D method. This can also be seen in Figure \ref{fig:OTdiagonalcuts}, where cuts along the main diagonal at different times are shown. As time evolves and the solution becomes more complex, the Pad\'e-[4/4] scheme provides significantly the best approximations, followed by Int-3 2D.

\begin{table}[!htb]
\centering
\begin{tabular}{llllll}
 Mesh size & Scheme & CPU & CPU (projection) \\ \hline
 $100\times 100$ & HLL 2D   & 1.00 &  1.27   \\
 			& Int-3 2D   & 	  5.02 &  5.29   \\
 			& Padé-[4/4] 2D & 9.40 &  9.91   \\ 
 			\hline
 $200\times 200$ & HLL 2D   & 1.00 &  1.74   \\
 			& Int-3 2D   & 	  5.13 &  5.62   \\
 			& Padé-[4/4] 2D & 9.78 &  10.50   \\ 
 			\hline
 $300\times 300$ & HLL 2D   & 1.00 &  1.53   \\
 			& Int-3 2D   & 	  5.72 &  5.94   \\
 			& Padé-[4/4] 2D & 10.37 &  10.70   \\ 
 			\hline
 $400\times 400$ & HLL 2D   & 1.00 &  1.35   \\
 			& Int-3 2D   & 	  5.22 &  5.40   \\
 			& Padé-[4/4] 2D & 9.86 &  10.32   \\ 
 			\hline
 $500\times 500$ & HLL 2D   & 1.00 &  1.21   \\
 			& Int-3 2D   & 	  5.16 &  5.27   \\
 			& Padé-[4/4] 2D & 9.84 &  10.14   \\ 
 			\hline
\end{tabular}
\caption{Test \ref{test:OT}. Relative CPU times. Comparison between the divergence cleaning technique in Section \ref{sec:MHD} (left) and the projection method (right) on several meshes.}
\label{table:OT}
\end{table}

\begin{figure}[!ht]
\begin{center}
\includegraphics[width=0.48\textwidth]{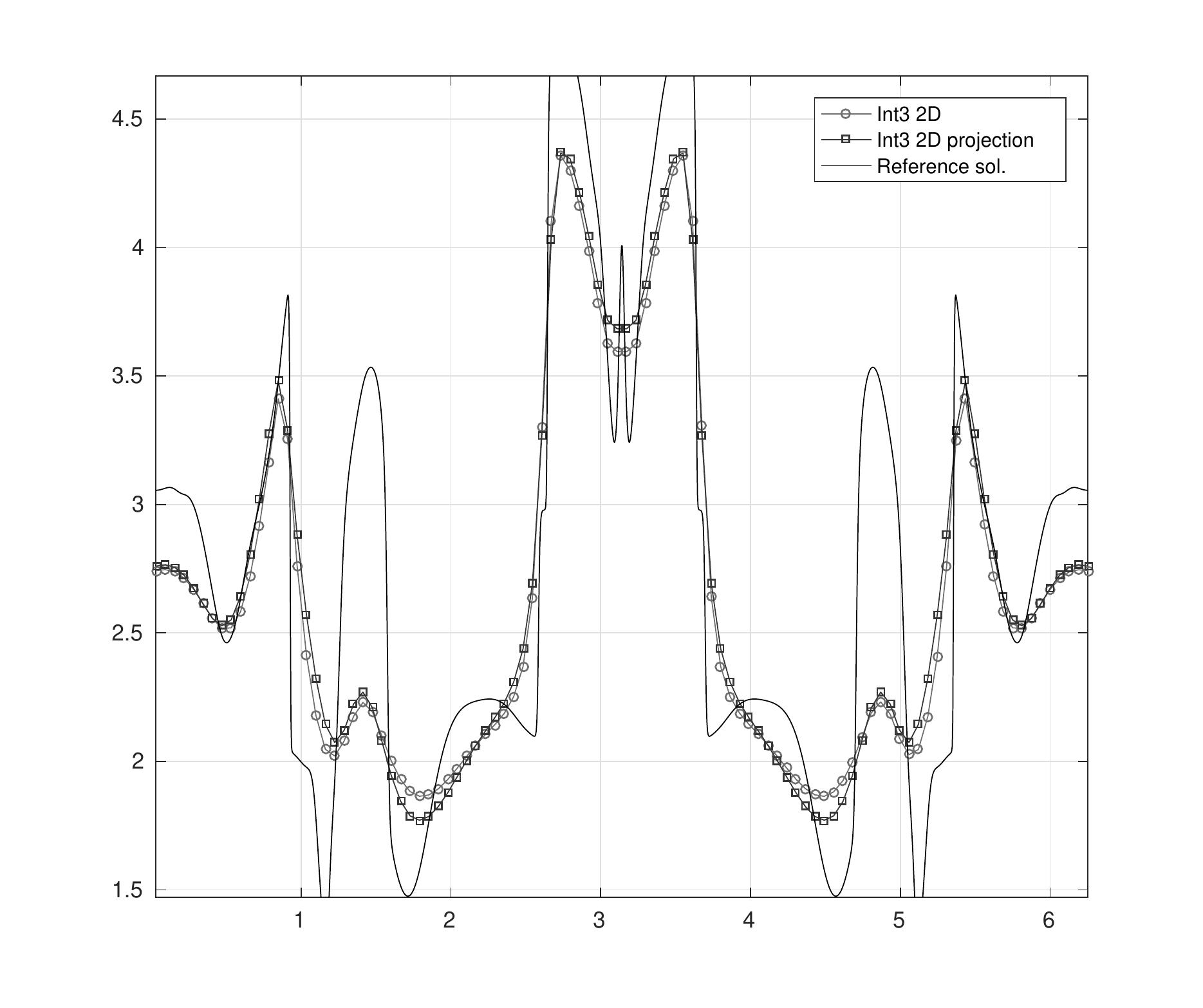}
\includegraphics[width=0.48\textwidth]{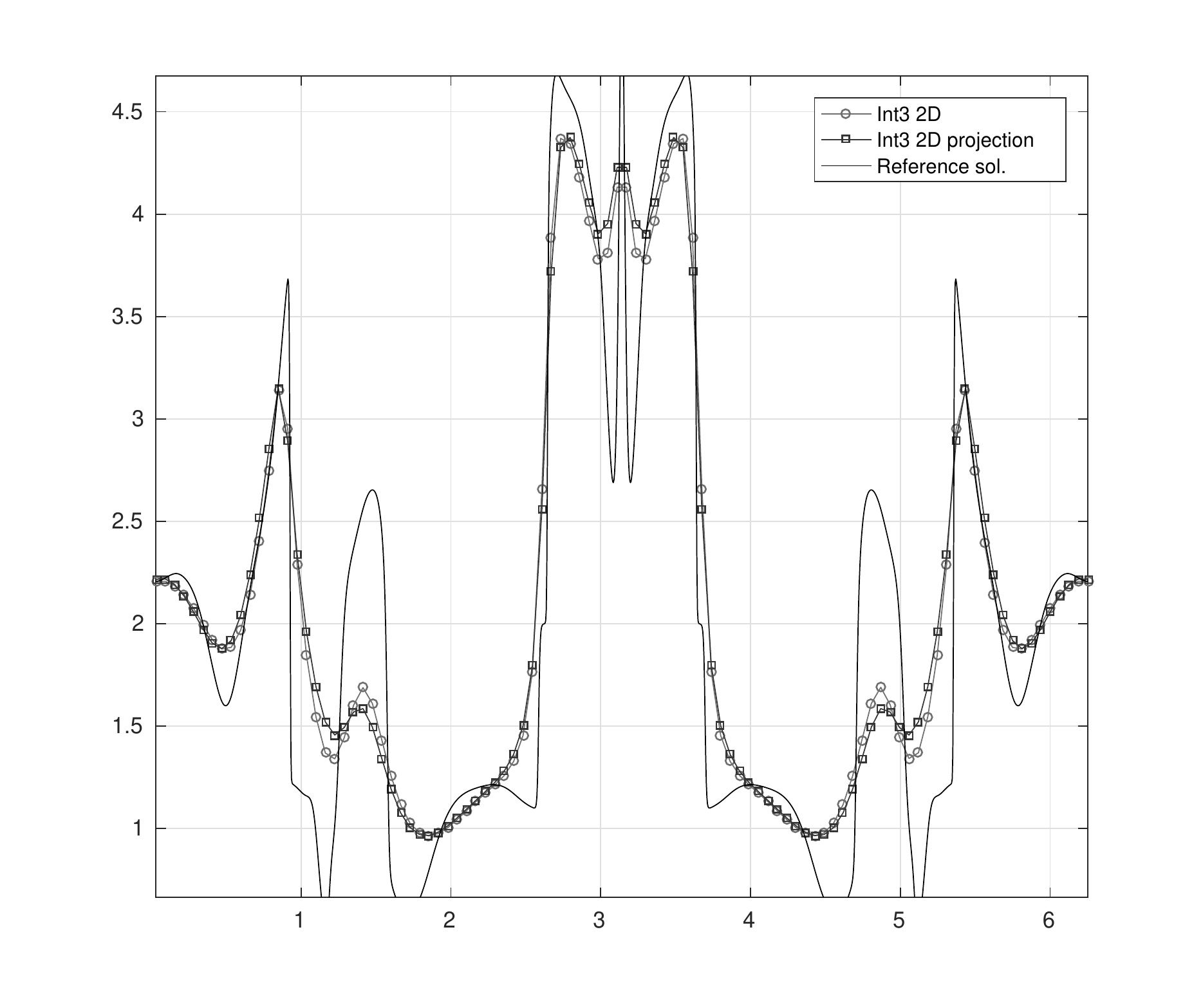}
\caption{Test \ref{test:OT}. Comparison of the divergence cleaning technique of Section \ref{sec:MHD} with the projection method using the Int-3 2D scheme. Cut along the main diagonal: density (left) and pressure (right) computed at time $t=\pi$ on a $100\times 100$ mesh.}
\label{fig:comparisonOT}
\end{center}
\end{figure}

On the other hand, we have compared the technique for imposing the divergence constraint defined in Section \ref{sec:MHD} with the well-known projection method introduced in \cite{BB1980}. Figure \ref{fig:comparisonOT} shows cuts along the main diagonal of the solutions obtained with the Int-3 2D scheme at time $t=0.5$ with both divergence cleaning techniques. As it can be observed, the results are comparable. Also, Table \ref{table:OT} shows the relative CPU times obtained when using both divergence cleaning methods. The results are again comparable, although our method seems to be slighly faster.

%\begin{figure}[!ht]
%\begin{center}
%\includegraphics[width=0.8\textwidth]{OT_comparison.pdf}
%\caption{Test \ref{test:OT}. $L^1$ norm of $\nabla\cdot\BB$ on the time interval $[0, \pi]$. Comparison between the divergence cleaning technique in Section \ref{sec:MHD} and the projection method \cite{BB1980}, using the Int-3 2D scheme.}
%\label{fig:divB}
%\end{center}
%\end{figure}

%\begin{table}[ht]
%\caption{Test \ref{test:OT}: Relative CPU times with respect to the first-order OS solver. Final time: $t=0.2$.}
%\begin{center}
%\begin{tabular}{lcc}
%Method & CPU (first order) & CPU (third order)\\ 
%\hline
%OS & $1.00$ & $5.82$ \\
%OS-Cheb-4 & $0.16$ & $1.04$ \\
%OS-Newman-4 & $0.38$ & $2.32$ \\
%OS-Halley-2 & $0.50$ & $2.79$ \\
%HLL & $0.05$ & $0.36$ \\
%\hline
%\end{tabular}
%\end{center}
%\label{table:OT}
%\end{table}

%%

\subsection{The rotor problem} \label{test:rotor}

The rotor problem, initially proposed in \cite{Balsara-Spicer}, is considered in this section. The initial condition consists of a dense rotating disk at the center of the domain, which is connected by means of a taper function with the ambient fluid at rest. The rotor is not in equilibrium, as the centrifugal forces are not balanced. Eventually, the rotating dense fluid will be confined into an oblate shape due to the action of the magnetic field.

Given $r_0=0.1$, $r_1=0.115$, $f=(r_1-r)/(r_1-r_0)$ and $r=[(x-0.5)^2+(y-0.5)^2]^{1/2}$, the initial conditions are given by
\[
(\rho, v_x, v_y)=
\begin{cases}
(10, -(y-0.5)/r_0, (x-0.5)/r_0) & \text{if } r<r_0, \smallskip \\
(1+9f, -(y-0.5)f/r, (x-0.5)f/r) & \text{if } r_0<r<r_1, \smallskip \\
(1, 0, 0) & \text{if } r>r_1,
\end{cases}
\]
with $B_x=2.5/\sqrt{4\pi}$, $B_y=0$ and $P=0.5$. We take $\gamma=5/3$. The computational domain is $[0, 1]\times [0, 1]$ and transmissive boundary conditions are considered.

\begin{figure}[!ht]
\begin{center}
\includegraphics[height=0.24\textheight]{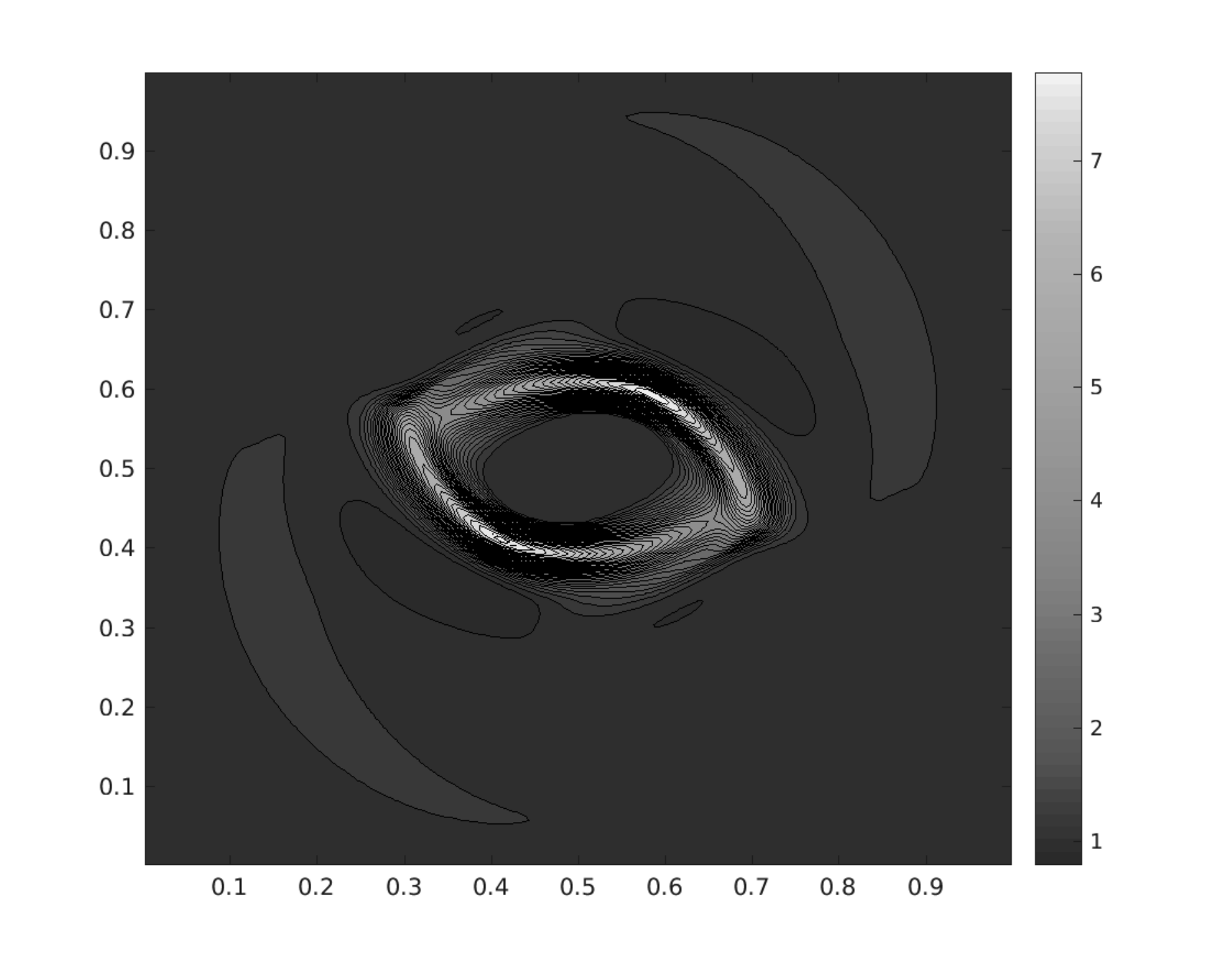}
\includegraphics[height=0.24\textheight]{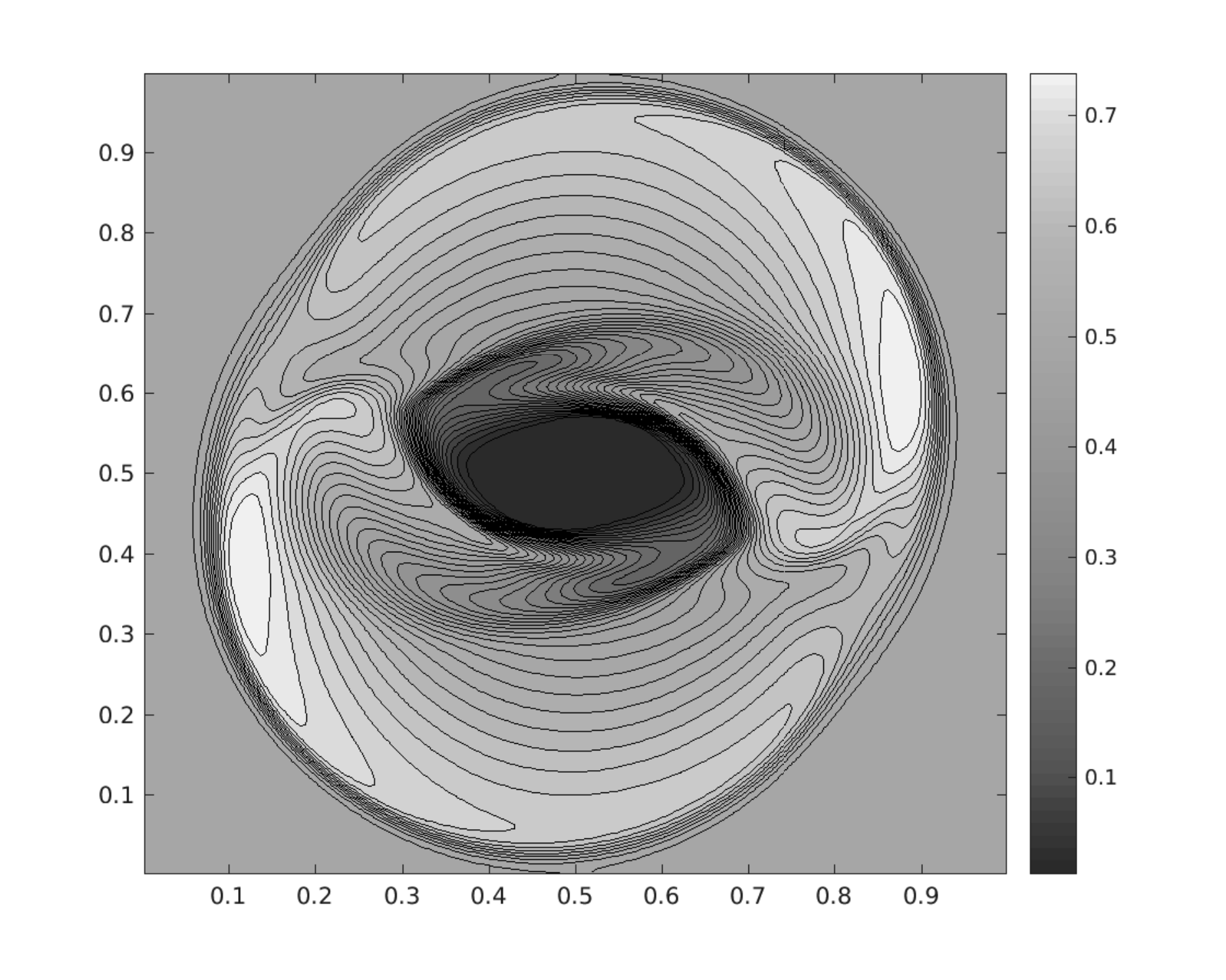} \\
\includegraphics[height=0.24\textheight]{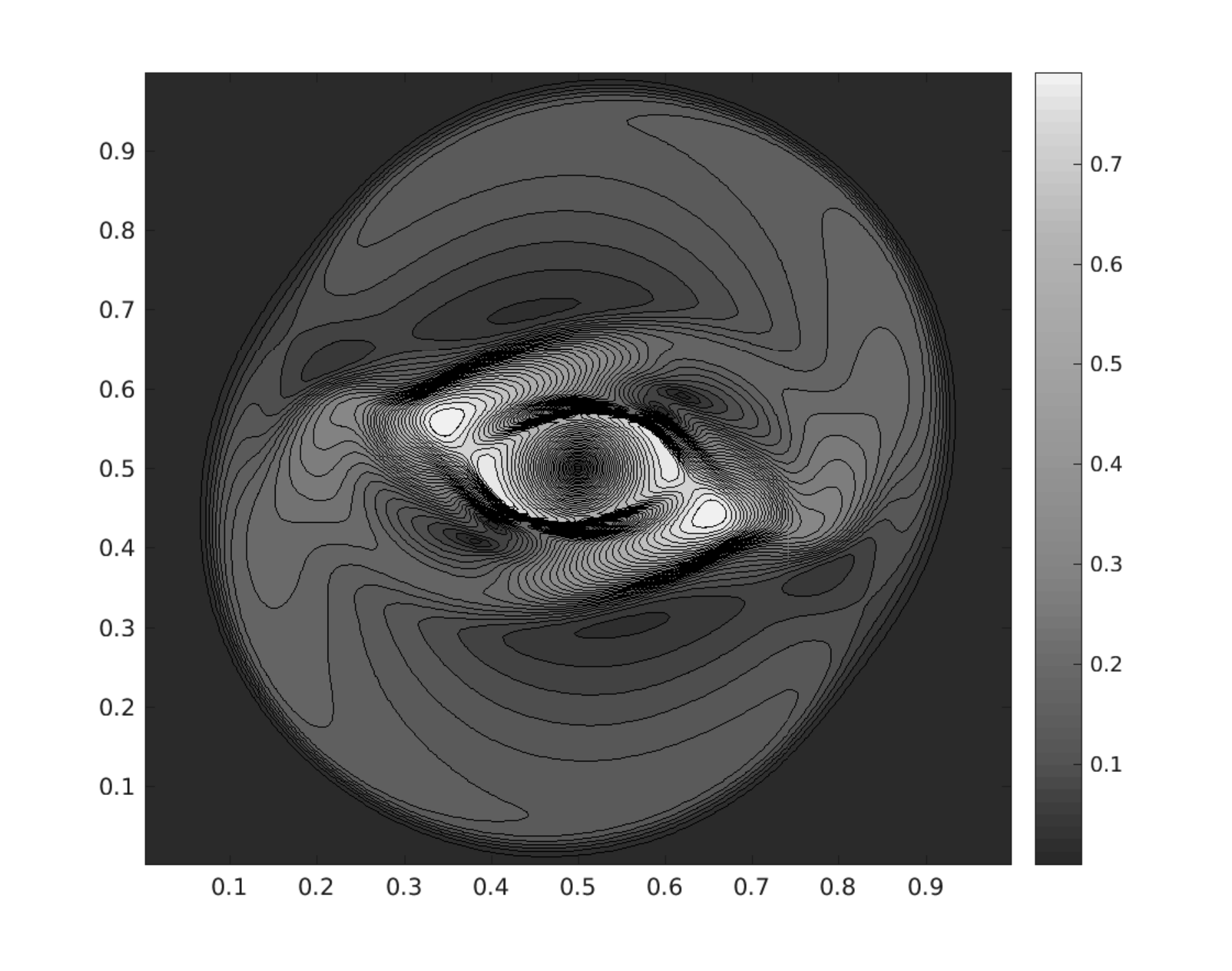}
\includegraphics[height=0.24\textheight]{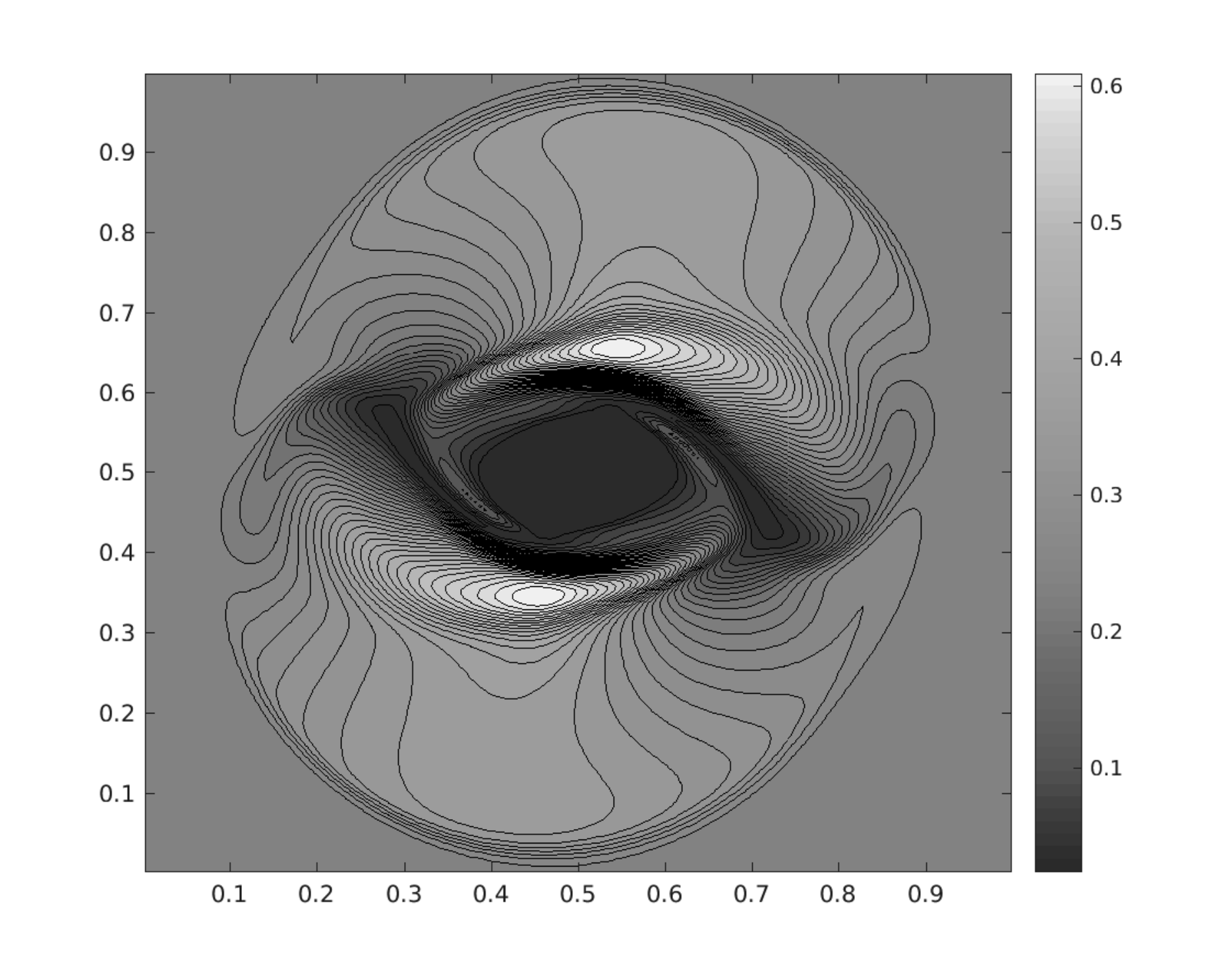}
\caption{Test \ref{test:rotor}. Density $\rho$ (top left), pressure $P$ (top right), Mach number $|\mathbf{v}|/a$ (bottom left) and magnetic pressure $|\mathbf{B}|^2/2$ (bottom right) computed at time $t=0.295$. Results obtained with the Pad\'e-[4/4] 2D scheme with $200\times 200$ cells.}
\label{fig:rotor}
\end{center}
\end{figure}

Figure \ref{fig:rotor} shows the solution obtained with the Pad\'e-[4/4] 2D scheme at time $t=0.295$ on a $200\times 200$ mesh with CFL$=0.8$. The results are in good agreement with those presented in \cite{Balsara-Spicer,RVM_Osher2016,Li-Shu,Toth2000}. In particular, it is worth noticing that the results produced with the first-order Pad\'e-[4/4] 2D scheme have a similar quality as those obtained in \cite{RVM_Osher2016} with a one-dimensional third-order approximate Chebyshev-DOT scheme on the same mesh.

\subsection{Two-dimensional Riemann problem} \label{test:2d-RP}

In this section we consider a two-dimensional Riemann problem originally proposed in \cite{Dedner}. The initial condition, given in Table \ref{table:2d-RP}, is chosen so that the solutions of three of the four one-dimensional Riemann problems are simple waves. Denoting the quadrants by Roman numbers as in Table \ref{table:2d-RP}, there is a  rarefaction wave for problem I$\leftrightarrow$II and shocks for II$\leftrightarrow$III and III$\leftrightarrow$IV.

\begin{table}[ht]
\caption{Test \ref{test:2d-RP}. Initial data for the 2d Riemann problem.}
\begin{center}
\begin{tabular}{rllllllll}
\multicolumn{1}{c}{Quadrant} & \multicolumn{1}{c}{$\rho$} & \multicolumn{1}{c}{$\rho v_x$} & \multicolumn{1}{c}{$\rho v_y$} & \multicolumn{1}{c}{$\rho v_z$} & \multicolumn{1}{c}{$B_x$} & \multicolumn{1}{c}{$B_y$} & \multicolumn{1}{c}{$B_z$} & \multicolumn{1}{c}{$E$} \smallskip \\ 
\hline
I: $x>0$, $y>0$ & 0.9308 & 1.4557 & -0.4633 & 0.0575 & 0.3501 & 0.9830 & 0.3050 & 5.0838 \smallskip \\
II: $x<0$, $y>0$ & 1.0304 & 1.5774 & -1.0455 & -0.1016 & 0.3501 & 0.5078 & 0.1576 & 5.7813 \smallskip \\
III: $x<0$, $y<0$ & 1.0000 & 1.7500 & -1.0000 & 0.0000 & 0.5642 & 0.5078 & 0.2539 & 6.0000 \smallskip \\
IV: $x>0$, $y<0$ & 1.8887 & 0.2334 & -1.7422 & 0.0733 & 0.5642 & 0.9830 & 0.4915 & 12.999
\end{tabular}
\end{center}
\label{table:2d-RP}
\end{table}

\begin{figure}[!ht]
\begin{center}
\includegraphics[height=0.24\textheight]{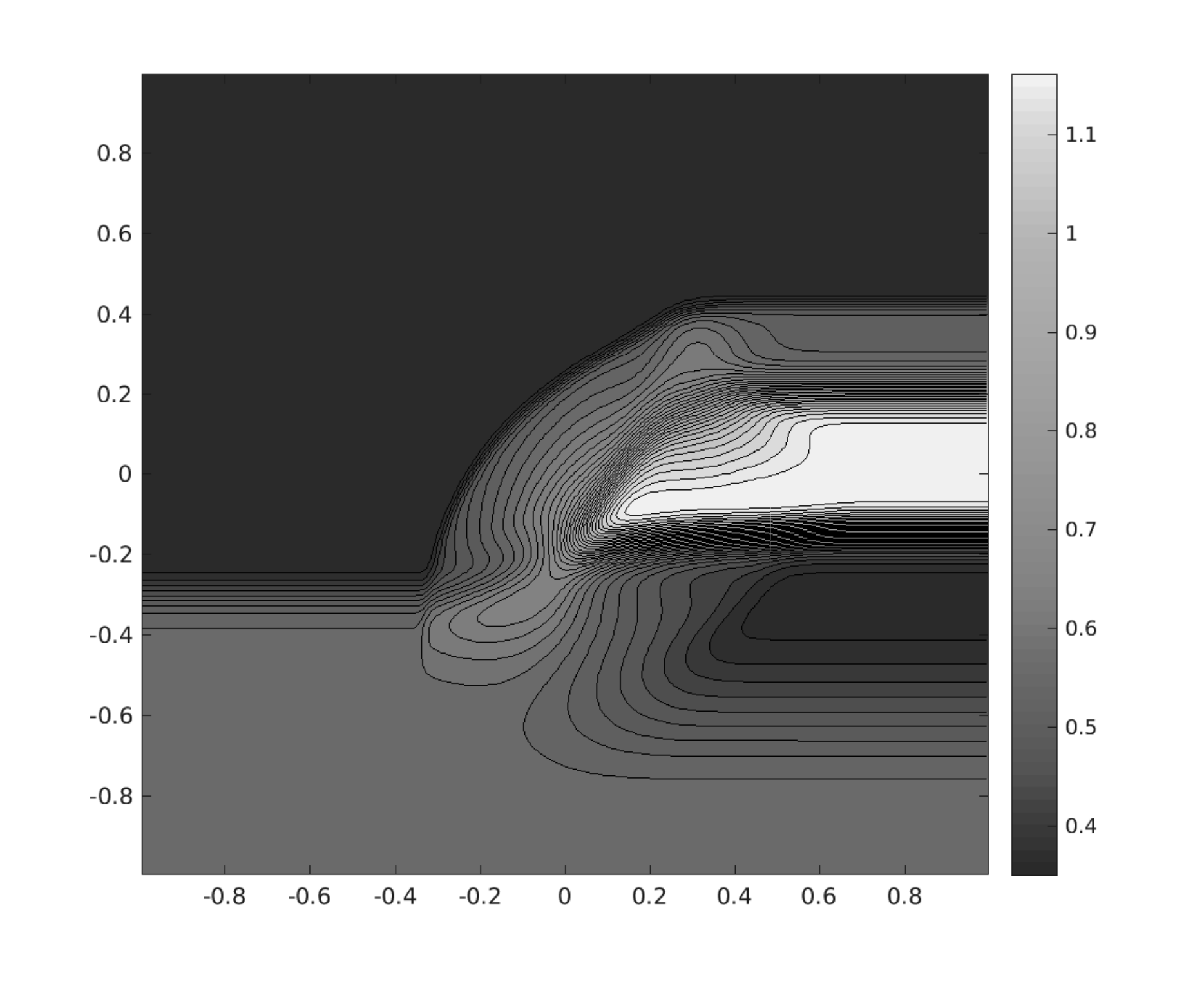}
\includegraphics[height=0.24\textheight]{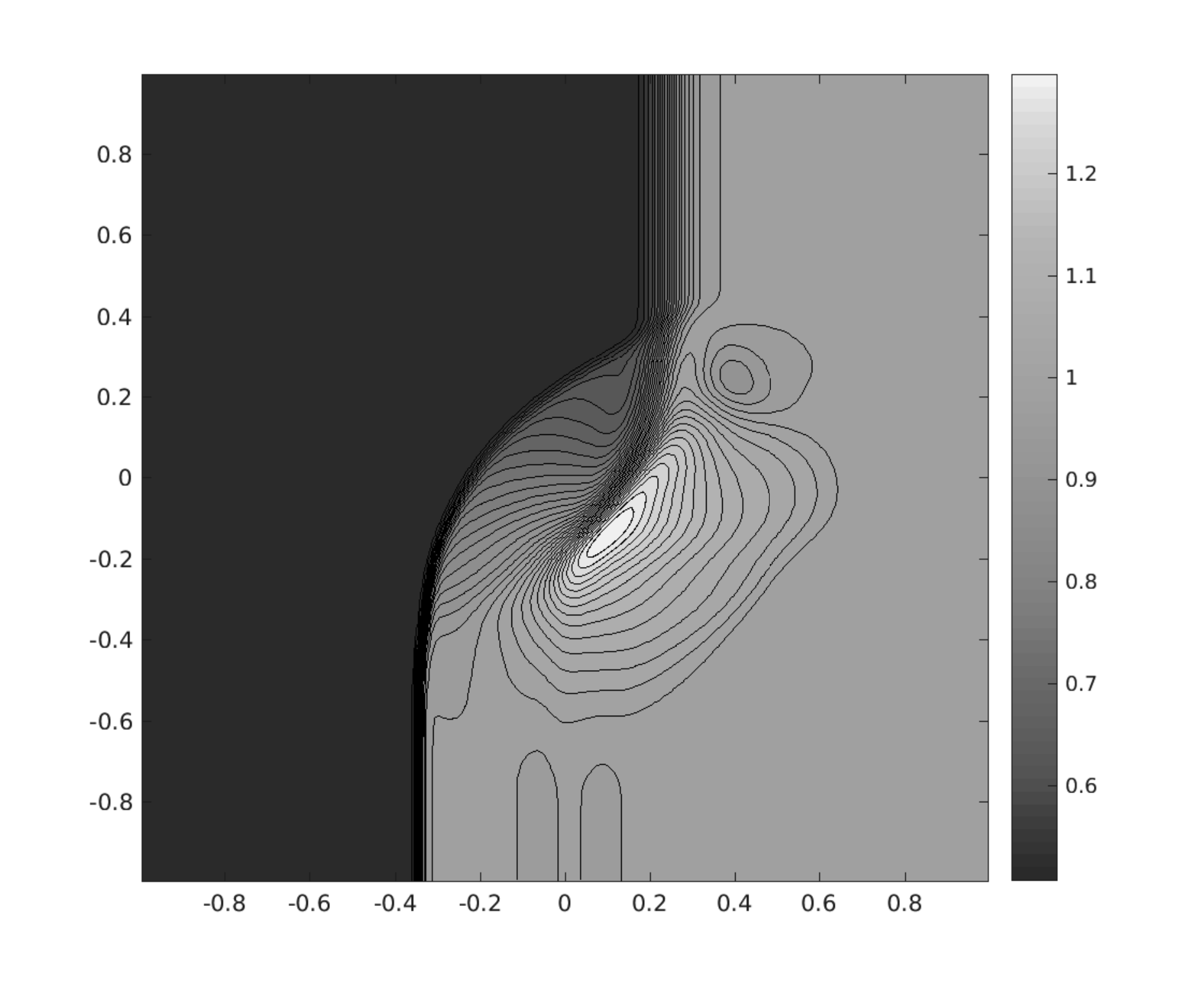} 
\caption{Test \ref{test:2d-RP}. Contours at time $t=0.2$ obtained with the Pad\'e-[4/4] 2D scheme on a $200\times 200$ mesh. Left: $B_x$. Right: $B_y$.}
\label{fig:2d-RP}
\end{center}
\end{figure}

The problem has been solved in the computational domain $[-1, 1]\times [-1, 1]$ using a $200\times 200$ mesh until a final time $t=0.2$, with $\gamma=5/3$. Following \cite{Li-Shu}, we have considered Neumann boundary conditions.

Figure \ref{fig:2d-RP} shows the contours of $B_x$ and $B_y$ computed with the Pad\'e-[4/4] 2D scheme using a CFL of $0.8$; for this test, the Int-3 2D scheme gives similar results. It was reported in \cite{Dedner} that some schemes present problems to keep $B_y$ constant across the shock in II$\leftrightarrow$III, and they can also produce strong distortions in $B_x$ and $B_y$ behind the rarefaction wave in I$\leftrightarrow$II. As it can be observed in Figure \ref{fig:2d-RP}, none of those pathologies appear in the computed solutions, that can be directly compared with those presented in \cite{Dedner,Li-Shu}.

\begin{figure}[!ht]
\begin{center}
\includegraphics[width=0.48\textwidth]{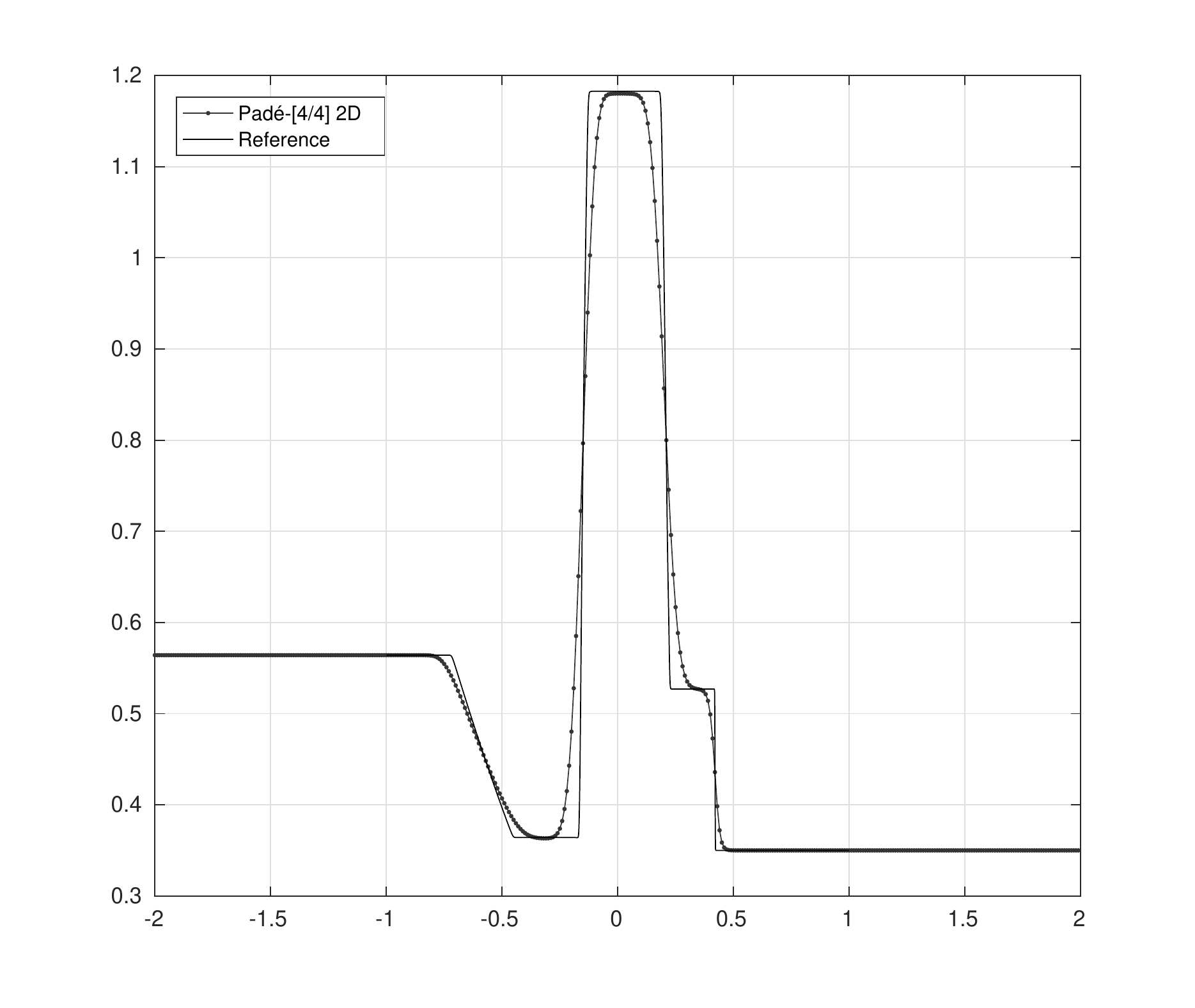}
\includegraphics[width=0.48\textwidth]{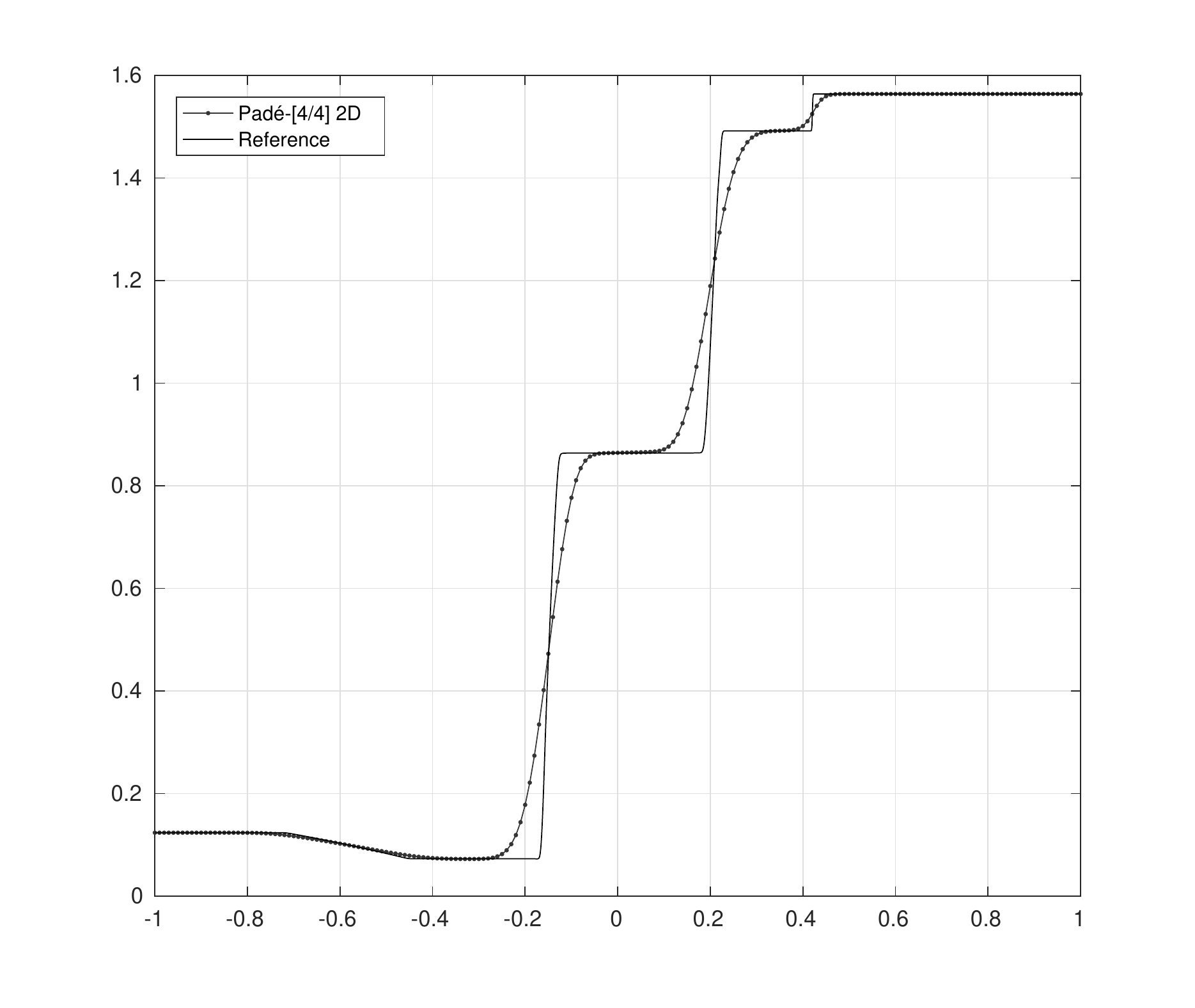}
\caption{Test \ref{test:2d-RP}. Cuts at $x=0.93$ and $t=0.2$. Solid line: reference solution. Dots: Pad\'e-[4/4] 2D. Left: $B_x$. Right: $v_x$.}
\label{fig:2d-RP_cuts}
\end{center}
\end{figure}

Finally, Figure \ref{fig:2d-RP_cuts} shows a comparison between the solution of the two-dimensional Riemann problem obtained with the Pad\'e-[4/4] 2D scheme and a one-dimensional reference solution. This serves to compare the quality of the solution of the I$\leftrightarrow$IV Riemann problem at $x=0.93$, in a similar way as it was done in \cite{Dedner,Li-Shu}.

\subsection{Spherical explosion} \label{test:explosion}

This test concerns the evolution of an explosion produced by an overpressured spherical region located at the center of the domain. It was originally proposed in \cite{Zachary}, although we have considered here the data proposed in \cite{Tang-Xu}. Specifically, the initial pressure is $P=100$ inside a circle of radius $r=10$ and $P=1$ outside. The density is $\rho=1$ on the whole domain and the velocity field is set to zero. As in \cite{Tang-Xu}, three different values of the initial magnetic field are considered, with a common value of the adiabatic constant $\gamma=2$. The computations have been performed with the Pad\'e-[4/4] 2D scheme on the computational domain $[-50, 50]\times [-50, 50]$, using a $400\times 400$ uniform mesh and CFL=$0.8$.

\begin{figure}[!ht]
\begin{center}
\includegraphics[height=0.24\textheight]{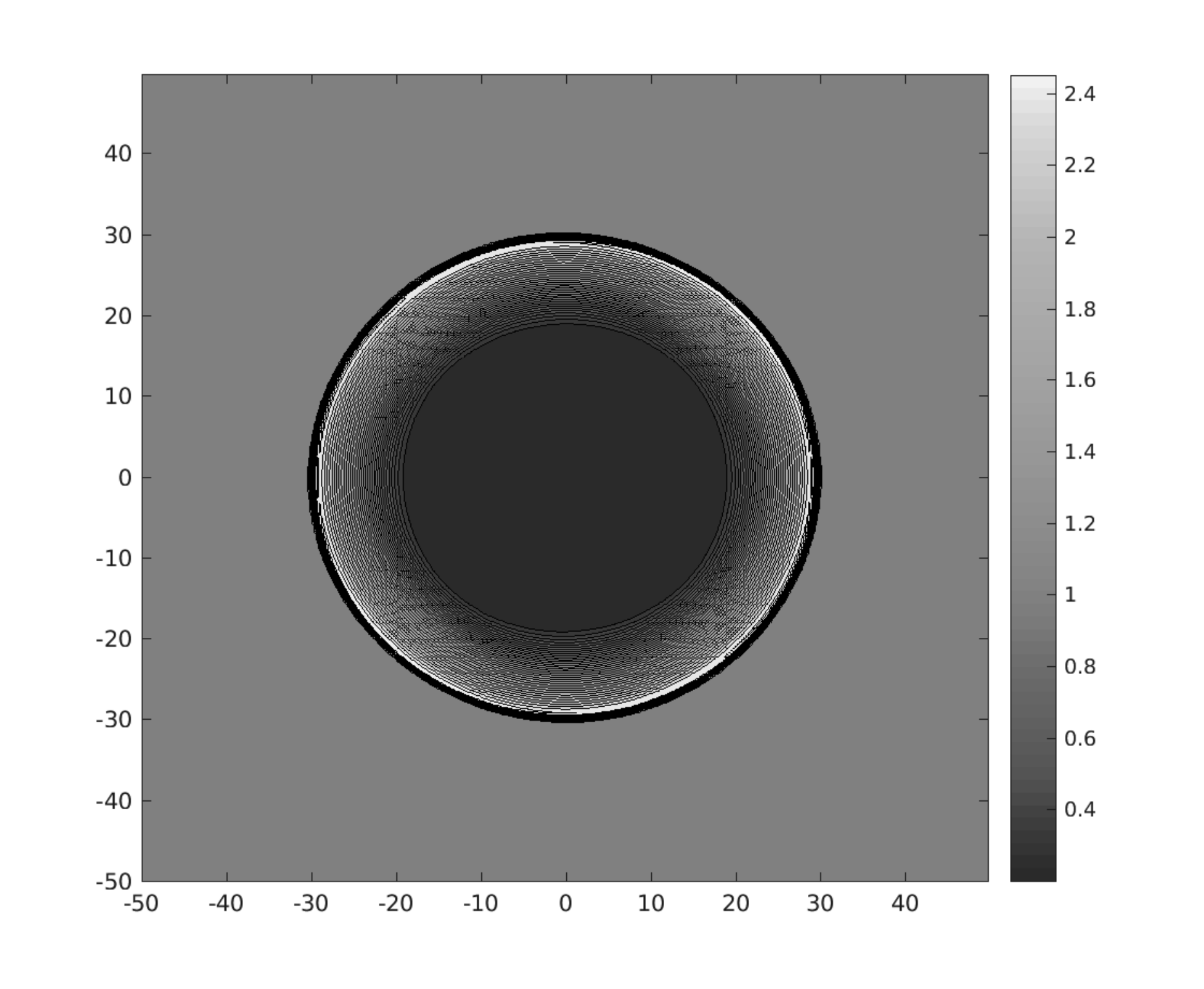}
\includegraphics[height=0.24\textheight]{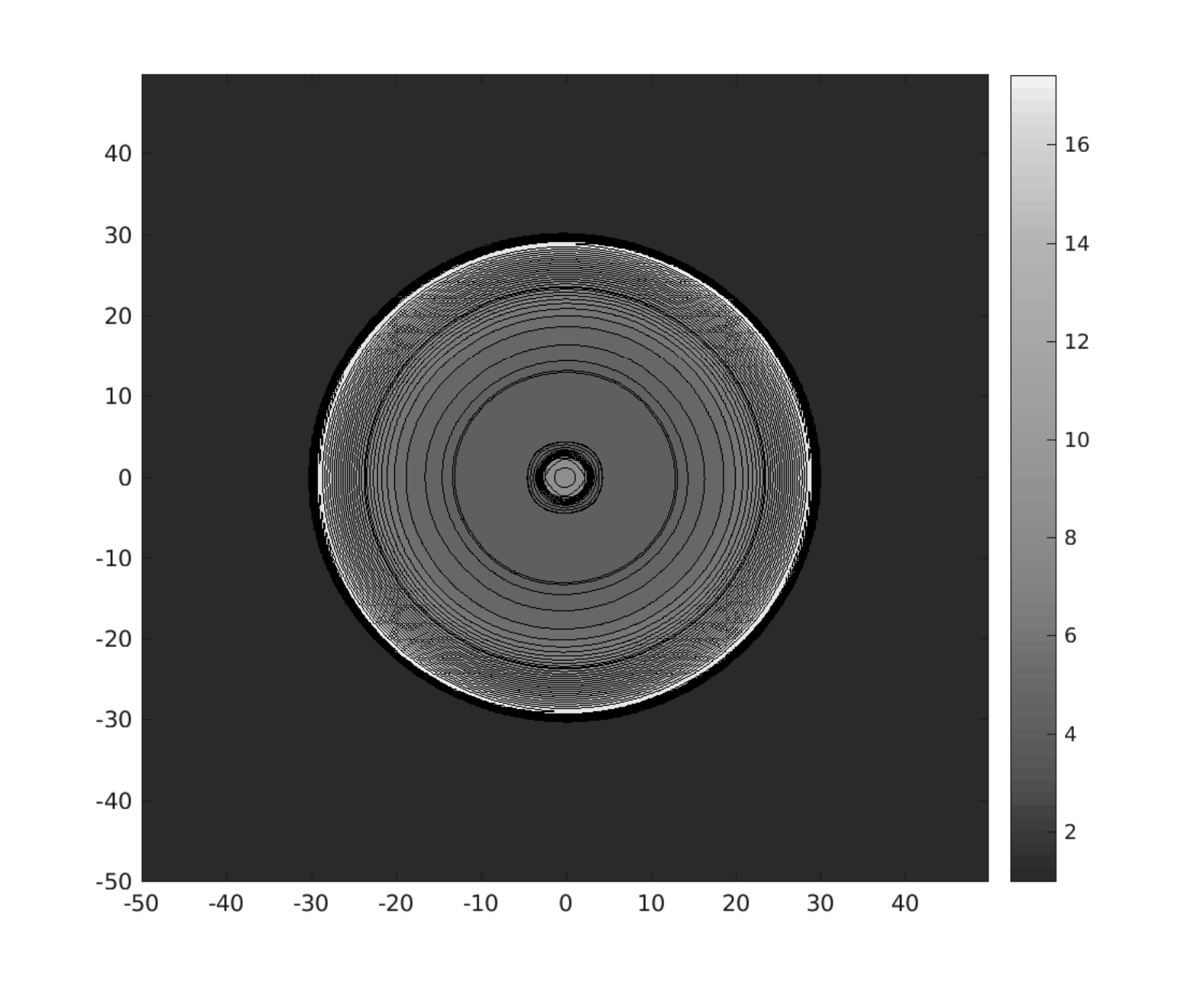} \\
\includegraphics[height=0.24\textheight]{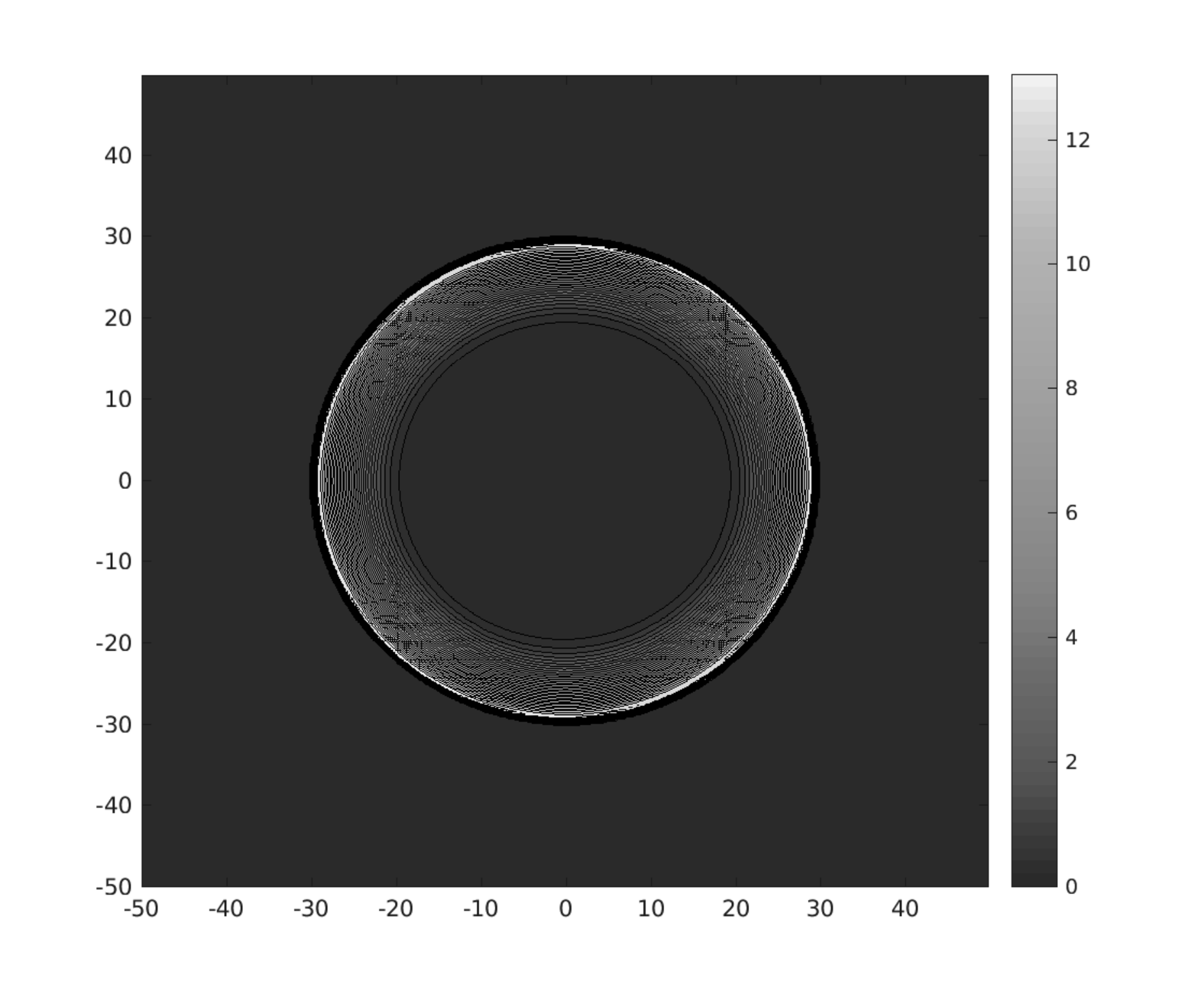}
\includegraphics[height=0.24\textheight]{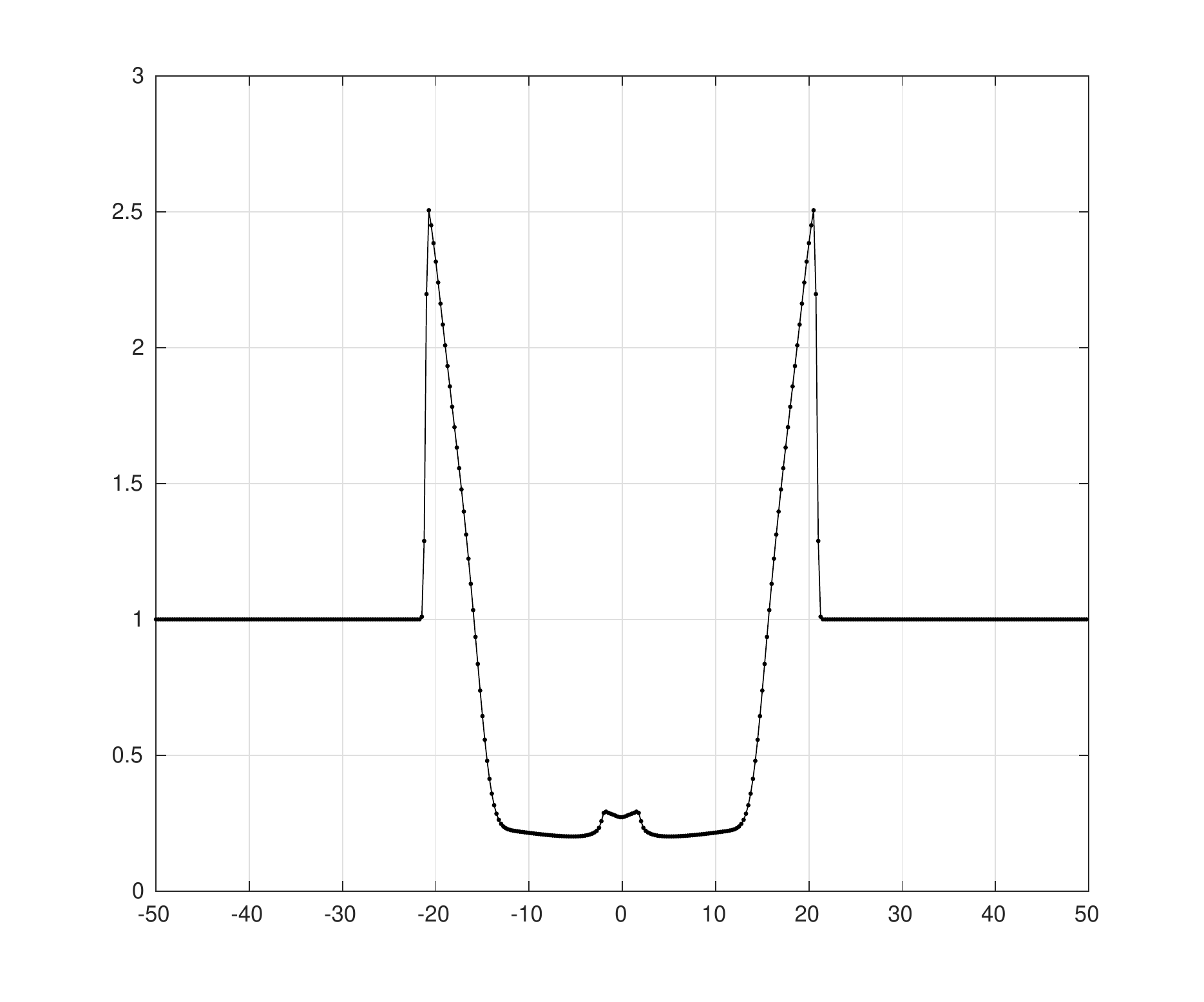}
\caption{Test \ref{test:explosion}. Solution obtained with initial magnetic field $(B_x, B_y, B_z)=(0, 0, 0)$. Density $\rho$ (top left), pressure $P$ (top right), kinetic energy (bottom left) and density along $y=0.5$ (bottom right) computed at time $t=3$.}
\label{fig:explosion1}
\end{center}
\end{figure}

Taking the initial magnetic field as $(B_x, B_y, B_z)=(0, 0, 0)$ the solution propagates symmetrically in the radial direction, as expected for the Euler equations. This behavior can be seen in Figure \ref{fig:explosion1}, which shows the results of the simulation at time $t=3$.

\begin{figure}[!ht]
\begin{center}
\includegraphics[height=0.24\textheight]{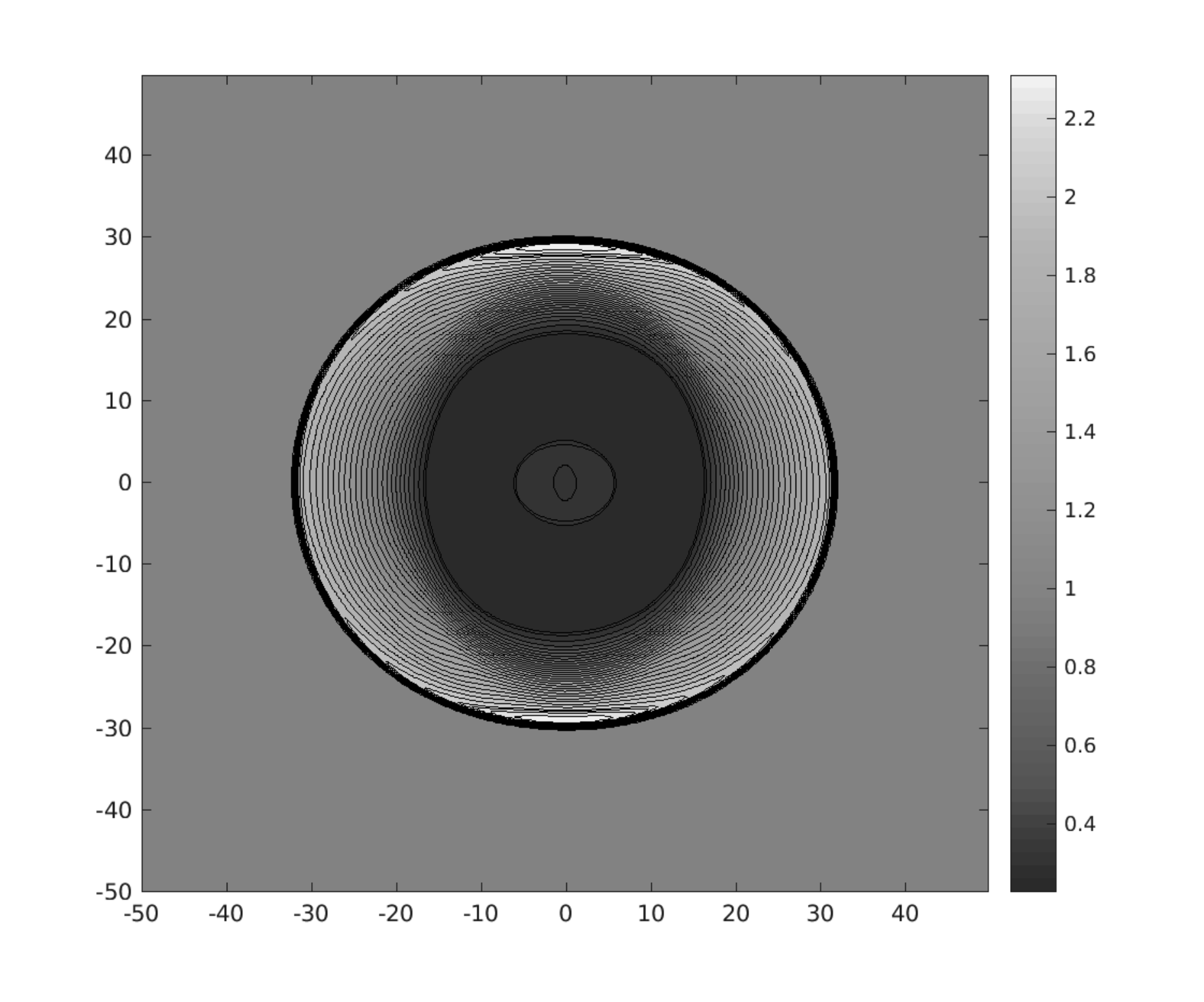}
\includegraphics[height=0.24\textheight]{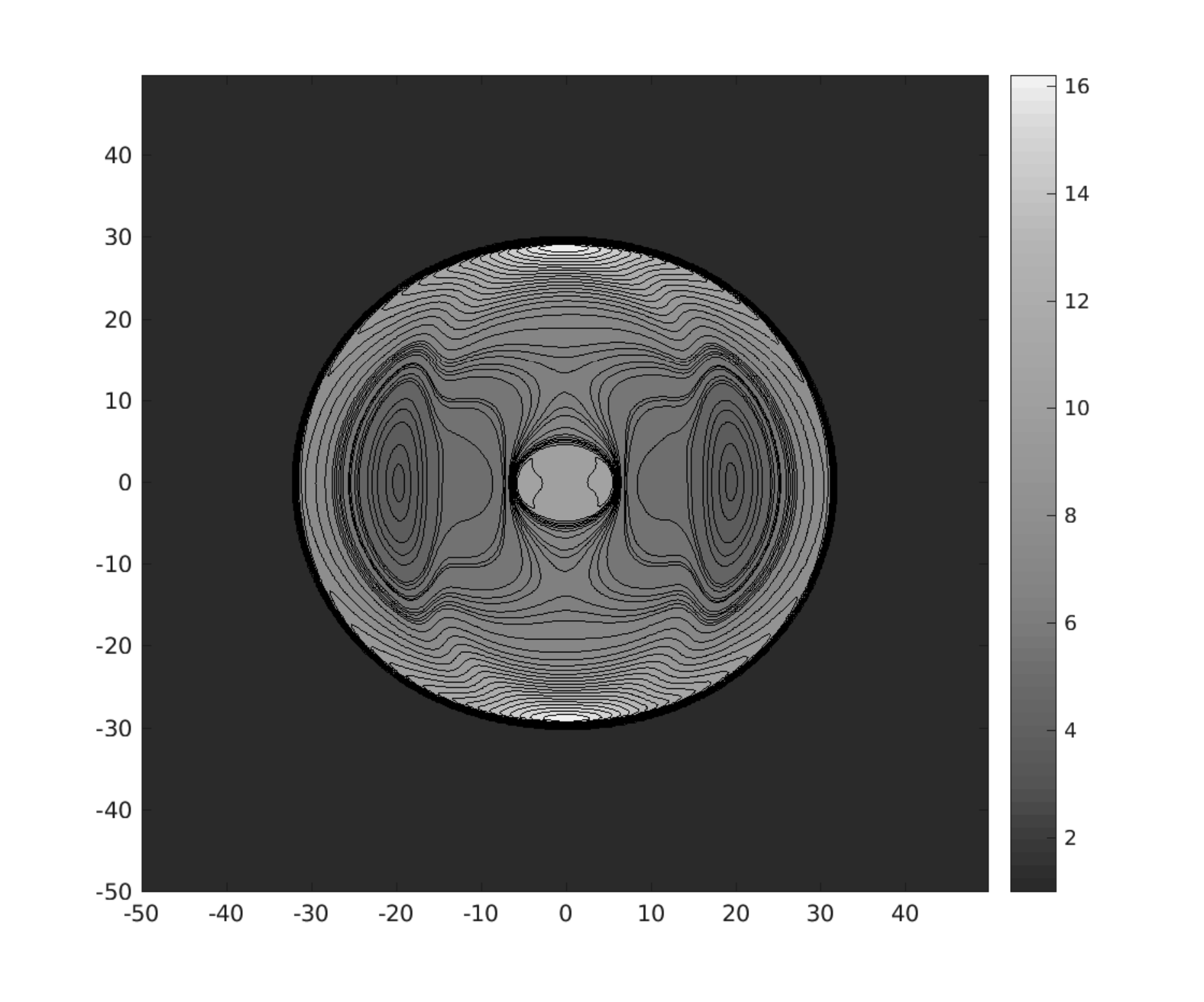} \\
\includegraphics[height=0.24\textheight]{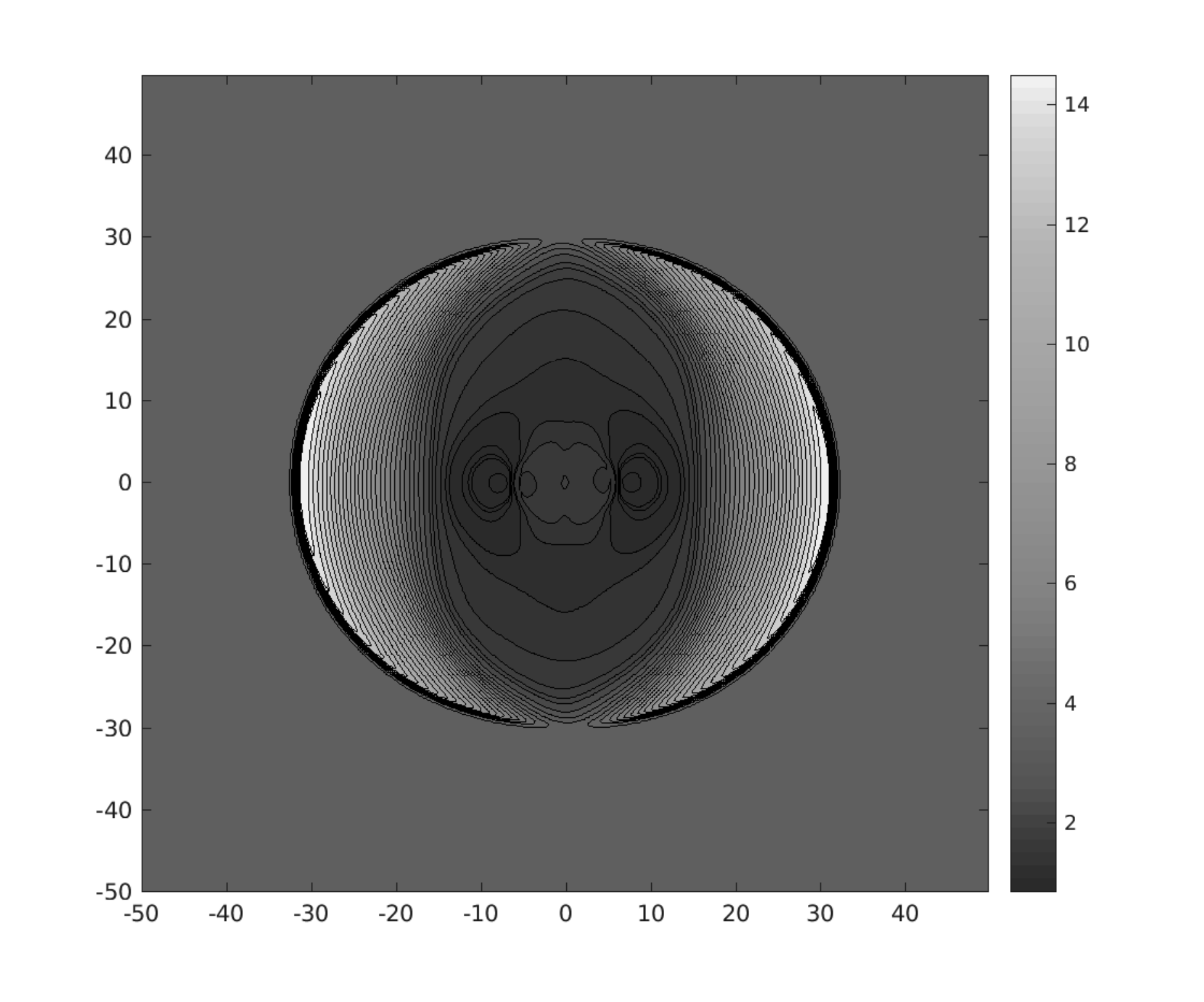}
\includegraphics[height=0.24\textheight]{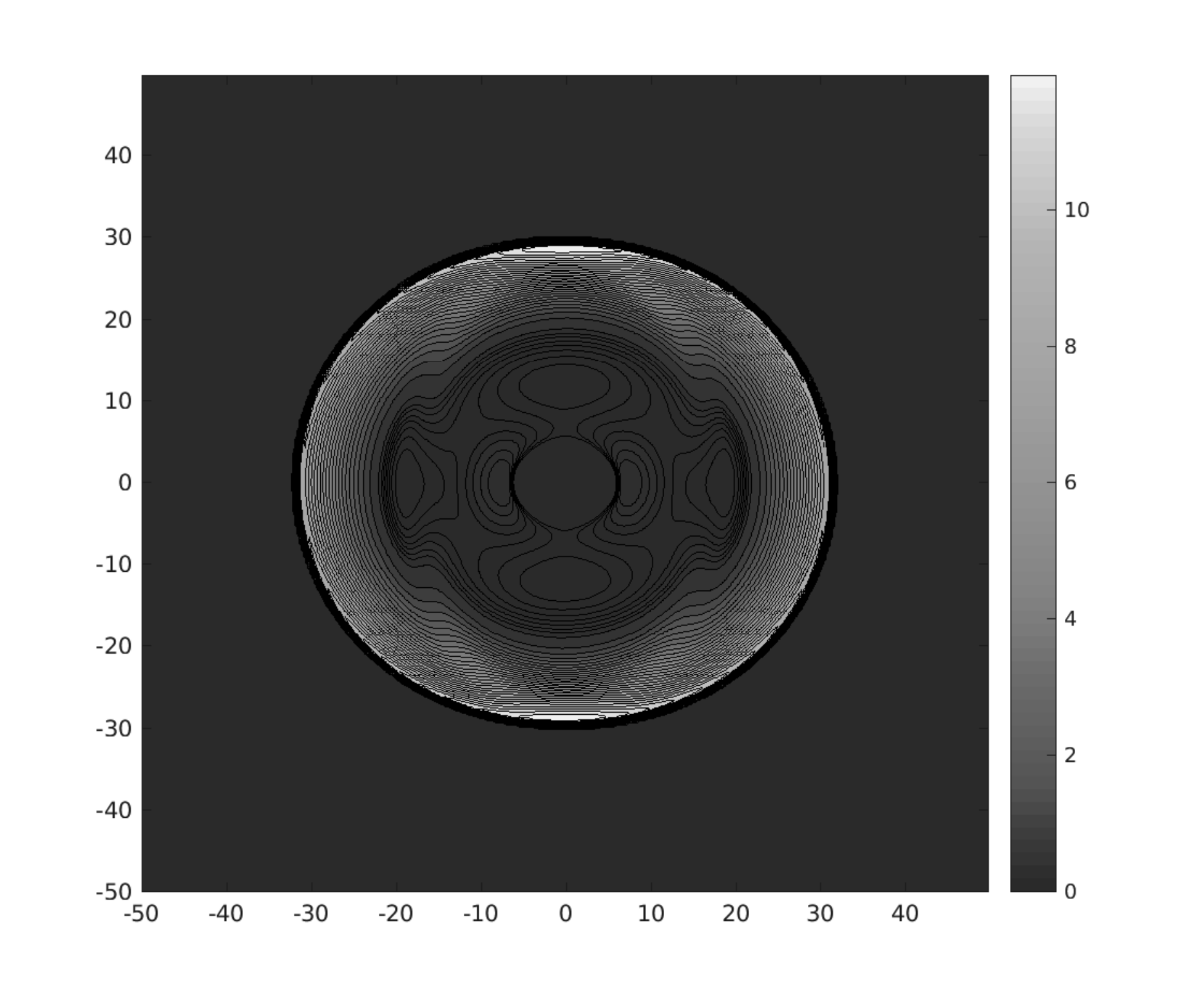}
\caption{Test \ref{test:explosion}. Solution obtained with initial magnetic field $(B_x, B_y, B_z)=(0, 5/\sqrt{\pi}, 0)$. Density $\rho$ (top left), pressure $P$ (top right), magnetic pressure (bottom left) and kinetic energy (bottom right) computed at time $t=3$.}
\label{fig:explosion2}
\end{center}
\end{figure}

If the initial magnetic field is taken as $(B_x, B_y, B_z)=(0, 5/\sqrt{\pi}, 0)$, the increasing of the strength of the magnetic field produces a slightly elongated shock in the $y$-direction: see Figure \ref{fig:explosion2}.

\begin{figure}[!ht]
\begin{center}
\includegraphics[height=0.24\textheight]{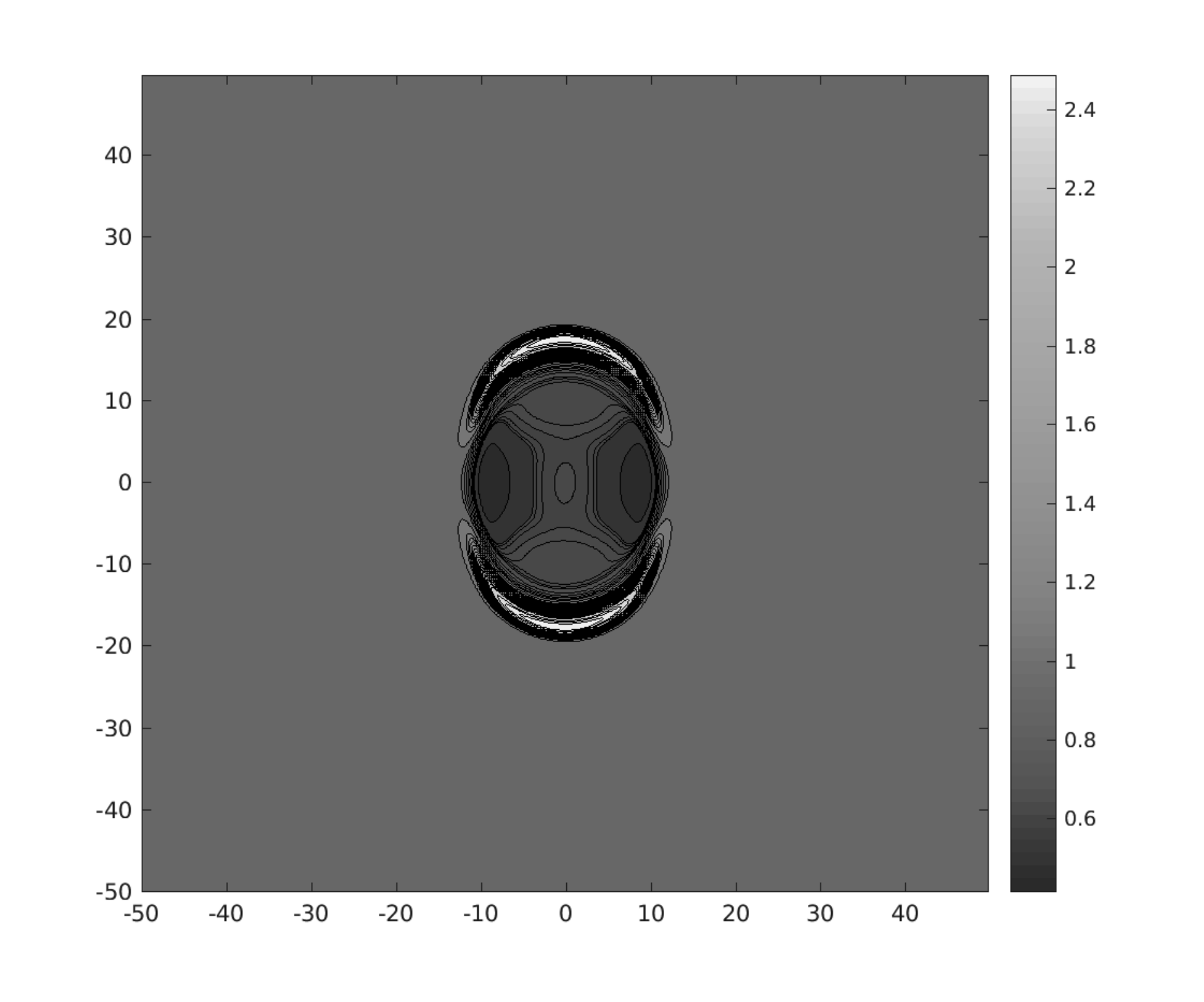}
\includegraphics[height=0.24\textheight]{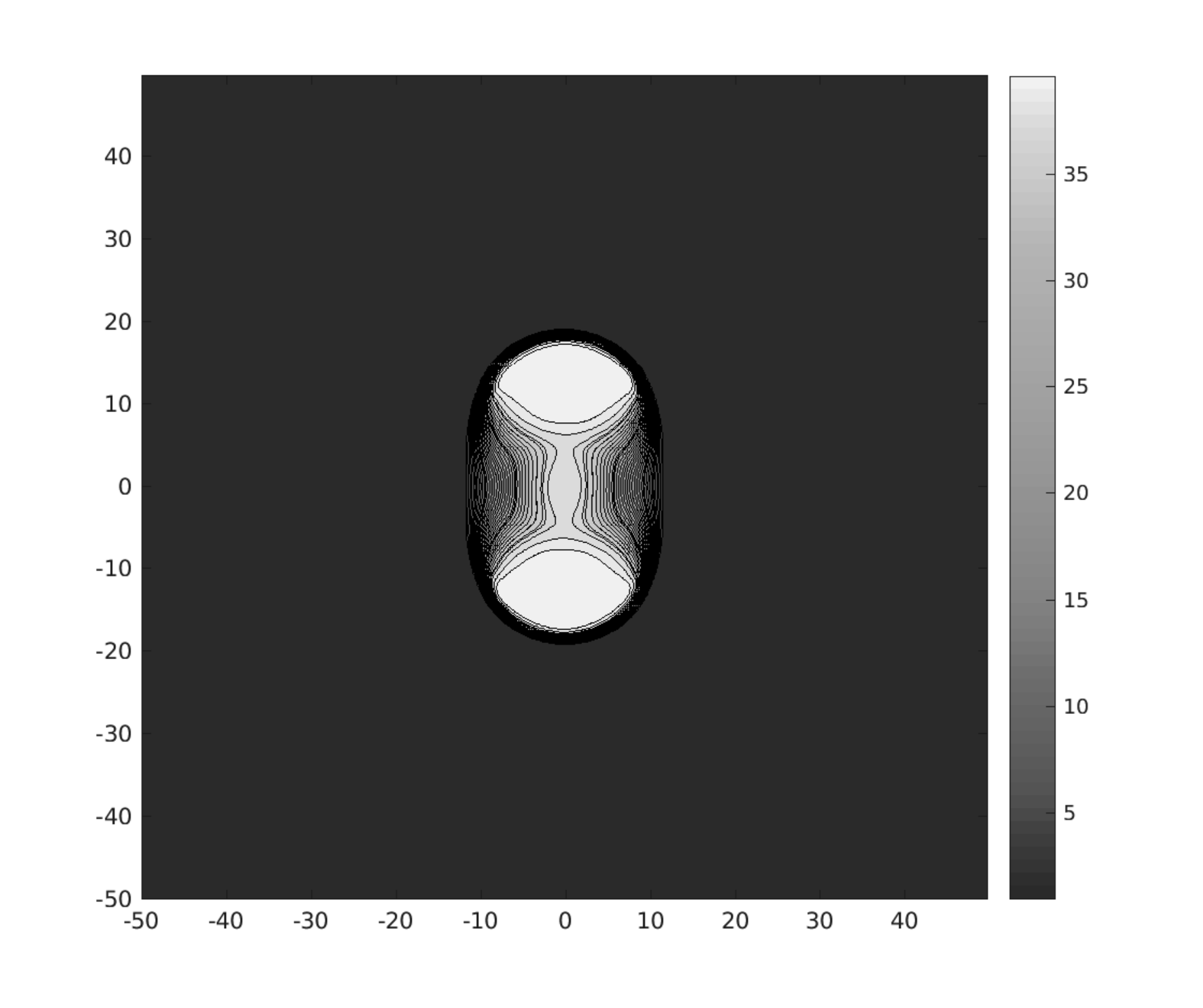} \\
\includegraphics[height=0.24\textheight]{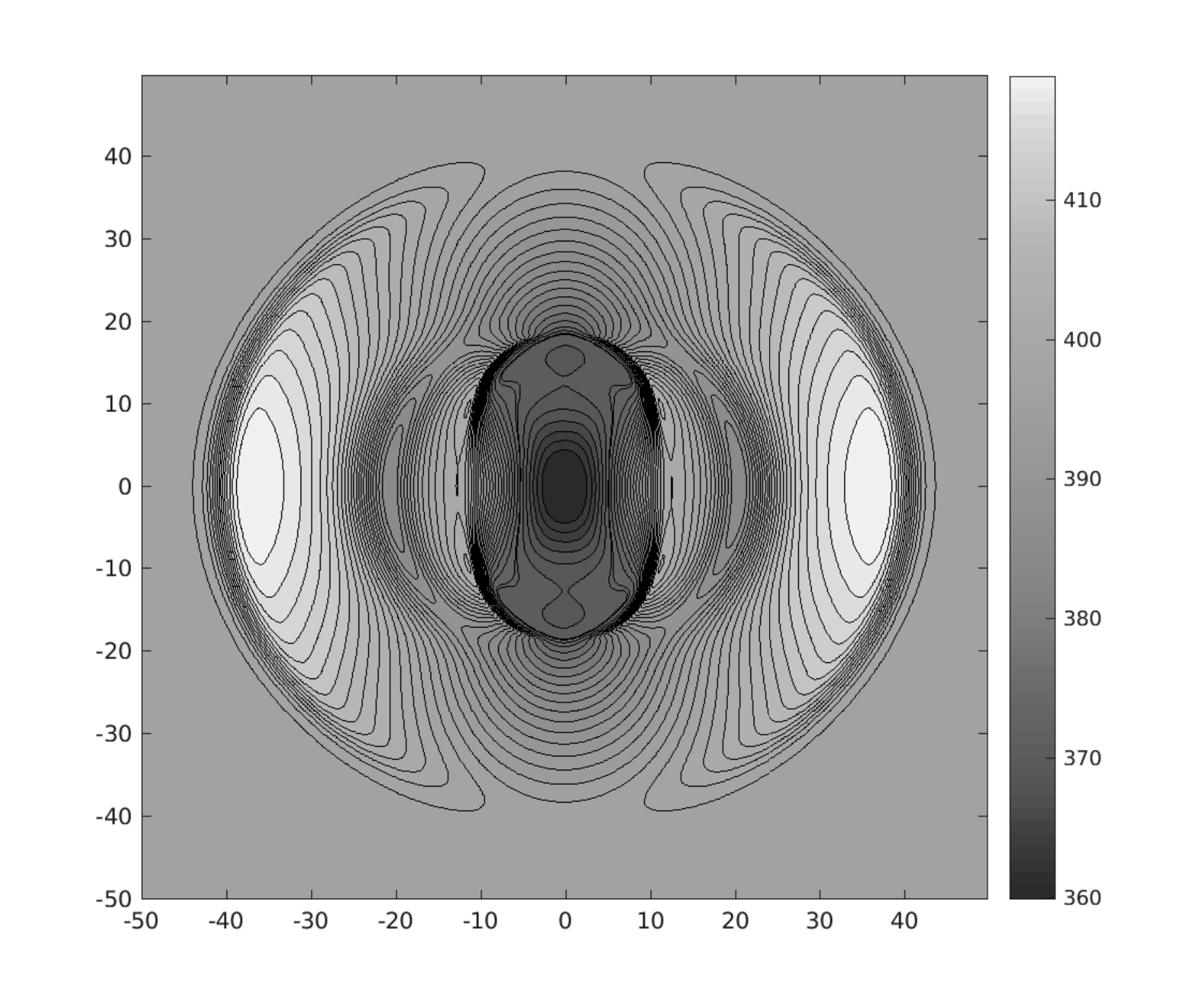}
\includegraphics[height=0.24\textheight]{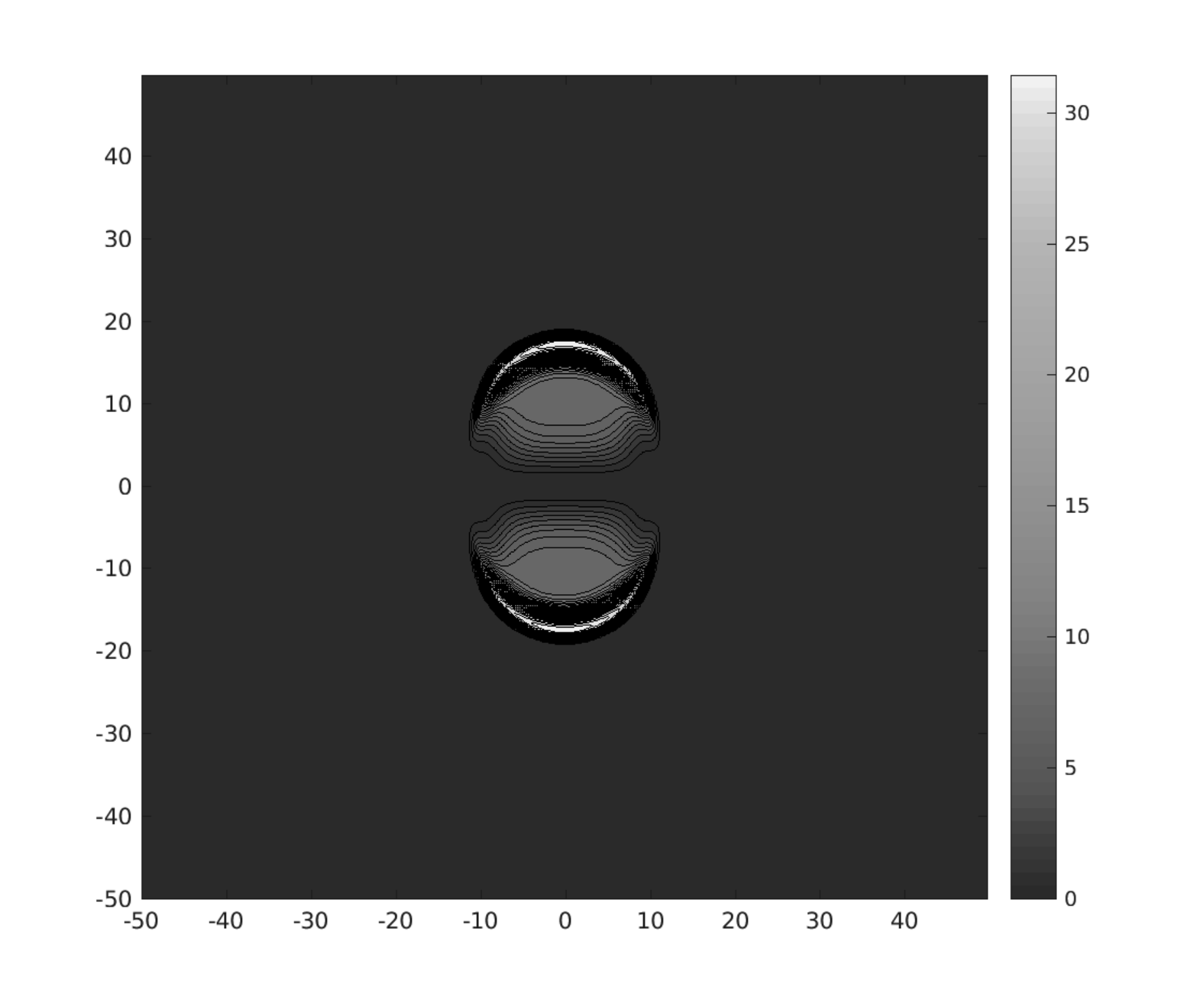}
\caption{Test \ref{test:explosion}. Solution obtained with initial magnetic field $(B_x, B_y, B_z)=(0, 50/\sqrt{\pi}, 0)$. Density $\rho$ (top left), pressure $P$ (top right), magnetic pressure (bottom left) and kinetic energy (bottom right) computed at time $t=1.05$.}
\label{fig:explosion3}
\end{center}
\end{figure}

For an even stronger initial magnetic field $(B_x, B_y, B_z)=(0, 50/\sqrt{\pi}, 0)$, the solution becomes highly anisotropic, with essentially no displacement of fluid in the orthogonal direction to the magnetic field. Figure \ref{fig:explosion3} shows the results obtained at time $t=1.05$, which is the time considered in \cite{Tang-Xu} for which the perturbation has not yet reached the boundary.

In all the three cases so far considered, our results are in good agreement with those presented in \cite{Tang-Xu,Zachary}.

%%%%

\section{Conclusions} \label{sec:conclusions}

We have introduced a general class of genuinely two-dimensional incomplete Riemann solvers for systems of conservation laws. Departing from a reinterpretation of Balsara's two-dimensional HLL scheme, we have extended it by substituting the underlying one-dimensional HLL solver with an arbitrary AVM-type numerical flux. This allows, in particular, to control the numerical diffusion of the scheme by using suitable aproximations to the absolute value function used as basis of the AVM flux. The numerical flux at an edge is built as a linear convex combination of a one-dimensional AVM flux and two multidimensional AVM corrections at the corners which take into account transversal features of the flow. In principle, the two-dimensional AVM solvers could be extended to three-dimensional problems, although this point is not addressed in this paper.

The proposed multidimensional AVM methods can be applied to arbitrary hyperbolic systems. In this paper we have focused on applications to the MHD equations, in order to make comparisons with other existing methods in the literature. Additionally, we have proposed a new technique for imposing the divergence-free condition on the magnetic field, which gives comparable results as the standard projection method. The performances of the proposed schemes have been assessed in a number of challenging problems in MHD.

Future work in this line of research include the extension of the multidimensional AVM solvers to the case of nonconservative hyperbolic systems, the construction of higher order versions of the schemes, as well as applications to different hyperbolic systems.

%%%%%

\end{document}